\documentclass[12 pt]{article}
\usepackage{amsfonts,amssymb,amscd}
\newtheorem{teo}{Theorem}[section]

\newtheorem{prop}[teo]{Proposition}

\newtheorem{lemma}[teo]{Lemma}

\newtheorem{obs}[teo]{Remark}

\newtheorem{defnc}[teo]{Definition}

\newtheorem{coro}[teo]{Corollary}

\newcommand{\qed}{\hfill \fbox{}}

\newcommand{\C}{{\mathbb C}}

\newcommand{\N}{{\mathbb N}}

\newcommand{\R}{{\mathbb R}}

\newcommand{\Z}{{\mathbb Z}}

\newcommand{\dif}{{\rm Diff}^{\omega} (S^1)}

\newcommand{\difro}{{\rm Diff}^{\omega} (\R,0)}

\newcommand{\difroCC}{{\rm Diff}\, (\C,0)}

\newcommand{\calb}{{\mathcal B}}

\newcommand{\Leb}{{\rm Leb}}

\newcommand{\fol}{{\mathcal{F}}}

\newcommand{\Length}{{\mathcal{L}\,}}

\begin{document}

\title{Global rigidity of conjugations for locally non-discrete subgroups of $\dif$}
\author{Anas Eskif \, \, \,  \& \, \, \, Julio C. Rebelo}
\date{}
\maketitle

\begin{abstract}
We prove a global topological rigidity theorem for locally $C^2$-non-discrete subgroups of $\dif$.
\end{abstract}

\noindent \hspace{1.2cm} {\small AMS-Classification (2010): 37C85, \, 22F05, \, 37E10}

\bigskip

\noindent \hspace{1.2cm} {\small Keywords: locally non-discrete group, topological rigidity, stationary measures}

\section{Introduction}

In this paper we establish a topological rigidity theorem for a large class of subgroups of the group $\dif$
consisting of (orientation-preserving) real analytic diffeomorphisms of the circle $S^1$. Indeed,
the primary object studied in this paper are finitely generated, locally $C^2$-non-discrete subgroups of $\dif$.
As is often the case, our choice of restricting attention to finitely generated
groups of orientation-preserving diffeomorphisms is made only to help us to focus on the main difficulties of the problem;
straightforward generalizations are left to the reader. On the other hand, the
regularity assumption ($C^{\omega}$) required from our diffeomorphisms is a far more important point, even if it can
substantially be weakened in several specific contexts. Possible extensions of our results to, say, smooth diffeomorphisms,
will briefly be discussed in the appendix of this article.

A group $G \subset \dif$ is said to be {\it locally $C^2$-non-discrete}\, if there is an open, non-empty interval
$I \subset S^1$ and a sequence $g_j$, $g_j \neq {\rm id}$ for every~$j \in \N$, of elements in $G$ whose restrictions to $I$
converge in the $C^2$-topology to the identity on $I$; see Section~2 for detail. Concerning these groups, for the time being,
it suffices to know that they form a rather large class of finitely generated subgroups of $\dif$. After stating
the main results of this paper, we shall resume the discussion of these groups and provide
further information on their nature.

Recall that two subgroups $G_1$ and $G_2$ of $\dif$ are said to be {\it topologically conjugate}\, if there is a
homeomorphism $h : S^1 \rightarrow S^1$ such that $G_2 = h^{-1} \circ G_1 \circ h$, i.e. to every element $g_{(1)} \in G_1$
there corresponds a unique element $g_{(2)} \in G_2$ such that $g_{(2)} = h^{-1} \circ g_{(1)} \circ h$ and conversely. Now
we have:

\vspace{0.1cm}

\noindent {\bf Theorem A}. {\sl Consider two finitely generated non-abelian topologically conjugate subgroups $G_1$ and
$G_2$ of $\dif$. Suppose that these groups are locally $C^2$-non-discrete.
Then every homeomorphism $h : S^1 \rightarrow S^1$ satisfying $G_2 = h^{-1} \circ G_1 \circ h$
coincides with an element of $\dif$.}

\vspace{0.1cm}

Theorem~A answers one of the questions raised in \cite{RebMMJ}. When this theorem is
combined with Theorem~\ref{PoissonBoundaryApplication}, we also obtain:

\vspace{0.1cm}

\noindent {\bf Theorem B}. {\sl Suppose that $\Gamma$ is a finitely generated hyperbolic group which is neither finite nor
a finite extension of $\Z$ and consider two topologically
conjugate faithful representations $\rho_1 : \Gamma \rightarrow \dif$ and $\rho_2 : \Gamma \rightarrow \dif$
of $\Gamma$ in $\dif$. Assume that $G_1 = \rho_1 (\Gamma) \subset \dif$ is locally $C^2$-non-discrete.
Then every (orientation-preserving) homeomorphism $h : S^1 \rightarrow S^1$ conjugating the representations
$\rho_1$ and $\rho_2$ coincides with an element of $\dif$.}

\vspace{0.1cm}

The main assumptions of Theorems~A and~B, namely the fact that our groups are locally $C^2$-non-discrete,
cannot be dropped. Indeed, counterexamples for the previous statements in the context of discrete groups
can be obtained in a variety of ways. For example, two cocompact representations in ${\rm PSL}\, (2, \R)$ of the fundamental
group of the genus~$g$ compact surface ($g\geq 2$) are always topologically conjugate. However these representations
are not $C^1$-conjugate unless they define the same point in the Teichmuller space.
Similarly a standard Schottky (free) group on two generators acting on $S^1$ gives rise to an action that is {\it structural
stable}\, in $\dif$. Thus, by perturbing the generators inside $\dif$, we obtain numerous actions that are
topologically but not $C^1$ conjugate to the initial Schottky group (cf. \cite{sullivan1} and references therein).

Theorem~B also becomes false if the group $\Gamma$ is $\Z$. The classical example due
to Arnold \cite{arnold} of an analytic diffeomorphism of $S^1$ topologically conjugate to an irrational rotation
by a singular homeomorphism provides an interesting example since the cyclic group generated by an irrational
rotation clearly verifies the condition of being (locally) $C^2$-non-discrete. Concerning the possibility of
generalizing Theorems~A and~B to higher rank abelian groups,
the reader is referred to the discussions in \cite{moser} and \cite{yoccoz}. On the other hand,
by virtue of the work of Kaimanovich and his collaborators, Theorem~B still
holds true for other type of groups including relatively hyperbolic ones; cf.
\cite{connel} and its references.

The above theorems also have consequences of considerable interest in the theory of secondary characteristic classes
of (real analytic) foliated $S^1$-bundles. For example, we state:

\vspace{0.1cm}

\noindent {\bf Corollary~C}. {\sl Let $(M_1, \fol_1)$ and $(M_2, \fol_2)$ be two analytic foliated
$S^1$-bundles over a same hyperbolic manifold~$N$. Assume
that these foliated $S^1$-bundles are topologically conjugate and that the holonomy
group of $(M_1, \fol_1)$ is locally $C^2$-non-discrete. Assume moreover that $\fol_1$ has at least one
simply connected leaf. Then the Godbillon-Vey classes of $(M_1, \fol_1)$ and $(M_2, \fol_2)$
coincide.}

\vspace{0.1cm}

Concerning Corollary~C, it is well known that Godbillon-Vey classes are invariant by homeomorphisms that
are transversely of class $C^2$; see \cite{candel}. Also, in the case in question,
the holonomy group $G_1 \subset \dif$ associated with the foliation
$\fol_1$ is, by assumption, locally $C^2$-non-discrete. Furthermore, since $N$ is hyperbolic and $\fol_1$ has at
least one leaf simply connected, there follows that $G_1$ is of hyperbolic type (isomorphic as abstract group to
the fundamental group of $N$). Thus Theorem~B ensures the transverse regularity of the conjugating homeomorphism.

The rest of this introduction contains an overview of our approach to the proofs of Theorems~A and~B including the
main connections with previous works and some interesting examples.

Very roughly speaking, the results in this paper are obtained by blending the technique of ``vector fields in the
closure of groups'', developed in \cite{russo} and \cite{na} for subgroups of ${\rm Diff}\, (\C,0)$ and in \cite{reb1}
for subgroups of $\dif$, with results related to stationary measures on $S^1$, see \cite{dkn}, \cite{antonov},
\cite{kleptsynnalski} and with measure-theoretic boundary theory for groups \cite{deroinETDS}, \cite{KaimanovichPoisson},
and \cite{connel}. We shall follow a chronological order to
explain the various connections between these works.

First, Shcherbakov and Nakai \cite{russo}, \cite{na}
have independently studied the dynamics of non-solvable subgroups of ${\rm Diff}\, (\C,0)$
and they observed the existence of certain vector fields whose local flows were ``limits'' of actual elements in
the pseudogroup (see Section~2 for detail). Slightly later, Ghys \cite{ghys} noted that non-solvable subgroups
of ${\rm Diff}\, (\C,0)$ always contain a non-trivial sequence of elements converging to the identity.
By analogy with the case of finite dimensional Lie groups, he suggested that the existence of vector fields with
similar properties should be a far
more general phenomenon and he went on to discuss the topological dynamics of the corresponding groups on the circle.

In the case of the circle, the program proposed by Ghys was fairly accomplished in \cite{reb1}. In this paper, vector fields
whose local flow is a limit of actual elements in a given group are said to belong to the {\it closure of the group}\, (see Section~2
for proper definitions). The role of ``locally non-discrete subgroups of $\dif$'' was emphasized and it was shown that
these locally non-discrete subgroups of $\dif$ admit non-zero vector fields in their closure. As an application
of these vector fields, the following theorem was also proved in \cite{reb1}:

\vspace{0.1cm}

\noindent {\bf Theorem} ({\bf [R1]}). {\sl There exists a neighborhood $\mathcal{U}$ of the identity in
$\dif$ with the following property. Assume that $G_1$ (resp. $G_2$) is a subgroup of $\dif$ generated
by diffeomorphisms $g_{1,1}, \ldots , g_{1,N}$ (resp. $g_{2,1}, \ldots , g_{2,N}$) lying in $\mathcal{U}$.
If $h :S^1 \rightarrow S^1$ is a homeomorphism satisfying $g_{2,i} = h^{-1} \circ g_{1,i} \circ h$ for every
$i=1, \ldots ,N$, then $h$ coincides with an element of $\dif$.}

\vspace{0.1cm}

\noindent This theorem can be thought of as a local version of Theorem~A.
In fact, the assumption that $h$ takes a generating set formed by elements ``close to the identity'' to elements
that are still close to the identity gives the statement in question an intrinsic {\it local character}.
For example, the above theorem from \cite{reb1} is satisfactory for
deformations/pertubations problems but falls short of answering the same question for general groups admitting generating sets
in the fixed neighborhood $\mathcal{U}$ unless the mentioned sets are, in addition, conjugated by $h$.
This type of difficulty was pointed out
and discussed in \cite{RebMMJ} and the method of \cite{reb1} suggests that these rigidity phenomena should hold
for general {\it locally non-discrete subgroups of $\dif$}\, (again see Section~2 for accurate definitions). From here one sees
that the main motivation of the present work was to clarify these issues.

It is mentioned in \cite{reb1} that the main example of locally non-discrete subgroups of $\dif$ is provided by
non-solvable groups admitting a finite generating set contained in $\mathcal{U} \subset \dif$, as follows
from Ghys's results in \cite{ghys}. Conversely the main examples of groups that are {\it locally discrete}\,
are provided by Fuchsian groups. The problem about understanding how the subgroups of $\dif$ are split in locally
discrete and locally non-discrete ones is then unavoidably raised.

Soon it became clear that locally discrete groups were, indeed, very common (see for example \cite{TAMS}). The
problem of finding {\it locally discrete}\, subgroups of $\dif$ beyond the context of Fuchsian groups, however, proved to be
much harder. Recently, however, much progress has been made towards the understanding of the structure of
locally discrete groups
thanks to the works of Deroin, Kleptsyn, Navas, and their collaborators, see \cite{advances} and the survey \cite{dknTokyo}
for some up-to-date information.
Meanwhile it was also observed in \cite{rebproceedings} that the Thompson-Ghys-Sergiescu
subgroup of ${\rm Diff}^{\infty} (S^1)$ is locally discrete. Whereas this example is only smooth, as opposed to real
analytic, the observation in question connects with the fundamental notion of {\it non-expandable point}\, and this
requires a more detailed explanation.

Fix a group $G$ of diffeomorphisms of $S^1$. A point $p \in S^1$
is said to be {\it expandable}\, (for the group $G$) if there is an element $g \in G$ such that $\vert g'(p) \vert >1$.
Among ``large'' (e.g. non-solvable) subgroups of $\dif$ all of whose orbits are dense, ${\rm PSL}\, (2, \Z)$
constitutes the simplest example of group with a {\it non-expandable}\, point. In turn, it is
observed in \cite{rebproceedings} that a {\it locally non-discrete} group all of whose orbits are dense cannot have
a non-expandable point (see also Lemma~\ref{FinishingJob-11.1}) and this implies
the above mentioned conclusion concerning the Thompson-Ghys-Sergiescu group.
Hence, a method to produce locally discrete groups consists of finding groups with non-expandable points and this does
not depend on whether or not the group arises from a Fuchsian group. Very recently, V. Kleptsyn and his collaborators
have found {\it free subgroups}\, of $\dif$ with non-expanding points and which are not conjugate to Fuchsian groups.
This very interesting result also justifies our Theorem~B: note that the groups in question are free and hence
they are also hyperbolic so that Theorem~B applies to them. In other words, they are never topologically
conjugate to a locally $C^2$-non-discrete subgroup of $\dif$.

Nonetheless, the full understanding of locally discrete subgroups of $\dif$ was not yet reached (see \cite{advances} and
\cite{dknTokyo} for further information). To continue
our discussion, we shall then restrict ourselves to the related problem of understanding ``rigidity'' of topological
conjugations between subgroups of $\dif$ which, ultimately, constitutes the actual purpose of this paper.
In the sequel, we then consider
two topologically conjugate subgroups $G_1$ and $G_2$ of $\dif$. Since topological rigidity is targeted, the examples
provided in the beginning of the introduction indicate that one of the groups, say $G_1$, should be assumed to be
$C^2$-locally non-discrete. At this level, Theorem~A fully answers the question provided that $G_2$ is
locally $C^2$-non-discrete as well. Thus, to make further progress, we need to investigate whether a
locally $C^2$-non-discrete group $G_1$ can be topologically conjugate to a locally $C^2$-discrete subgroup $G_2$.
Following our above stated results, the state-of-art of this problem can be summarized as follows.

First we assume once and for all that $G_1$ (and hence $G_2$) is {\it minimal}\, i.e. all of its orbits are dense
in $S^1$. Moreover these groups are also assumed to be non-abelian. The material
presented in Sections~2 and~3 of the present paper shows that this assumption can be made
without loss of generality (otherwise topological rigidity is always verified). Theorem~B also settles the question
when the groups are of hyperbolic type. Also, if $G_2$ is conjugate to a Fuchsian group, then a conjugating homeomorphism
$h$ between $G_1$ and $G_2$ cannot exist as pointed out in \cite{RebMMJ}. These general statements apart, the existence
of non-expanding points plays again a role in the problem. Thus we may consider the obvious alternative
\begin{itemize}
  \item All points in $S^1$ are expandable for $G_2$.

  \item $G_2$ has at least one non-expandable point.
\end{itemize}
In the first case, an unpublished result of Deroin asserts that the (locally $C^2$-discrete group)
$G_2$ is essentially a Fuchsian group. Therefore the preceding implies that a conjugating homeomorphism between $G_1$
and $G_2$ cannot exist (cf. \cite{RebMMJ}).
Alternatively, the non-existence of topological conjugation between $G_1$ and $G_2$ can directly be derived from
Theorem~\ref{ConjugatignEpxanding} in Section~5. In fact, the argument in Section~5
relies only on the following assumptions:
\begin{enumerate}
  \item $G_1$ is locally $C^2$-non-discrete.

  \item $G_1$ is minimal and non-abelian.

  \item Every point in $S^1$ is expandable for $G_2$.
\end{enumerate}
In this respect, the examples obtained by Kleptsyn et al. of free subgroups of $\dif$ require a few additional comments.
Namely, these groups contain diffeomorphisms exhibiting three fixed points (two of them hyperbolic and the remaining one
parabolic) so that the group is not topologically conjugate to a Fuchsian group. Moreover,
these groups possess non-expandable points (so that they are locally $C^2$-discrete). Our results also show, or recover the
information that might be read off their construction, that these examples by Kleptsyn et al.
form an independently class of groups as they are topologically conjugate neither to
Fuchsian groups nor to locally $C^2$-non-discrete groups (Theorem~B).

In closing, recall that a classical problem that lends further interest to regularity properties of homeomorphisms
conjugating groups actions is the possibility of having different Godbillon-Vey characteristic classes. In the case
of (global) groups acting on $S^1$, our results are satisfactory provided that one of the groups is
locally $C^2$-non-discrete. On the other hand, in the locally discrete case, this problem is difficult even if the
groups in question arise from Fuchsian groups and we refer the reader to \cite{ghys2} and its references for further
information.

To finish the introduction, let us provide an overview of the structure of this paper. Section~2 contains accurate definitions
for most of the notions relevant for this paper. The section then continues by revisiting results related to Shcherbakov-Nakai
theory in a form adapted to our needs. The second half of this section, namely Subsection~2.2, provides
a description of the topological dynamics associated with a locally $C^2$-non-discrete subgroup of $\dif$. This
description faithfully parallels the corresponding results established in \cite{ghys}
for the case of groups admitting a generating set ``close to the identity''.

Section~3 is devoted to proving Theorem~A in different types of special situations.
These include the case where the groups  $G_1$
and $G_2$ have finite orbits as well as the case in which these groups are solvable but non-abelian. The results of Section~3 are
implicitly used throughout the paper since, from Section~4 on, they allow us
to restrict our discussion to a sort of ``generic case'' for the group $G_1$. Roughly speaking, this generic situation
is such that we can fix and interval $I \subset S^1$ and, for every $\varepsilon >0$, we can find
a finite collection of elements in $G_1$ satisfying the following conditions:
\begin{itemize}
  \item Diffeomorphisms in this collection are $\varepsilon$-close to the identity in the $C^2$-topology on $I$.
  \item The collection of these diffeomorphisms generated a non-solvable subgroup of $\dif$.
\end{itemize}

In Section~4, we construct an explicit sequence of diffeomorphisms in $G_1$ converging to the identity in the
$C^2$-topology on the above mentioned interval~$I$. As explained in the beginning of Section~4, this construction
is necessary to control the convergence rate of this sequence to the identity while also having a bound on
the growing rate of the sequence formed by the corresponding higher order derivatives.
Note that the existence of a sequence converging to the identity in the $C^2$-topology on some interval
is the very definition of a locally $C^2$-non-discrete group, yet this definition does not provide any estimate on
the $C^3$-norm (for example) of the diffeomorphisms in this sequence (see Section~4 for a detailed discussion).
In the construction of the mentioned desired sequence, we take
advantage of the fact that we can select finitely many elements of $G_1$ generating a non-solvable
group and being arbitrarily close to the identity on a fixed interval~$I$.

In Section~5, we shall prove Theorem~A modulo Proposition~\ref{FinishingJob-11.3} whose proof is deferred
to Section~6. In fact, in this section
Theorem~\ref{ConjugatignEpxanding} will be proved and this theorem provides a statement fairly stronger than
than what is strictly needed to derive Theorem~A in the ``generic case'' (under discussion since Section~4). As
to Proposition~\ref{FinishingJob-11.3}, the reader will note that its use can be avoided modulo working with bounded distortion
estimates for iterates of diffeomorphisms possessing parabolic fixed points, as similarly done in \cite{reb1}. The
interest of Proposition~\ref{FinishingJob-11.3} lies primarily in the fact that it saves significant space by
allowing us to focus exclusively on hyperbolic fixed points which, in turn, are linearizable \cite{Stenrberg}.

In Section~6 we collect essentially all the results in this paper for which Ergodic theory appears to be an indispensable
tool. The two main statements in this section are Proposition~\ref{FinishingJob-11.3} and Theorem~\ref{PoissonBoundaryApplication}.
Since the role played by Proposition~\ref{FinishingJob-11.3} in our paper was already discussed, for the time
being it suffices to mention that Theorem~\ref{PoissonBoundaryApplication} reduces Theorem~B to Theorem~A. We also note
that the proof of Proposition~\ref{FinishingJob-11.3} relies heavily on \cite{dkn} and, in fact, this proposition
is a straightforward consequence of the proof of ``Th\'eor\`eme F'' in \cite{dkn}.
The proof of Theorem~\ref{PoissonBoundaryApplication} is more involved as it combines standard facts about hyperbolic groups
with results from \cite{deroinETDS} and from \cite{KaimanovichPoisson} and still depends in a crucial way of
Proposition~\ref{regularstationary-addedinApril} which deserves additional comments.

Roughly speaking, given a ``generic'' (finitely generated) locally non-discrete group $G \subset \dif$,
Proposition~\ref{regularstationary-addedinApril} ensures the existence of a measure $\mu$ on $G$ giving rise to
an absolutely continuous stationary measure. Whereas this statement has interest in its own right and deserves
to be called theorem, as opposed to proposition, its status is not clear to us. Indeed,
this results seems to be known to some experts as as a consequence of the techniques introduced in \cite{connel}.
In view of this, we have decided not to include a fully detailed proof of it which, in addition,
would have made this paper a few pages longer. However, we present a clear sketch
of the corresponding proof, referring to \cite{connel} for estimates and convergence details only.
It may also be interesting to note that, albeit
our discussion relies strongly on \cite{connel}, it is conceivable that Furstenberg's methods \cite{furst}
might as well be enough to imply the desired statement.

Finally the appendix (Section~7) contains a partial answer in the analytic category
to a question raised in \cite{deroinETDS}. The argument exploits the construction carried out in Section~4.
The appendix then ends with a summary of the role played by the regularity assumption ($C^{\omega}$) in this
paper. In particular, we highlight some specific problems whose solutions would lead to non-trivial generalizations of our statements
to less regular groups of diffeomorphisms.

\noindent {\bf Acknowledgements}. The authors are grateful to V.Kleptsyn for interesting discussions and for
suggesting us the statement of Proposition~\ref{hisHoldercontinuous}.
A good part of this work was carried out during a long term
visit of the second author to IMPA and he wishes to warmly thank P.Sad for his invitation.

\section{Locally non-discrete groups: vector fields and topological dynamics}
The definition of locally non-discrete groups is implicit in \cite{reb1} and formulated in \cite{TAMS} and \cite{deroinETDS}.
In the analytic case, it reads as follows:

\begin{defnc}
\label{localnondiscreteness}
A subgroup $G$ of $\dif$ is said to be locally $C^m$-non-discrete if there is a non-empty open interval $I \subseteq S^1$
and a sequence of elements $\{ g_i \} \subset G$, $g_i \neq {\rm id}$ for every~$i \in \N$, whose restrictions $\{ g_{i \vert I} \}$
to $I$ converge to the identity in the $C^m$-topology (as maps from $I$ to $S^1$).
\end{defnc}

Naturally, a group $G \subset \dif$ is called {\it locally $C^m$-discrete}\,
if it fails to satisfy the conditions of Definition~\ref{localnondiscreteness}.
Unless otherwise stated, the terminology used in this paper is such that every {\it interval}\, is open,
connected and non-empty. In the course of this paper, we shall mainly work with locally $C^2$-non-discrete subgroups
of $\dif$.

It also useful to adapt Definition~\ref{localnondiscreteness} to the context of pseudogroups.

\begin{defnc}
\label{localnondiscretenesspseudogroup}
Consider an open set $U \subset \R$ along with a pseudogroup $\Gamma$ of analytic diffeomorphisms from open subsets of $U$ to $\R$.
The pseudogroup $\Gamma$ is said to be locally $C^m$-non-discrete if
there is an interval (open, connected and non-empty) $I \subset U$ and a sequence of maps $\{ g_i \} \subset \Gamma$
satisfying the following conditions:
\begin{enumerate}
  \item For every $i \in \N$, the interval $I$ is contained in the domain of definition of $g_i$ viewed as element of the
  pseudogroup $\Gamma$.
  \item The restriction $g_{i \vert I}$ of $g_i$ to $I$ does not coincide with the identity map.
  \item The sequence $\{ g_{i \vert I} \}$ formed by the restrictions of the $g_i$ to $I$ converge to the identity
  in the $C^m$-topology (as maps from $I$ to $\R$).
\end{enumerate}
\end{defnc}

\subsection{Vector fields in the closure of pseudogroups}

Vector fields whose local flow can be approximated by elements in the initial group (pseudogroup)
constitute a very important tool to investigate the dynamics associated with {\it locally non-discrete}\,
groups (pseudogroups). The idea of approximating a flow by elements in a group/pseudogroup
is made accurate by the following definition.

\begin{defnc}
\label{flowclosurepseudogroup}
Consider an open set $U \subset \R$ along with a pseudogroup $\Gamma$ of maps from open subsets of $U$ to $\R$.
Consider also a vector field $X$ defined on an
interval $I \subset U$ and let $\Psi_X$ denote its local flow.
The vector field $X$ is said to be (contained) in the $C^m$-closure of $\Gamma$ if, for every interval $I_0 \subset I$
and for every $t_0 \in \R_+$ such that $\Psi_X^t$ is defined on $I$ for $0 \leq t  \leq t_0$,
there exists a sequence of maps $\{ g_i \} \subset \Gamma$ satisfying the conditions below:
\begin{itemize}
  \item For every $i \in \N$, the interval $I_0$ is contained in the domain of definition of $g_i$ viewed as element of the
  pseudogroup $\Gamma$.

  \item The sequence $\{ g_{i \vert I_0} \}$ formed by the restrictions of the $g_i$ to $I_0$ converge to
  $\Psi_X^{t_0} : I_0 \rightarrow \R$ in the $C^m$-topology (where $k \in \N \cup \{ \infty \}$).
\end{itemize}
\end{defnc}

Unless otherwise mentioned, whenever we mention a vector field $X$ belonging to the closure of a pseudogroup $\Gamma$ it is
implicitly assumed that this vector field does not vanish identically.
It is clear from the definitions that a pseudogroup containing
some (non-identically zero) vector field in its $C^r$-closure cannot be {\it locally $C^r$-discrete}.

Before going further into the structure of the topological dynamics of locally $C^2$-non-discrete subgroups
of $\dif$, let us quickly revisit some results established by Shcherbakov and Nakai for pseudogroups of holomorphic
diffeomorphisms fixing $0 \in \C$; see \cite{na}, \cite{russo}. The discussion below is slightly simplified
by the fact that only local diffeomorphisms having real coefficients will be considered. Let $\difro$ denote the group
of germs of orientation-preserving analytic diffeomorphisms, i.e. if $g \in \difro$ then we assume that $g' (0) >0$.

First, we have:

\begin{lemma}
\label{preliminarieslemma1}
Let $\Gamma$ be a pseudogroup generated by finitely many elements of $\difro$ and denote by
$\Gamma_{(0,0)}$ the group of germs at $(0,0)$ corresponding to $\Gamma$. Then the following are equivalent:
\begin{enumerate}
  \item $\Gamma_{(0,0)}$ is an infinite cyclic group unless it is reduced to the identity.
  \item $\Gamma$ is locally $C^r$-discrete, for every $r \in \N \cup \{ \infty \}$.
  \item $\Gamma$ does not contain vector fields in its $C^r$-closure, for every $r \in \N \cup \{ \infty \}$.
\end{enumerate}
\end{lemma}

\noindent {\it Proof}. Since the elements in $\Gamma \subset \difro$ are assumed to preserve the orientation of
$\R$, it follows at once that every element different from the identity in $\Gamma_{(0,0)}$ has infinite order. Assuming
once and for all that $\Gamma_{(0,0)}$ is not reduced to the identity, consider an
element $g \neq {\rm id}$ in $\Gamma_{(0,0)}$. Modulo replacing $g$ by its inverse $g^{-1}$, we can assume that
$g' (0) \leq 1$. Let us then split the discussion into two cases.

\vspace{0.1cm}

\noindent {\it Case 1}. Suppose that $g' (0) = \lambda <1$. In this case, there are local (analytic) coordinates
where $g(x) =\lambda \, x$; see \cite{Stenrberg}. Now suppose that $\Gamma_{(0,0)}$ is an abelian group. Then, in the above coordinates,
every element of $\Gamma_{(0,0)}$ coincides with a linear map of type $x \mapsto c \, x$ for a constant $c \in \R^{\ast}_+$.
In other words, $\Gamma_{(0,0)}$ is naturally identified to a multiplicative subgroup of $\R_+^{\ast}$. The
equivalence of the three statements above then becomes clear.

Suppose now that $\Gamma_{(0,0)}$ is not abelian and consider another element $g_1 \neq {\rm id}$ belonging in $D^1 \Gamma_{(0,0)}$.
Though $g_1 \neq {\rm id}$, the derivative of $g_1$ at $0 \in \R$ equals~$1$ since $g_1 \in D^1\Gamma_{(0,0)}$ is
a product of commutators. Now,
by repeating the standard argument of Shcherbakov-Nakai with elements of
$\Gamma_{(0,0)}$ having the form $\lambda^{-N(k)} . g(\lambda^{N(k)} x)$, it is well-known that
a suitable choice of the integers $N(k)$ leads to an analytic vector field $X$ in the $C^{\infty}$-closure $\Gamma$
(the reader will note that this vector field $X$ is defined around $0 \in \R$ which in general does not happen for
Shcherbakov-Nakai vector fields). The proof of the lemma is therefore completed in Case~1.

\vspace{0.1cm}

\noindent {\it Case 2}. Suppose that $g' (0) = 1$. In view of the preceding discussion, we can assume without lost
of generality that every element in $\Gamma_{(0,0)}$ is actually tangent to the identity. Under this assumption, denote
by $Y$ the formal vector field whose induced time-one map coincides with~$g$.

Whereas the vector field $Y$ is only formal, its (formal) real flow $\Psi_Y^t$ contains the germs of all elements in $\difro$
commuting with $g$. Let $\mathcal{T}$ be the sets of those values of $t \in \R$ for which
$\Psi_Y^t$ actually defines an element of $\difro$. Clearly $\mathcal{T}$ is an additive subgroup of $\R$. Moreover,
it is well-known that the formal power series defining
$Y$ will be convergent provided that the set $\mathcal{T}$ is not discrete in $\R$, see \cite{Baker}, \cite{Ecalle}.
We will split the discussion into two cases according to whether or not $\Gamma_{(0,0)}$ is abelian.

Suppose first that $\Gamma_{(0,0)}$ is abelian so that it embeds in the $1$-parameter group generated by the formal flow
of $Y$. Clearly $\Gamma_{(0,0)}$ is infinite cyclic if and only if
$\mathcal{T}$ is a discrete subgroup of $\R$. In this case, there also follows that $\Gamma$ is locally $C^m$-discrete
and that $\Gamma$ contains only trivial vector fields in its $C^m$-closure (for every $r \in \N \cup \{ \infty \}$). Conversely,
if $\mathcal{T}$ is not discrete in $\R$, then it must be dense. Furthermore the
formal vector field $Y$ turns out to be analytic (\cite{Baker}, \cite{Ecalle}). It is now immediate
to check that $Y$ itself is contained in the $C^{\infty}$-closure of $\Gamma$.

It only remains to consider the case where $\Gamma_{(0,0)}$ is not abelian. Since $\Gamma_{(0,0)}$ is also contained in the (normal)
subgroup of $\difro$ consisting of elements tangent to the identity, being non-abelian implies the existence of elements
$g_1, \, g_2 \in \Gamma_{(0,0)}$, $g_1, \, g_2 \neq {\rm id}$, having different orders of contact with the identity. These two
elements can then be used to produce a vector field of Shcherbakov-Nakai in the $C^{\infty}$-closure of $\Gamma$.
The lemma is proved.\qed

Shcherbakov-Nakai vector fields for non-solvable subgroups of $\difroCC$
was the first genuinely non-linear situation where vector fields in the closure of (countable) groups were proven to exist.
Subgroups of $\difro$ (or even of $\difroCC$) are obviously special, as opposed to groups of $\dif$, in the sense
that their elements share a same fixed point namely, the origin. In addition to the existence of free discrete subgroups
in $\dif$, the absence of a common fixed point for elements in free
subgroups of $\dif$ is the main obstacle to extend to this context the results obtained in \cite{na}, \cite{russo}.
This difficulty was overcome for the first time in \cite{reb1}. The following lemma
singles out the key point that is common to all constructions of vector fields having similar properties
(for detailed explanations see \cite{RebMMJ}).

\begin{lemma}
\label{C1vectorfield}
Suppose that the pseudogroup $\Gamma$ consisting of local diffeomorphisms from open sets of an open (non-empty) interval $J
\subset \R$ to $\R$ contains
a sequence of elements $\{ \tilde{g}_i \}$ satisfying the following conditions:
\begin{enumerate}
  \item For every~$i$, $\tilde{g}_i$ is defined on a fixed non-empty open interval $I$.
  Moreover, again for every~$i$, the restriction $\tilde{g}_{i\vert_I}$ of $\tilde{g}_i$ to $I$ is different from the identity.

  \item The sequence of local diffeomorphisms $\tilde{g}_{i\vert_I}$ converges
  to the identity in the $C^m$-topology.

  \item There is a uniform constant $C$ such that
  $$
  \Vert \tilde{g}_i - {\rm id} \Vert_{m,I} \leq C\, \Vert \tilde{g}_i - {\rm id} \Vert_{m-1,I}
  $$
where $\Vert \, . \, \Vert_{m,I}$ (resp. $\Vert \, . \, \Vert_{m-1,I}$) stands for the $C^r$-norm (resp. $C^{m-1}$-norm) of
$\tilde{g}_i - {\rm id}$ on~$I$.
\end{enumerate}
Then there is a (non-identically zero) vector field $X$ contained in the $C^{r-1}$-closure of $\Gamma$.
\end{lemma}

\noindent {\it Proof}. For every $i$, we consider the vector field $X_i$ defined on $I$ by the
formula
$$
X_i = \frac{1}{\Vert \tilde{g}_i -{\rm id} \Vert_{m,I}} \, (\tilde{g}_i (x) -x) \partial /\partial x \, .
$$
It follows at once that the $C^m$-norm of $X_i$ on $I$ is bounded by~$1$ and, in addition, that this bound is attained
in the closure of $I$. In turn, condition~3 above shows that the $C^{m-1}$-norm
of $X_i$ is bounded from below by a positive constant. In fact, we have
$$
0 < \frac{1}{C} \leq \Vert X_i \Vert_{m-1,I}
$$
for every~$i \in \N$. Owing to Ascoli-Arzela theorem, and modulo passing to a subsequence, the sequence of vector fields
$\{ X_i \}$ converges in the $C^{r-1}$-topology towards a $C^{m-1}$-vector field $X$. Furthermore, $X$ is not identically zero
since it must verify $\Vert X \Vert_{m-1,I} \geq 1/C >0$. Now a standard application of Euler polygonal method shows that the
vector field $X$ is contained in the $C^{m-1}$-closure of $\Gamma$ in the sense of Definition~\ref{flowclosurepseudogroup}.
The lemma is proved.\qed

Lemma~\ref{C1vectorfield} will often be used in the context where $m=2$.
The method originally put forward in \cite{reb1} is summarized by Proposition~\ref{correctingasequence} below; see also
\cite{RebMMJ} and \cite{advances}.

\begin{prop}
\label{correctingasequence}
Consider a pseudogroup $\Gamma$ consisting of maps from an interval $J \subset \R$ to $\R$ and satisfying
the two conditions below.
\begin{itemize}
  \item There is a sequence of elements $g_i \in \Gamma$, $g_i \neq {\rm id}$ for every~$i \in \N$, all of whose elements
  are defined on $J$. Moreover, this sequence converges to the identity on the $C^m$-topology on $J$.
  \item There is an element $f \in \Gamma$ possessing a hyperbolic fixed point $p \in J$.
\end{itemize}
Then there is an open interval $I \subset J$ containing $p$ and a sequence of elements $\{ \tilde{g}_i \}$ in $\Gamma$
satisfying the conditions of Lemma~\ref{C1vectorfield} (in particular, all the diffeomorphisms $\tilde{g}_i$ are defined
on $I$ and none of them coincides with the identity on $I$).
\end{prop}

\noindent {\it Proof}. We shall sketch the argument since extensions of this basic idea will play an
important role in Sections~4 and~5. It suffices to consider the case $m=2$.
By assumption, we have $f(p) =p$ and $f'(p) = \lambda \in (0,1)$. Since $f$ is analytic,
there is a local coordinate~$x$ around $p$ where $f(x) = \lambda \, x$ \cite{Stenrberg}.
Let then $I \subset J$ be an interval containing $p$
whose closure is contained in the domain of definition of the coordinate~$x$. First, we have the following:

\noindent {\it Claim 1}. Without loss of generality, we can assume that $g_i (p) \neq p$ for every~$i$.

\noindent {\it Proof of Claim~1}. Suppose that $g_i (p) =p$ for all but finitely many~$i$. If, for some large enough~$i$,
we have $g_i' (p) =1$ then by considering elements of the form $\{ \lambda^{-N} \, g_i (\lambda^N x) \}$ (with $i$ fixed),
we can obtain a Shcherbakov-Nakai vector field defined
{\it on a neighborhood of $p$}\, and contained in the $C^{\infty}$-closure of $\Gamma$.
The existence of this vector field actually suffices for our purposes, yet we point out that the sequence
of elements $\{ \lambda^{-N} \, g_i (\lambda^N x) \}$ verifies the conditions of our statement.

It follows from the preceding that the proposition holds provided that there is some $g_i$ not commuting with $f$
and satisfying $g_i (p) =p$. Hence it only remains to consider the possibility of having all the diffeomorphisms
$g_i$ commuting with $f$ and satisfying $g_i (p) =p$ (modulo dropping finitely many terms of the initial sequence). Since $g_i$
commutes with $f$, it must be given on $I$ and in the coordinate~$x$ by $g_i (x) = \lambda_i \, x$. However the sequence
$\{ \lambda_i \}$ converges to~$1$ since $\{ g_i \}$ converges~$C^2$ (in fact $C^{\infty}$) to the identity. In other words,
the sequence $\{ g_i \}$ satisfies the requirements in our statement.\qed

Considering the last possibility discussed in the proof of the above claim, the reader will note that the $C^1$-closure
of $\Gamma$ contains a flow consisting of linear maps $x \mapsto \Lambda \, x$
for every $\Lambda \in \R^{\ast}$. Indeed, for every~$i$, $\lambda_i \neq 1$ since $g_i \neq {\rm id}$. There follows that
the multiplicative group of $\R^{\ast}$ generated by the collection of all $\lambda_i$ is dense in $\R^{\ast}$ what, in turn,
ensures that the mentioned vector field lies in the $C^1$-closure (indeed in the $C^{\infty}$-closure) of $\Gamma$.

Going back to the proof of our proposition, in what follows we assume that $g_i (p) \neq p$ for
every~$i \in \N$. Next, let $\kappa_i$ be a sequence of
positive integers going to infinity to be determined later. Set
$$
\tilde{g}_i = f^{-\kappa_i} \circ g_i \circ f^{\kappa_i} = \lambda^{-\kappa_i} \, g_i(\lambda^{\kappa_i} x) \, .
$$
Note that the second derivative $\tilde{g}_i''$ of $\tilde{g}_i$ at a point $x$ is simply
$\tilde{g}_i'' (x) = \lambda^{\kappa_i} \, g'' (\lambda^{\kappa_i} x )$
provided that both sides are defined. This simple formula shows that $\sup_{x \in I} \Vert \tilde{g}_i'' (x) \Vert$
decreases as $\kappa_i$ increases. On the other hand
the absolute value of $\lambda^{-n} \, g_i (0)$ increases monotonically with $n$ and becomes unbounded as
$n \rightarrow \infty$ since $g_i (0) \neq 0$.
Therefore the $C^1$-norm of $\tilde{g}_i - {\rm id}$ on $I$ also increases with~$n$.
Thus, for every~$i$ fixed, we can find $\kappa_i \in \N^{\ast}$ so that the following estimate holds:
$$
\sup_{x \in I} \Vert \tilde{g}_i'' (x) \Vert \leq \sup_{x \in I} \{ \Vert \tilde{g}_i - {\rm id} \Vert
+ \Vert \tilde{g}_i' - 1 \Vert \} \, .
$$
For these choices of $\kappa_i$ we immediately obtain
$$
\Vert \tilde{g}_i - {\rm id} \Vert_{2,I} < 2 \Vert \tilde{g}_i - {\rm id} \Vert_{1,I}
$$
proving the proposition.\qed

\subsection{Topological dynamics of locally non-discrete subgroups}

The material presented in this section is very closely related to the description of the topological
dynamics associated with groups generated by diffeomorphisms close to the identity obtained in \cite{ghys}. In fact,
our purpose is to prove the following:

\begin{prop}
\label{locallynondiscret}
Let $G \subset \dif$ be a locally $C^2$-non-discrete group. Then either $G$ has a finite orbit or every orbit of
$G$ is dense in $S^1$. Moreover, the set of points in $S^1$ having finite orbit under $G$ is itself finite. Finally,
if $I$ is a connected interval in the complement of this set and $G_I$ denote the subgroup of $G$ consisting of those
diffeomorphisms fixing $I$, then the action of $G_I$ on $I$ has all orbits dense in $I$.
\end{prop}

Since our assumptions are slightly more general than those used in \cite{ghys}, we shall provide below a detailed
proof for Proposition~\ref{locallynondiscret}. We begin by recalling a well-known proposition;
see for example \cite{candel}, \cite{navas}.

\begin{prop}
\label{dynamicswithcantorset}
Denote by ${\rm Homeo}\, (S^1)$ the group of homeomorphisms of the circle and consider a subgroup $G \subset
{\rm Homeo}\, (S^1)$. Then one of the following holds:
\begin{enumerate}
  \item The group $G$ possesses a finite orbit in $S^1$.
  \item The $G$-orbit of every point $p \in S^1$ is dense in $S^1$.
  \item There is a Cantor set $K \subset S^1$ invariant by $G$ and such that the $G$-orbit of every
  point $p \in K$ is dense in $K$. This set is unique and contained in the closure of the $G$-orbit
  of every point $p \in S^1$.\qed
\end{enumerate}
\end{prop}

Consider now a subgroup $G$ of $\dif$. If $G$ possesses a finite orbit then the statement of
Proposition~\ref{dynamicswithcantorset} can be made more accurate as follows. The presence of a finite
orbit implies that the rotation numbers of the elements in $G$ take values in some finite set. In turn,
the subgroup $G_0$ of $G$ consisting of those diffeomorphisms fixing every point in the mentioned finite orbit
has finite index in $G$. Thus, $G_0$ is not reduced to the identity unless $G$ is a finite group. Assuming
that $G$ is not finite and choosing $g \in G_0$, $g \neq {\rm id}$, there follows that
the set {\it of all points in $S^1$ possessing finite orbit under $G$}\, must be finite since it is contained
in the set of {\it fixed points of $g$}. Hence, we have proved:

\begin{coro}
\label{finiteorbitsfiniteset}
Assume that the group $G$ is infinite but has a finite orbit $\mathcal{O}_p$. Denote by
$\mathcal{O}_{G}^{< \infty} \subset S^1$ the set consisting of those points $q \in S^1$ whose orbit
under $G$ is finite. Then $\mathcal{O}_{G}^{< \infty}$ is a finite set.
In particular, $G$ possesses a finite index subgroup $G_0$ whose elements fixes
every single point in $\mathcal{O}_{G}^{< \infty}$.\qed
\end{coro}

Dealing with subgroups of $\dif$ having finite orbit will naturally involve
groups of analytic diffeomorphism of the interval $[0,1]$ (i.e. the group of diffeomorphisms from $[0,1]$ to
$[0,1]$ fixing the endpoints $0$ and $1$). In this direction, the following statement is also well-known
and attributed to G. Hector (see \cite{ghys} for a proof).

\begin{prop}
\label{Hector}
{\rm ({\bf G. Hector})}\,
Let $G_I$ denote a group consisting of orientation-preserving real analytic diffeomorphisms of $[0,1]$.
Suppose that the only points in $[0,1]$ that are fixed for every element in $G_I$ are precisely
the endpoints~$0$ and~$1$. Suppose also that $G$ is neither trivial nor an infinite cyclic group.
Then the orbit of every point $p \in (0,1)$ is dense in $(0,1)$.\qed
\end{prop}

We are now able to prove Proposition~\ref{locallynondiscret}.

\vspace{0.2cm}

\noindent {\it Proof of Proposition~\ref{locallynondiscret}}. Most of the proof amounts to showing that a subgroup
$G \subset \dif$ leaving invariant a Cantor set $K \subset S^1$ must be locally $C^2$-discrete. We begin by
proving this assertion.

Consider then a group $G \subset \dif$ preserving a Cantor set $K$. Proposition~\ref{dynamicswithcantorset} ensures that $K$
is the unique minimal set of $G$ in $S^1$. Furthermore $K$ and the whole of $S^1$ are the only non-empty closed subsets of $S^1$
that are invariant by $G$.

Suppose for a contradiction that $G$ is locally $C^2$-non-discrete. In other words, suppose the existence of an
interval $I \subset S^1$ along with a sequence of elements in $\{ g_i \} \subset G$
satisfying the following:
\begin{enumerate}
  \item $g_i \neq {\rm id}$ for every~$i \in \N$. Since $G$ is constituted by analytic diffeomorphisms, this condition
  also implies that the restriction $g_{i \vert_I}$ of $g_i$ to $I$ does not coincide with the identity on $I$, for every $i \in \N$.
  \item The sequence of restricted maps $g_{i \vert_I} : I \rightarrow S^1$ converges to the identity on the $C^2$-topology over $I$.
\end{enumerate}

Now we have:

\noindent {\it Claim}. The intersection $I \cap K$ is not empty.

\noindent {\it Proof of the Claim}. Suppose that $I \cap K = \emptyset$ and denote by $\widetilde{I}$ the connected component
of $S^1 \setminus K$ containing $I$. The endpoints of $\widetilde{I}$ belong to $K$ and are automatically fixed by every element
of the subgroup $G_I$ of $G$ defined by
$$
G_I = \{ g \in G \; ; \; \; g(\widetilde{I}) \cap \widetilde{I} \neq \emptyset \, \} \, .
$$
Thus, modulo dropping finitely many terms of the sequence $\{ g_i \}$, we can assume that every $g_i$ fixes a chosen
endpoint $p$ of $\widetilde{I}$. Consider a neighborhood $U$ of $p$ and the pseudogroup $\Gamma_U$ induced on $U$ by restrictions of
elements in $G_I$. Note also that $G_I$ fixes $p$ so that we can also consider the group $\Gamma_p$ of germs at $p$ of elements
in $\Gamma_U$. Since $p \in K$ and $K$ is invariant by $G$, it follows that the $C^r$-closure of $\Gamma_U$ contains neither
(standard) Shcherbakov-Nakai vector fields (asymptotically defined on an one-sided interval starting at $p$) nor vector fields
defined on neighborhood of $p$. Owing to Lemma~\ref{preliminarieslemma1},
we conclude that $\Gamma_p$ is infinite cyclic since it cannot be trivial
(the germs at $p$ of the diffeomorphisms $g_i$ belong to $\Gamma_p$). On a fixed neighborhood of~$p$, and for every $i \in \N$,
the diffeomorphism $g_i$ is locally obtained as the map induced by a certain (possibly formal) local flow $\Psi$ at a corresponding
time~$t_i$. The additive subgroup of $\R$ generated by the times~$t_i$ must be discrete, otherwise the local flow $\Psi$
is actually defined for all~$t \in \R$ and the associated analytic vector field is in the closure of~$\Gamma_p$. This is impossible
since $p \in K$ and $K$ is a Cantor set invariant by~$\Gamma_p$. Being discrete, the subgroup of $(\R, +)$ generated by the
times~$t_i$ has a generator~$t_0 >0$. Thus, the dynamics of the group $G_I$ on $\widetilde{I}$ consists of the iterations of
a single diffeomorphism having the endpoints of $\widetilde{I}$ fixed. In particular, the orbit of every point in
$\widetilde{I}$ by the diffeomorphism in question converges to a fixed point of this diffeomorphism.
This contradicts the existence of a sequence of elements in $G_I$ converging to the identity on~$I$.\qed

To complete the proof of the proposition, we proceed as follows. According to a classical theorem due to Sacksteder \cite{candel}, \cite{navas},
there is a point $p \in K$ and a diffeomorphism $f \in G$ such that $f(p) =p$ and $0 < \vert f' (p) \vert <1$. Since $I \cap K$ is
not empty and the dynamics of $G$ on $K$ is minimal, there is no loss of generality in supposing that $p \in I \cap K$. Now, by
considering the pseudogroup $\Gamma$ generated on $I$ by $f$ and by the sequence of maps $g_{i \vert_I}$, Proposition~\ref{correctingasequence}
ensures the existence of a nowhere zero vector field $X$ defined about $p$ and contained in the $C^1$-closure of $\Gamma$.
This yields a contradiction since $K$ is a Cantor set supposed to be invariant by $G$ and, hence, by $\Gamma$.
The resulting contradiction then proves our claim that a locally $C^2$-non-discrete group $G \subset \dif$ cannot leave a Cantor
set $K \subset S^1$ invariant.

To complete the proof of the proposition, we just need to further discuss the case in which $G$ has a finite orbit.
The very assumption that $G$ is locally $C^2$-non-discrete implies that $G$ cannot be finite. Thus the set
$\mathcal{O}_{G}$ of Corollary~\ref{finiteorbitsfiniteset} is finite. Let $I$ be a connected component
of $S^1 \setminus \mathcal{O}_{G}$ and consider the subgroup $G_I$ of $G$ consisting of diffeomorphisms fixing~$I$.
To finish the proof of Proposition~\ref{locallynondiscret} is suffices to check that the action of $G_I$ on $I$
has all orbits dense. Owing to Proposition~\ref{Hector}, if this does not happen then $G_I$ must be infinite cyclic.
Assuming that $G_I$ is infinite cyclic, this group is also locally non-discrete. Lemma~\ref{preliminarieslemma1} then
allows us to conclude that the orbits of $G_I$ on $I$ are still dense. Proposition~\ref{locallynondiscret} is established.\qed

\section{Rigidity in the presence of points with large stabilizers and related cases}

The purpose of this section is to prove Theorem~A in some specific cases related, for example, to the existence of finite
orbits for a non-solvable group (say $G_1$). We shall also settle the case in which $G_1$ is an actual
solvable group. This material will reduce the proof of Theorem~A to a {\it generic situation}\, where, roughly speaking, the
group $G_1$ is not solvable and every point in $S^1$ has cyclic (possibly trivial) stabilizer. This generic situation
is, however, substantially harder and will be detailed in the subsequent sections of this paper.

In the sequel, consider a locally $C^2$-non-discrete subgroup $G_1$ of $\dif$. Then fix
an interval $I \subseteq S^1$ and a sequence $\{ g_{1,i} \}$ of elements in $G_1$
whose restrictions $\{ g_{{1,i} \vert I} \}$ to $I$ converge to the identity in the $C^2$-topology
(with $(g_{1,i} \neq {\rm id})$ for every $i$).
Next, let $G_2$ be another subgroup of $\dif$ that happens to be {\it topologically conjugate}\, to $G_1$. The reader
is reminded that $h$ is supposed to be orientation-preserving.

Having fixed the sequence $\{ g_{{1,i} \vert I} \}$, for every $r \in \N$ we
consider the subgroup $G_{1,r} \subset G_1$ generated by the elements
$g_{1,1} , \ldots , g_{1,r}$ (notation: $G_{1,r} = \langle g_{1,1} , ... , g_{1,r}\rangle$). In
the subsequent discussion, we shall be allowed to ``redefine'' the sequence $\{ g_{{1,i} \vert I} \}$ by dropping
finitely many terms of it and then setting $g_{1,i} = g_{1, i+i_0}$ for every $i \in \N$ and for a certain $i_0 \in \N$.

Throughout this section the group $G_1$ is assumed to be non-abelian. Furthermore,
unless otherwise mentioned, the following condition is also supposed to hold:

\noindent ($\ast$) {\it For every $r \in \N$, the group $G_{1,r}$ possesses a finite orbit while these groups are not finite themselves.}

\begin{lemma}
Modulo redefining the sequence $\{ g_{{1,i} \vert I} \}$ by dropping finitely many terms of it, there is a
finite set $P=\{p_1, ... , p_l\} \subset S^1$ whose points
$p_j$, $i=1 ,\ldots , l$, are fixed points for all the groups $G_{1,r}$.
\end{lemma}

\noindent {\it Proof}.
Let $P_1 \subset S^1$ be the set of points having finite orbit for $G_{1,1}$. Owing to Corollary~\ref{finiteorbitsfiniteset},
the set $P_1$ consists of finitely many points. Naturally, for every $r \geq 1$, the set of points with finite orbit under the
group $G_{1,r}$ is contained in $P_1$ since $G_{1,1} \subset G_{1,r}$. Denoting by $P_r \subset S^1$ the set of points having
finite orbit under $G_{1,r}$, we have $P_1 \supset P_2 \supset \cdots$ so that the intersection
$$
P = \bigcap_{r=1}^{\infty} P_r
$$
is contained in $P_1$. Furthermore this intersection is not empty since our assumption ensures that none of the sets $P_r$ is empty.
Thus, to prove the lemma, it suffices to show that for $i$ sufficiently large, the diffeomorphisms $g_{1,i}$ fixes all points in
$P$. For this let $I_1$ denote a connected component of $S^1 \setminus P$ having non-empty intersection with the open interval
$I$. Since $\{ g_{1,i} \}$ converges to the identity on $I$, for $i$ large enough we must have $g_{1,i} (I_1) \cap I_1 \neq \emptyset$.
Since, on the other hand, the set $P$ is invariant under $g_{1,i}$, it follows at once that $g_{1,i}$ fixes every point in $P$.
The lemma is proved.\qed

Next, let us also consider the group $G_2$ along with the homeomorphism $h$. We begin by letting
$g_{2,i} = h^{-1} \circ g_{1,i} \circ h$ for every $i \in \N$.
We also pose $G_{2,r} = \langle g_{2,1} , \ldots , g_{2,r} \rangle$. Next recall that $P = \{p_1, ... , p_l\}$ and let
$q_j = h^{-1} (p_j)$, for $j=1, \ldots ,l$.
It is cleat that the set $Q = \{ q_1, \ldots ,q_l \}$ is constituted by fixed points of $G_{2,r}$ for every $r \in \N$.

Now let $p_1 \in P$ and $q_1 = h^{-1} (p_1) \in Q$ be fixed. From what precedes, the stabilizer  of $p_1$ (resp. $q_1$)
contains all of $G_{1,r}$ (resp. $G_{2,r}$) for every~$r \in \N$. Now we shall consider a few different possibilities
involving the algebraic structure of the groups $G_{1,r}$.

\begin{prop}
\label{nonsolvablewithfixedpoints}
Modulo choosing $r \in \N$ very large, suppose that the group $G_{1,r}$ is not solvable.
Then the conjugating homeomorphism $h$ coincides with a real analytic diffeomorphism of $S^1$.
\end{prop}

\noindent {\it Proof}.
Let $\Gamma_1$ (resp. $\Gamma_2$) denote the germ of $G_{1,r}$ (resp. $G_{2,r}$) about $p_1$ (resp. $q_1$). Naturally both groups
$\Gamma_1, \, \Gamma_2$ can be identified with non-solvable subgroups of $\difro$ which are (locally) topologically conjugate by a
homeomorphism induced by the restriction of $h$. A result due to Nakai \cite{na2} ensures then that $h$ is real analytic on
a neighborhood of $0 \simeq p_1$. Since $p_1$ is an arbitrary point in $P$, we conclude that $h$ is analytic on a neighborhood
of every point in $P$. Finally, up to choosing $n$ even larger if needed, we can assume that $G_{1,n}$ has dense orbits on the
connected components of $S^1 \setminus P$, cf. Proposition~\ref{Hector}. From this it promptly follows that the local analytic
character of $h$ about points in $P$ extends to all of~$S^1$. The proof of our proposition is over.\qed

Let us now assume that for every $r \in \N$, the group $G _{1,r}$ (and thus $G _{2,r}$) is solvable. We begin with a
general and well-known lemma concerning solvable subgroups of $\dif$.

\begin{lemma}
\label{solvablegroups}
Let $G \subset \dif$ be a solvable subgroup of $\dif$. Then either $G$ has a finite orbit or it is topologically conjugate to a
group of rotations.
\end{lemma}

\noindent {\it Proof}.
Since $G$ is solvable, its action on $S^1$ preserves a probability measure~$\mu$.
Hence the support ${\rm Supp}\, (\mu)$ of $\mu$ is a closed subset of $S^1$ invariant by $G$. Consider a minimal set
$\mathcal{M}$ for $G$ contained in ${\rm Supp}\, (\mu)$. In view of Proposition~\ref{dynamicswithcantorset},
$\mathcal{M}$ must be of one of the following types:
the entire circle, a finite orbit or a Cantor set. Suppose first that $\mathcal{M}$ coincides with
all of~$S^1$. Then by parameterizing the circle by the integral of the Radon-Nikodym derivative, a topological conjugation between
$G$ and a rotation group of $S^1$ can be constructed (in particular $G$ is abelian). In turn, if the support of $\mu$ is a
finite set, then this set is invariant by $G$ so that this group has finite orbits. Hence to finish the proof of the lemma it
suffices to check that $\mathcal{M}$ cannot be a Cantor set. This last
assertion follows from Sacksteder's theorem; see \cite{candel}, \cite{navas}. In other words, if $\mathcal{M}$
is a Cantor set, then there is an element $g \in G$ and a point $p \in \mathcal{M}$ such that $p$ is a hyperbolic fixed point for $g$.
Now, since $g$ preserves $\mu$ and $p \in {\rm Supp}\, (\mu) = \mathcal{M}$, it follows that the point $p$ must
have strictly positive $\mu$-mass. However the measure $\mu$ is invariant by $G$ and finite which, in turn, forces the orbit of
$p$ to be finite itself thus completing the proof of the lemma.\qed

\begin{prop}
\label{something1}
Up to dropping finitely many terms of the sequence $\{ g_{1,i} \}$, suppose
that $G _{1,r}$ is an infinite solvable group for every~$r \in \N$. Suppose, in addition, the existence of $r_0$ for which
$G_{1,r_0}$ has a finite orbit
for some $r \in \N$.
Then the homeomorphism $h$ conjugating $G_1$ to $G_2$ coincides with an analytic diffeomorphism of~$S^1$.
\end{prop}

\noindent {\it Proof}.
Let $p \in S^1$ be a point whose orbit under $G_{1,r_0}$ is finite. Consider then the stabilizer $\Gamma_p$ of $p$ in $G _{1,r_0}$.
The group $\Gamma_p$ is naturally identified to an infinite solvable subgroup of $\difro$. Besides, with suitable
identifications, the restriction of $h$ to a
neighborhood of~$p$ topologically conjugates $\Gamma_p$
to another subgroup $\Gamma_q$ of $\difro$. Again the proof
of the proposition becomes reduced to checking that the homeomorphism (still denoted by $h$) conjugating $\Gamma_p \subset \difro$
to $\Gamma_q \subset \difro$ must be analytic on a neighborhood of $0 \in \R$. For this, let us consider the following possibilities:

\noindent {\it Case 1}.
Suppose that $\Gamma_p$ (an thus $\Gamma_q$) is not abelian. From the description of solvable subgroups of $\difro$, there follows that
solvable non-abelian subgroups of $\difro$ have elements $f_1,g_1$ satisfying the following conditions:
\begin{itemize}
\item $f_1$ has a hyperbolic fixed point at $0 \in \R$.
\item $g_1$ is tangent to the identity at $0 \in \R$ (though $g \neq {\rm id}$).
\end{itemize}
According to Section~2, the local diffeomorphisms $f_1,g_1$ can be combined to construct a (non-identically zero) analytic vector field
$X_1$ defined {\it on a neighborhood of $0 \in \R$}\, and contained in the closure of $\Gamma_p$. A similar vector field $X_2$ can be
defined by means of the elements $f_2 = h^{-1} \circ f_1 \circ h$ and $g_2 = h^{-1} \circ g_1 \circ h$ of $\Gamma_q$. By using
the fact that $h$ conjugates the actions of $\Gamma_p, \, \Gamma_q$, it follows from the indicated constructions that $h$
conjugates $X_1$ to $X_2$ in a time-preserving manner. Thus $h$ must be analytic about $0 \in \R$ and this establishes the
proposition in the first case.

\noindent {\it Case 2}.
Suppose now that $\Gamma_{p}$ (and thus $\Gamma_{q}$) is an infinite abelian subgroup. Thus so is $G_{1,r_0}$. Therefore,
owing to Case~1 and up to dropping finitely many terms from the sequence $\{ g_{1,i} \}$,
we can assume that $G_{1,r}$ is abelian for every $r \in \N$. We then define the abelian group
$G_{1,\infty}$ be letting
$$
G_{1,\infty} = \bigcup_{r=1}^{\infty} G_{1,r} \; .
$$
Next fix an element $f \in G_{1,r_0} \subset G_{1,\infty}$,
$f \neq {\rm id}$, having a fixed point. Since $G_{1,\infty}$ is abelian, the finite set consisting of fixed points of $f$ is invariant
under $G_{1,\infty}$. In particular, for $i$ large enough, every diffeomorphism in the sequence $\{ g_{1,i} \}$ must
fix every point fixed by $f$. In particular, the stabilizer $\Gamma_{p, \infty}$ of $p$ in $G_{1,\infty}$ is an abelian
group contained the diffeomorphisms $g_{1,i}$ for large enough~$i$. In particular $\Gamma_{p, \infty}$ is
non-discrete since $\{ g_{1,i} \}$ converges
to the identity on~$I$. Given that $\Gamma_{p, \infty}$ is abelian non-discrete, we can resort to Lemma~\ref{preliminarieslemma1} in Section~2
to produce vector fields $X_1, \, X_2$ defined around $p$ and $q$ respectively that are
conjugated by $h$ in a time-preserving manner. Hence, it follows
again that $h$ must be analytic about $0 \in \R$. The proposition is proved.\qed

To finish this section we shall establish a last reduction to the proof of Theorem~A in the form of
Proposition~\ref{reductiontothelastcase}. To state this reduction, recall that $I \subset S^1$ is a fixed interval
for which $G_1$ contains a sequence of elements
$\{ g_{1,i} \}$, $(g_{1,i} \neq {\rm id})$,
whose restrictions $\{ g_{{1,i} \vert I} \}$ to $I$ converge to the identity in the $C^2$-topology.

\begin{prop}
\label{reductiontothelastcase}
To prove Theorem~A, we can assume without loss of generality that the following hold:
\begin{itemize}
  \item There is $N \in \N$ for which the group generated by $\{ g_{1,1}, \ldots , g_{1,N} \}$
  is not solvable.

  \item For every given $\varepsilon >0$ (and up to dropping finitely many terms from the sequence $\{ g_{1,i} \}$),
  all the diffeomorphisms $g_{1,1}, \ldots , g_{1,N}$ are $\varepsilon$-close to the identity in
  the $C^2$-topology on the interval~$I$.

  \item No point $p \in S^1$ is simultaneously fixed by all the diffeomorphisms $g_{1,1}, \ldots , g_{1,N}$.

  \item In general, every finite subset generating a non-solvable subgroup of $G_1$ cannot have a common
  fixed point.
\end{itemize}
\end{prop}

To prove Proposition~\ref{reductiontothelastcase}, let us assume once and for all that $\varepsilon >0$ is given.
We also assume without loss of generality that the sequence of diffeomorphisms $\{ g_{1,i} \}$ is constituted by
diffeomorphisms $\varepsilon$-close to the identity in the $C^2$-topology on the interval~$I$.

Assume there is $r_0 \in \N$ such that the group $G_{1,r_0}$ is not solvable. Owing to Proposition~\ref{nonsolvablewithfixedpoints},
Theorem~A holds provided that the non-solvable group $G_{1,n_0}$ possesses a finite orbit. More generally,
Proposition~\ref{nonsolvablewithfixedpoints} also justifies the last assertion in the statement of
Proposition~\ref{reductiontothelastcase}.

Summarizing the above discussion, to establish Proposition~\ref{reductiontothelastcase} it suffices to show that
Theorem~A holds provided that all the groups $G_{1,r}$, $r \in \N^{\ast}$, are solvable. This will be our aim in the sequel.

To begin with, recall the general fact that every finite subgroup of $\dif$ is analytically conjugate to a rotation group.
On the other hand, we also know that every {\it infinite solvable group}\, having no finite orbit
is topologically conjugate to a rotation group, cf. Lemma~\ref{solvablegroups}. By virtue of Proposition~\ref{something1},
we can therefore assume that each group $G_{1,r}$ is abelian and topologically conjugate to a group of rotations.

Consider again the group $G_{1, \infty} = \bigcup_{r=1}^{\infty} G_{1,r}$.
Clearly $G_{1, \infty}$ is an infinite locally non-discrete abelian group all of whose orbits are infinite.
Although it is infinitely generated, the
action of $G_{1, \infty}$ still preserves a probability measure $\mu_{\infty}$. In fact, let $\mu_r$ be a probability
measure invariant by $G_{1,r}$ and take $\mu_{\infty}$ as an accumulation point of the sequence $\{ \mu_r \}$.
The fact that $G_{1,r} \subset G_{1,r+1}$ promptly implies that $\mu_{\infty}$ must be invariant by $G_{1,r}$ for
every $r \in \N$. Since $G_{1, \infty}$ has no finite orbit, it follows that $\mu_{\infty}$ has no atomic component
so that its support must coincide with all of~$S^1$ (recall that the support cannot be a Cantor set thanks to
Sacksteder theorem \cite{candel}, \cite{navas}). Hence,
the Radon-Nikodym derivative of this measure allows us to construct a topological coordinate $H$ on $S^1$ in which
$G_{1,\infty}$ is a group of rotations. Next we have:

\begin{lemma}
\label{densenessofrotations}
In the topological coordinate $H$, the group $G_{1,\infty}$ is a dense subgroup of the group
of all rotations of~$S^1$.
\end{lemma}

\noindent {\it Proof}. Consider the map $\rho : G_{1,\infty} \rightarrow \R/\Z$ assigning to an element
$g \in G_{1,\infty}$ its rotation number. Because $G_{1,\infty}$ is an abelian group, the map $\rho$ is
a homomorphism so that its image $\rho \, (G_{1,\infty}) \subset S^1$ is a dense set of $S^1$ viewed as a multiplicative
group. Moreover, the homomorphism $\rho$ is injective since, in the coordinate $H$, the rotation corresponding to
an element $g\in G_{1,\infty}$ is nothing but the rotation of angle equal to the rotation number of~$G$. The lemma
then follows from the fact that the subgroup $\rho (G_{1,\infty})$ is clearly infinite.\qed

The next lemma is also elementary.

\begin{lemma}
\label{commutingisrotation}
Suppose that $g: S^1 \rightarrow S^1$ is a homeomorphism of the circle that commutes with a dense set $E$ of rotations.
Then $g$ is itself a rotation.
\end{lemma}

\noindent {\it Proof}.
Consider the circle equipped with the standard euclidean metric~$d$ induced from $\R$ by the identification
$S^1 = \R /\Z$. To show that $g$ is a rotation amounts
to check that $g$ is an isometry of~$d$. Hence, chosen an interval $J$ with endpoints $x,y$, we need to show that
the length of $g(J)$ equals to the length of~$J$. If this were not true, then there would exist $J \subset S^1$
such that the length $L(J)$ of $g(J)$ would be strictly smaller than the length $L (g(J))$ of $g(J)$.
Now, since $E$ is a dense set of rotations, we can find an element $\sigma \in E$ such that $\sigma (g(J)) \subset J$.
Thus the map $\sigma \circ g$ maps $J$ strictly inside itself and must therefore have
a fixed point $p \in J \subset S^1$.
Furthermore $\sigma \circ g$ commutes with all rotations in $E$ so that the orbit of $p$ by elements in $E$
must consist of fixed points for $\sigma \circ g$. However, since the orbit of $p$ by all rotations in $E$ is clearly
dense in $S^1$, there follows that $\sigma \circ g$ coincides with the identity. The resulting contradiction
proves the lemma.\qed

Let us close this section with the proof of Proposition~\ref{reductiontothelastcase}

\vspace{0.1cm}

\noindent {\it Proof of Proposition~\ref{reductiontothelastcase}}.
The proof amounts to showing
that the initial sequence of diffeomorphisms $\{ g_{1,i} \} \subset G_1$
can be chosen so as to ensure that for large enough~$r \in \N$ the group $G_{1,r}$ cannot be topologically conjugate to
a group of rotations. For this, consider a finite
generating set $f_{1,1}, \ldots ,f_{1,s}$ for $G_1$. Given the initial sequence $\{ g_{1,i} \} \subset G_1$, we consider all diffeomorphisms
of the form $g_{1,j,i} = f_{1,j}^{-1} \circ g_{1,i} \circ f_{1,j}$ where $j \in \{ 1, \ldots ,s \}$ and $i \in \N$.
Next, the indices $j,i$ can be reorganized to ensure that all the
diffeomorphisms $g_{j,i}$ are actually contained in the initial sequence $\{ g_{1,i} \}$. With this new definition
of the sequence $\{ g_{1,i} \}$, the following holds:

\noindent {\it Claim}. The group $G_{1, \infty}$ is no longer topologically conjugate to a group of rotations.

\noindent {\it Proof}.
By construction the group $G_{1, \infty}$ consists of elements having the form $f_{1,j}^{-1} \circ \tilde{g}_k \circ f_{1,j}$,
where $\tilde{g}_k \in \dif$ is a certain sequence of diffeomorphisms converging to the identity on~$I$ (in the $C^2$-topology).
Suppose for a contradiction that $G_{1, \infty}$ is abelian without finite orbits.
Now fixed $k$, the elements $g_k$ and $f_{1,j}^{-1} \circ g_k \circ f_{1,j}$, $j=1, \ldots ,s$ have all the same rotation number.
What precedes then ensures that all these elements are the same. Indeed, it was seen that the ``rotation number
homomorphism'' from $G_{1, \infty}$ to~$S^1$ is one-to-one. In other words, for every $k \in \N$ and every~$j=1, \ldots ,s$
the diffeomorphisms $g_k$ and $f_{1,j}$ do commute.

Now recall the existence of a topological coordinate $H$
where $G_{1, \infty}$ is identified to a group of rotations that happens to be dense in the group of all rotations of~$S^1$.
Let $\Gamma$ be the subgroup of $G_{1, \infty}$ generated by all the elements $g_k$, $k \in \N$ and note that $\Gamma$
is dense the group of all rotations of~$S^1$ as well. Finally, always working in the coordinate~$H$, the generators $f_{1,1},
\ldots ,f_{1,s}$ of $G_1$ commute with all elements in $\Gamma$. Lemma~\ref{commutingisrotation} then ensures that every $f_{1,j}$
is itself another rotation in the coordinate~$H$. Hence the group $G_1$ must be abelian and this yields the desired contradiction.\qed

Now the proposition results from the repeating word-by-word the preceding discussion.\qed

\section{Convergence estimates for sequences of commutators}\label{section.convergence}

This section is devoted to providing an algorithmic way to construct diffeomorphisms converging to the identity
on a suitably fixed interval. These constructions will allow for a more effective use of the assumption that
our groups are locally non-discrete and this deserves further comments. Consider a locally $C^2$-non-discrete
group $G \subset \dif$. In other words, there is an interval $I \subset S^1$ and a sequence of elements
$\{ g_i \} \subset G$ satisfying the conditions of Definition~\ref{localnondiscreteness}. An inconvenient point
concerning this definition lies in the fact that the sequence $\{ g_i \}$ is {\it a priori given}\, and this prevents us
from having any additional control on the diffeomorphisms $g_i$. For example, we have no information whatsoever on the
higher order derivatives of $g_i$ and, in particular, no information on the growing rate of the sequence
$\Vert g_i \Vert_3$, where $\Vert \, . \, \Vert_3$ stands for the $C^3$-norm. In the context
of Theorem~A, if $\{g_{1,i}\}$ is a sequence as above for the group $G_1$, then the corresponding sequence
$g_{2,i} = h^{-1} \circ g_{1,i} \circ h$ of elements in $G_2$ is known to converge to the identity only in the $C^0$-topology.
With a purpose of deriving non trivial implications on the regularity of $h$, it is natural to look for sequences
as above such that $\{ g_{2,i} \}$ converges in topologies stronger than the $C^0$-topology. The main immediate virtue
of the construction presented below is to yield some control
on the growing rate of the sequence formed by higher order derivatives of $g_i$. A
further application of our construction will be given in the form of partial answer to some questions raised
in \cite{deroinETDS}, cf. Appendix (Section~7).

Considering then a locally $C^2$-non-discrete group $G \subset \dif$. We want to {\it explicitly construct}
sequences of elements in $G$ satisfying the conditions of Definition~\ref{localnondiscreteness} where,
by {\it explicitly constructing},
it is meant an algorithmic procedure beginning with a finite set of suitably chosen elements in $G$. The very
nature of the algorithm will immediately yield elementary estimates that will be needed in the proof of
Theorem~A.

Owing to Proposition~\ref{reductiontothelastcase}, we fix some interval $I \subset S^1$ and a collection $S \subset G$
of elements $\overline{g}_{1}, \ldots , \overline{g}_{N}$ generating a non-solvable subgroup.
The diffeomorphisms $\overline{g}_i$, $i=1, \ldots , N$
are also assumed to be $\varepsilon$-close to the identity in the $C^2$-topology on the interval~$I$, where the
value of $\varepsilon >0$ will be fixed only later. The idea to produce a sequence of diffeomorphisms
converging to the identity out of the finite set $S = \{ \overline{g}_{1}, \ldots , \overline{g}_{N} \}$
consists of iterating commutators. This is a slight refinement of the method employed by
Ghys in \cite{ghys} which relies on a fast iteration technique. Indeed, the difficulty in proving convergence to
the identity of sequences of iterated commutators lies in the fact that an estimation of the $C^m$-norm of a
commutator $[f_1 ,f_2] = f_1 \circ f_2 \circ f_1^{-1} \circ f_2^{-1}$ requires estimations on the $C^{m+1}$-norm
of $f_1,\, f_2$. To establish the convergence of a sequence of ``iterated commutators'' becomes therefore tricky as at each step
there is an intrinsic loss of one derivative. It is thus natural to try to overcome this difficulty by means of
some suitable fast iteration scheme and this is the idea
of Ghys \cite{ghys} who uses holomorphic extensions and the uniform convergence of those. In the holomorphic
context, however, Cauchy formula enables us to substitute the loss of one derivative by the loss of a portion
of the domain of definition: hence it is enough to ensure that the portion of domain lost at step of
the iteration scheme becomes smaller and smaller so that some uniform domain is kept
at the end.

Since we will work only with $C^2$-convergence the same fast iteration scheme is not available, albeit adaptations
are still possible. Yet, we prefer to introduce a slightly more elaborated iterative procedure which has the advantage
of avoiding fast convergence estimates.
The idea is to add a step of {\it renormalization}\, at each stage of the commutator iteration. This renormalization
step has a regularizing effect on derivatives of order two or greater. A simplified version of the same
idea was already used in the proof of Proposition~\ref{correctingasequence}. One advantage of our procedure is to
avoid the loss of derivatives; other advantages will become clear in the course of the discussion and these include
the convergence rate to the identity of the resulting sequence; see Remark~\ref{dominatingpolynomials}.

After this brief overview of the material to be developed below, we begin to provide accurate
definitions. We shall work with the pseudogroup generated by $S = \{ \overline{g}_1 ,\ldots ,\overline{g}_N \}$
on the interval~$I \subset \R$ where $\overline{g}_1 ,\ldots ,\overline{g}_N$ generate a non-solvable group.
Also, and whereas we shall primarily think of
$\overline{g}_1 ,\ldots ,\overline{g}_N$ as maps defined on $I$, it is sometimes useful to keep in mind
that all of these maps mere restrictions to $I$ of global analytic diffeomorphisms of $S^1$ (still denoted
by $\overline{g}_1 ,\ldots ,\overline{g}_N$, respectively).

Following Ghys~\cite{ghys}, let us associate to the set
$S= \{ \overline{g}_{1}, \ldots , \overline{g}_{N} \}$ a {\it sequence of sets}\, $S(k)$, $k = 1,2,\ldots$, inductively defined as follows:
\begin{itemize}
\item $S(0)=S$
\item $S(k)$ is the set whose elements are commutators of the form $[\tilde{f_i}^{\pm 1} , \tilde{f}_j^{\pm 1}]$ where
$\tilde{f}_i \in S(k-1)$ and
$\tilde{f}_j \in S(k-1)\cup S(k-2)$ ($\tilde{f}_j \in S(0)$ if $k=1$).
\end{itemize}
Still according to Ghys~\cite{ghys}, the resulting sequence of sets $S(k)$ is never reduced to the identity
since $S= \{ \overline{g}_{1}, \ldots , \overline{g}_{N} \}$ generates a non-solvable group.
This also yields the following:

\begin{lemma}
\label{SKK-1nonsolvable}
For every $k \in \N$, the subgroup generated by  $S(k) \cup S(k-1)$ is non-solvable.
\end{lemma}

\noindent {\it Proof}. Assume
there were $k \in \N$ such that $\Gamma =\langle S(k) \cup S(k-1) \rangle$ is solvable, where
$\langle S(k) \cup S(k-1) \rangle$ stands for the group generated by $S(k) \cup S(k-1)$.
Since $\Gamma \subset \dif$, there follows that $\Gamma$ is, indeed, metabelian, i.e. its derived group
$D^1 \Gamma$ is abelian. Recalling that $D^1 \Gamma$ is the group generated by all commutators of the
form $[\gamma_1, \gamma_2]$ where $\gamma_1, \, \gamma_2 \in \Gamma$, there follows that
the sets $S(k+1)$ and $S(k+2)$ are contained in $D^1 \Gamma$.
Since $D^1 \Gamma$ is abelian, the definition of the sequence of sets $\{ S(k) \}$ promptly implies that the set
$S(k+3)$ must coincide with $\{ {\rm id} \}$. Hence the initial group generated by
$\overline{g}_{1}, \ldots , \overline{g}_{N}$ must be solvable. The resulting contradiction proves the lemma.\qed

By virtue of Proposition~\ref{reductiontothelastcase}, we obtain the following corollary:

\begin{coro}
\label{againnofixedpoints}
In order to prove Theorem~A, we can assume that the elements in $S(k) \cup S(k+1)$ do not share
a common fixed point (and this for every~$k \in \N$).\qed
\end{coro}

From now on, we set $I = [-a,a] \subset \R$, $a >0$, with the obvious identifications. Given $\varepsilon >0$,
we permanently fix a set of diffeomorphisms $\overline{g}_{1}, \ldots , \overline{g}_{N}$ generating a non-solvable group and
$\varepsilon$-close to the identity in the $C^2$-topology on~$I$. The value of $\varepsilon >0$ convenient for
our purposes will only be fixed later. In the rest of the section, these conditions are assumed to hold
without further comments.

Unless otherwise mentioned, in what follows we shall say that $f : I' \subseteq I \subset \R \rightarrow \R$ is a diffeomorphism
meaning that $f$ is a diffeomorphism from $I' \subset \R$ to $f(I') \subset \R$.
Let us begin our discussion by stating a simple general lemma.

\begin{lemma}
\label{lossofaderivative}
Given $\epsilon_0 >0$ small and $m \geq 1$, there is a neighborhood $\mathcal{U}_0^m$ of the identity
in the $C^m$-topology on $I$ such that the commutator
$[f_1,f_2] = f_1 \circ f_2 \circ f_1^{-1} \circ f_2^{-1}$ satisfies the two conditions below
provided that $f_1, f_2$ belongs to $\mathcal{U}_0^m$:
\begin{itemize}
  \item Viewed as an element of the pseudogroup generated by $f_1,f_2$ on $I$, the map $[f_1,f_2]$ is well defined
  on $[-a+5 \epsilon_0 , a - 5\epsilon_0]$.

  \item There is a constant $C>0$ such that
  the $C^{m-1}$-distance $\Vert [f_1,f_2] - {\rm id} \Vert_{r-1, [-a+5 \epsilon_0 , a - 5\epsilon_0]}$
  from $[f_1,f_2]$ to the identity on the interval $[-a+5 \epsilon_0 , a - 5\epsilon_0]$
  satisfies the estimate
  $$
  \Vert [f_1,f_2] - {\rm id} \Vert_{m-1, [-a+5 \epsilon_0 , a - 5\epsilon_0]} < C \, \Vert f_1 -{\rm id} \Vert_{m, [-a,a]} \,
  \Vert f_2 -{\rm id} \Vert_{m, [-a,a]}
  $$
  where $\Vert f_1 -{\rm id} \Vert_{m, [-a,a]}$ (resp. $\Vert f_2 -{\rm id} \Vert_{m, [-a,a]}$) stands for the
  $C^m$-distance from $f_1$ (resp. $f_2$) to the identity on the interval $I = [-a,a]$.\qed
\end{itemize}
\end{lemma}

The reader will note that the constant $C$ in the above lemma depends only on the neighborhood $\mathcal{U}_0^m$. In particular
$C$ does not increase {\it when the neighborhood is reduced}.

We now focus on the case $m=2$ (see Appendix for a more general discussion). Since we can always
reduce $\varepsilon >0$, the neighborhood $\mathcal{U}_0^2$ can be chosen as
\begin{equation}
\mathcal{U}_0^2 = \{ f \in C^2 ([-a,a]) \; ; \; \; \Vert f - {\rm id} \Vert_{2, [-a,a]} < \varepsilon \} \label{definitionU20}
\end{equation}
where $C^2 ([-a,a])$ stands for the space of $C^2$-functions defined on $[-a,a]$ and taking values in $\R$.
For this neighborhood $\mathcal{U}_0^2$, the constant provided by Lemma~\ref{lossofaderivative} will be denoted
by $C$ and the value of $C$ does not increase when $\varepsilon$ decreases.

Now we state a simple complement to Lemma~\ref{lossofaderivative}:

\begin{lemma}
\label{estimatesecondderivative}
Up to reducing $\varepsilon >0$, for every pair $f_1 ,f_2 \in \mathcal{U}_0^2$ the second derivative $D^2[f_1,f_2]$ of the commutator
$[f_1,f_2]$ on the interval $[-a+5 \epsilon_0 , a - 5\epsilon_0]$ satisfies the estimate
$$
\sup_{x \in [-a+5 \epsilon_0 , a - 5\epsilon_0]} \left| D^2[f_1,f_2] \right| \leq \, 5 \, \max_{x \in (-a,a)} \{\left| D^2f_1 \right| ,
\left |D^2 f_2 \right| \}
$$
where $D^2f_j$ stands for the second derivative of $f_j$, $j=1,2$.
\end{lemma}

\noindent {\it Proof}. The proof is elementary and we shall summarize the argument. For $j=1,2$, the very definition
of $\mathcal{U}_0^2$ yields (see~(\ref{definitionU20}))
$$
1-\varepsilon \leq  \vert D^1_x f_j \vert \leq 1+\varepsilon \; \; \; {\rm and} \; \; \; \frac{1}{1+\varepsilon} \leq
\frac{1}{\vert D^1x f_j \vert}\leq \frac{1}{1-\varepsilon}
$$
for every $x \in [-a,a]$. Concerning the inverses of $f_1, f_2$, we also have
$$
D^1_x f_j^{-1}  =\frac{1}{D^1_{f_j^{-1} (x)} f_j} \; \; \; {\rm and} \; \; \; D^2_x f_j^{-1}
=-\frac{D^2_{f_j^{-1} (x)} f_j}{[D^1_{f_j^{-1} (x)} f_j]^3} \; .
$$
Next we compute the second derivative of $[f_1,f_2]$ at a point belonging to $[-a+5 \epsilon_0 , a - 5\epsilon_0]$.
In this calculation, the points at which the several derivatives are evaluated will be omitted: since $[f_1,f_2]$ is well
defined on $[-a+5\epsilon_0,a-5\epsilon_0]$, it suffices to know that all these points belong to the interval $(-a,a)$. Since
$D^1[f_1,f_2]=D^1 f_1.D^1 f_2.D^1 f_1^{-1}.D^1 f_2^{-1}$, we obtain
\begin{eqnarray*}
D^2[f_1,f_2] & = & D^2 f_1.(D^1 f_2)^2.(D^1 f_1^{-1})^2.(D^1 f_2^{-1})^2 + D^1 f_1.D^2 f_2.(D^1 f_1^{-1})^2.(D^1 f_2^{-1})^2 +\\
& = & \mbox{} + D^1 f_1.D^1 f_2.D^2 f_1^{-1}.(D^1 f_2^{-1})^2 + D^1 f_1.D^1 f_2.D^1 f_1^{-1}.D^2 f_2^{-1} \, .
\end{eqnarray*}
Therefore on $[-a+5\epsilon_0,a-5\epsilon_0]$, we have
\begin{eqnarray*}
\vert D^2[f_1,f_2] \vert & \leq & \frac{(1+\varepsilon)^2}{(1-\varepsilon)^4} \vert D^2 f_1 \vert +
\frac{(1+\varepsilon)^2}{(1-\varepsilon)^4} \vert D^2 f_2 \vert + \frac{(1+\varepsilon)^2}{(1-\varepsilon)^5} \vert D^2 f_1 \vert
+ \frac{(1+\varepsilon)^2}{(1-\varepsilon)^4} \vert D^2 f_2 \vert \\
& \leq & \left[ \frac{(1+\varepsilon)^2}{(1-\varepsilon)^4} + \frac{(1+\varepsilon)^2}{(1-\varepsilon)^4}
+ \frac{(1+\varepsilon)^2}{(1-\varepsilon)^5} + \frac{(1+\varepsilon)^2}{(1-\varepsilon)^4} \right] \, \max \{D^2f_1,D^2f_2\} \, .
\end{eqnarray*}
Up to choosing $\varepsilon$ sufficiently small, there follows that $\vert D^2[f_1,f_2] \vert
\leq \, 5 \, \max \{ \vert D^2f_1 \vert , \vert D^2f_2 \vert \}$
proving the lemma.\qed

Let us now begin the construction of a sequence of diffeomorphisms in $G$ converging to the
identity in the $C^2$-topology on $I = [-a,a]$.
First recall that non-solvable subgroups of $\dif$ are known to
have elements with hyperbolic fixed points (see for example \cite{argumentghys}).
Let then $F \in G$ be a diffeomorphism satisfying $F(0) =0$ and $F'(0) = \Lambda \in (0,1)$.
The next step is to define a new sequence $\{ \widetilde{S}(k) \}$ of subsets of $G$. The sequence
$\{ \widetilde{S}(k) \}$ will depend on a fixed integer $n \in \N^{\ast}$ which will be omitted in the notation.
To define the sequence $\{ \widetilde{S}(k) \}$ we proceed as follows:
\begin{itemize}
\item $\widetilde{S}(1)$ is the set formed by the commutators having the form
$[F^{-n} \circ \tilde{f}_1 \circ F^n , F^{-n} \circ \tilde{f}_2 \circ F^n ]$ where $\tilde{f}_1, \, \tilde{f}_2 \in S$.
Thus $\widetilde{S}(1) = F^{-n} \circ S(1) \circ F^n$.

\item $\widetilde{S}(k)$ is the set formed by the commutators $[F^{-n} \circ \tilde{f}_1 \circ F^n , F^{-n} \circ \tilde{f}_2 \circ F^n ]$
with $\tilde{f}_1, \tilde{f}_2 \in  \widetilde{S}(k-1)$ and by the
commutators $[F^{-n} \circ \tilde{f}_1 \circ F^n , F^{-2n} \circ \tilde{f}_2 \circ F^{2n} ]$
with $\tilde{f}_1 \in \widetilde{S}(k-1)$ and $\tilde{f}_2 \in \widetilde{S}(k-2)$.
\end{itemize}
In other words, the sequence $\{ \widetilde{S}(k) \}$ verifies
$\widetilde{S}(k) = F^{-kn} \circ S(k) \circ F^{kn}$ for every $k \in \N$.
Taking advantage of the fact that all our local diffeomorphisms have global realizations in $G$ as a group of diffeomorphisms
of the circle, we obtain the following:

\begin{lemma}
\label{nondegeneratingsequence}
The sequence of sets $\widetilde{S}(k)$ never degenerates into $\{ {\rm id} \}$.
\end{lemma}

\noindent {\it Proof}. When all the diffeomorphisms in question are globally viewed as diffeomorphisms
of the circle, the set $\widetilde{S}(k)$ is conjugate to the set $S(k)$, for every $k \in \N$. The
statement follows then from Ghys theorem claiming that the initial sequence
$S(k)$ cannot degenerate into $\{ {\rm id} \}$ provided that $G$ is non-solvable.\qed

The global realizations of our diffeomorphisms ensure that the domain
of definition of elements in $\widetilde{S} (k)$ are always non-empty as every diffeomorphism
is clearly defined on all of~$S^1$. However, going back to our local setting where
the initial $C^2$-maps $\overline{g}_{1} , \ldots , \overline{g}_{N}$ are defined on $[-a,a]$ and where the domains
of definition for their iterates are understood in the sense of pseudogroup, the content of the last statement becomes
unclear. In other words, in the context of pseudogroups, the statement of Lemma~\ref{nondegeneratingsequence} is only
meaningful for those elements in $\widetilde{S}(k)$ having non-empty domain of definition {\it when
viewed as elements of the pseudogroup in question}. In any event,
the estimates developed below will show that this is always the case provided that we start with a
sufficiently small $\varepsilon >0$.

Now consider a sequence $\{ f_j \}$ of elements in $G$ so that $f_j \in S (j)$ and $f_j \neq {\rm id}$
for every~$j \in \N$. Set
$$
g_j = F^{-jn} \circ f_j \circ F^{jn} \, .
$$
Owing to Corollary~\ref{againnofixedpoints}, we can
assume that $f_j (0) \neq 0$ for every $j \in \N$. The central result of this section reads as follows:

\begin{prop}
\label{convergeG1}
Up to starting with a sufficiently small $\varepsilon >0$ and an appropriately chosen $n \in \N^{\ast}$,
the above constructed sequence $\{ g_j \}$ of elements in $G$ converges
to the identity in the $C^2$-topology on all of the interval $[-a,a]$.
\end{prop}

Recall that $\Lambda = F'(0)$. Therefore, modulo changing coordinates, we can
assume that $F(x) = \Lambda \, x$ for every $x \in [-a,a]$,
\cite{Stenrberg}. In these coordinates, $g_j$ becomes $g_j = \lambda^{-jn} \, f_j (\lambda^n \, x)$. Fix $\epsilon_0 >0$ small
(for example $\epsilon_0 = a/20$). We choose
$\varepsilon >0$ and $n \in \N$ so as to fulfil all the conditions below.
\begin{itemize}
  \item[(A)] - The value of $n$ is chosen to be the smallest positive integer for which the following conditions are satisfied:
$$
0 < \lambda^n a < a -5\epsilon_0 \; \; \; {\rm and} \; \; \; \lambda^n < 1/20 \, .
$$

  \item[(B)] - Lemma~\ref{lossofaderivative} holds on $\mathcal{U}_0^2$ for some $C >0$.

  \item[(C)] - $\varepsilon >0$ is small enough to ensure that Lemma~\ref{estimatesecondderivative} holds and that
$$
\varepsilon \, \max \left\{ (\lambda^{-n} +1) C \, , \;   (\lambda^{-n} +1) \right\} < 1/10 \; .
$$
\end{itemize}

\vspace{0.1cm}

\noindent {\it Proof of Proposition~\ref{convergeG1}}. The proof is by induction. First
consider a diffeomorphism $g_1 \in \widetilde{S} (1)$. By assumption, $g_1 = \lambda^{-n} \, f_1 (\lambda^n x)$ for
some $f_1$ given as a commutator $[\overline{g}_i, \overline{g}_j]$ for some $i, j \in \{ 1, \ldots ,N \}$.
Owing to Lemma~\ref{lossofaderivative}, $f_1$ is defined on $[-a + 5\epsilon_0, a-5\epsilon_0]$ when viewed as element of the pseudogroup
generated by $\overline{g}_1, \ldots ,\overline{g}_N$
on $[-a,a]$. Furthermore, the $C^1$-norm of $f_1 - {\rm id}$ on
$[-a + 5\epsilon_0, a-5\epsilon_0]$ satisfies
\begin{equation}
\Vert f_1 - {\rm id} \Vert_{1, [-a + 5\epsilon_0, a-5\epsilon_0]} < C \, \varepsilon^2 \, . \label{firsttobereferred}
\end{equation}
Next observe that $g_1 = \lambda^{-n} \, f_1 (\lambda^n x)$ is defined on all of
$[-a,a]$ since $\lambda^n a < a -5\epsilon_0$. Moreover, we clearly have:
$$
\sup_{x \in [-a,a]} \vert g_1 (x) -x \vert = \sup_{x \in [-a,a]} \vert \lambda^{-n} \, f_1 (\lambda^n x) -x \vert
= \lambda^{-n} \sup_{y \in [-a + 5\epsilon_0, a-5\epsilon_0]} \vert f_1 (y) - y \vert \, .
$$
Similarly
$$
\sup_{x \in [-a,a]} \vert D^1_x g_1  -1 \vert  = \sup_{y \in [-a + 5\epsilon_0, a-5\epsilon_0]} \vert D^1_yf_1 - 1 \vert \, .
$$
In particular, we obtain
\begin{equation}
\sup_{x \in [-a,a]} \vert g_1 (x) -x \vert + \sup_{x \in [-a,a]} \vert D^1_x g_1  -1 \vert < (\lambda^{-n} +1) C \varepsilon^2 \, .
\label{tobereferred2}
\end{equation}
Finally, the second derivative of $g_1$ at a point $x \in [-a,a]$ is such that $D^2_xg_1 = \lambda^n \, D^2_{\lambda^n x} f_1$
so that
\begin{equation}
\sup_{x \in [-a,a]} \vert D^2_x( g_1 - {\rm id}) \vert = \sup_{x \in [-a,a]} \vert D^2_x g_1 \vert <
\lambda^n 5 \max_{x \in (-a,a)} \{\left| D^2_x \overline{g}_i \right| ,
\left |D^2_x \overline{g}_j \right| \} \leq 5 \lambda^n \varepsilon \, , \label{tobereferred3}
\end{equation}
where we have used Lemma~\ref{estimatesecondderivative}. Comparing Estimates~(\ref{tobereferred2}) and~(\ref{tobereferred3}),
there follows that
$$
\Vert g_1 - {\rm id} \Vert_{2, [-a,a]} \leq (\lambda^{-n} +1) C \varepsilon^2  +  5 \lambda^n \varepsilon \leq \frac{\varepsilon}{10}
+ \frac{\varepsilon}{10} + \frac{\varepsilon}{4} = \frac{\varepsilon}{2} \,
$$
where conditions~(A), (B), and~(C) concerning the choices of $\varepsilon$, $n$, and the constant $C$ were used.
In particular, we see that $g_1$ belongs to $\mathcal{U}_0^2$. Since $g_1$ is an arbitrary element of $\widetilde{S} (1)$,
we conclude that $\widetilde{S} (1) \subset \mathcal{U}_0^2$
so that the procedure can be iterated. Consider then
$g_2 = \lambda^{-n} [\tilde{f}_i, \tilde{f}_j] \circ ( \lambda^n \, x)$ where $\tilde{f}_i, \tilde{f}_j$ belong
to $\widetilde{S}(1) \cup \{ \overline{g}_1, \ldots , \overline{g}_N\}$. Repeating word-by-word,
the preceding argument we eventually obtain
$$
\Vert g_2 - {\rm id} \Vert_{2, [-a,a]} \leq  \frac{\varepsilon}{2} \,
$$
(in particular $g_2$ is defined on all of $[-a,a]$). However an element $g_3 \in \widetilde{S} (3)$ can be written as
$g_3 =  \lambda^{-n} [\tilde{f}_i, \tilde{f}_2] \circ \lambda^n \, x$ where $\tilde{f}_i,
\tilde{f}_j$ satisfy
$$
\max \{ \Vert \tilde{f}_i - {\rm id} \Vert_{2, [-a,a]} \, ; \, \Vert \tilde{f}_j - {\rm id} \Vert_{2, [-a,a]} \} < \varepsilon/2 \, .
$$
Therefore, what precedes allows us to conclude that
$$
\Vert g_3 - {\rm id} \Vert_{2, [-a,a]}  < \frac{\varepsilon}{2^2} \, .
$$
Now a straightforward induction shows that
\begin{equation}
\Vert g_{2k} - {\rm id} \Vert_{2, [-a,a]}  < \frac{\varepsilon}{2^k} \, \label{fellby2k}
\end{equation}
and completes the proof of Proposition~\ref{convergeG1}.\qed

\begin{obs}
\label{dominatingpolynomials}
{\rm Consider a sequence $g_1, g_2 , \ldots$ so that $g_k \in \widetilde{S} (k)$ as above. Consider also
the sequence of real numbers given by $\{ \Vert g_{k} - {\rm id} \Vert_{2, [-a,a]} \}$.
Estimate~(\ref{fellby2k}) shows that the subsequence of $\{ \Vert g_{k} - {\rm id} \Vert_{2, [-a,a]} \}$
formed by those $g_k$ with even order decays at least as $1/\sqrt{2}^k$. In fact,
it can be shown that the entire sequence $\{ \Vert g_{k} - {\rm id} \Vert_{2, [-a,a]} \}$ decays faster than
$\Theta^k$ for every a priori given $\Theta >0$. Indeed, the choice of $\varepsilon >0$ made in condition~(C) can be modified
by replacing the $1/10$ on the right side of the corresponding estimate by a sufficiently small $\delta >0$. Note that
this change does not affect either $n$ or the constant~$C$ whereas it allows us to obtain a finer estimate
than $\varepsilon /2$ for $\Vert g_1 - {\rm id} \Vert_{2, [-a,a]}$. A standard induction argument then yields
a faster exponential decay for the sequence $\{ \Vert g_{k} - {\rm id} \Vert_{2, [-a,a]} \}$.
In turn, we have show that every element $g_k$ in $\widetilde{S} (k)$ satisfies
$\Vert g_{2k} - {\rm id} \Vert_{2, [-a,a]}  < \varepsilon /2^k$ so that there is
$k_0 \in \N$ for which every element in $\widetilde{S} (k)$
satisfies the estimate in condition~(C) with $\delta$ in the place of $1/10$.
Thus, up to dropping finitely many terms, the sequence $\{ g_k \}$ converges to the identity faster than $\Theta^k$.
Since only finitely many terms have been dropped, there follows that the initial sequence $\{ g_k \}$
converges to the identity faster than $\Theta^k$. This simple observation will be useful in the next section.}
\end{obs}

\section{Expansion, bounded distortion and rigidity}

In this section we shall complete the proof of Theorem~A up to
Proposition~\ref{FinishingJob-11.3} whose proof is deferred to the next section. We begin
by recalling that the argument in \cite{ghystsuboi} reduces the proof of Theorem~A to checking that $h$
is a diffeomorphism of class~$C^1$. The remainder of the section will be devoted to proving this statement
under the conditions indicated below.

To make the discussion accurate, let $G_1$ and $G_2$ be two finitely generated subgroups
of $\dif$ that are conjugate by a homeomorphism $h : S^1 \rightarrow S^1$. By assumption, the group $G_1$
is locally $C^2$-non-discrete. In view of the material presented in the previous sections, the following
conditions can be assumed to hold without loss of generality.
\begin{enumerate}

  \item All the orbits of $G_1$ are dense in $S^1$ (in particular $G_1$ has no finite orbit). The same
  condition applies to $G_2$ since these groups are topologically conjugate.

  \item There is an interval $I = [-a,a] \subset \R \subset S^1$ ($a\neq 0$) and an element $F_1$ in
  $G_1$ satisfying $F_1(0) =0$ and $F_1' (0) =\lambda_1 \in (0,1)$.

  \item For every $\varepsilon >0$, we can find a finite set
  $\{ \overline{g}_{1,1} , \ldots , \overline{g}_{1,N} \} \subset G_1$ satisfying all the conditions below:
  \begin{itemize}
    \item $\overline{g}_{1,1} , \ldots , \overline{g}_{1,N}$ are $\varepsilon$-close to the identity in the $C^2$-topology
    on $I$ (where $I = [-a,a]$ is the above chosen interval).

    \item $\overline{g}_{1,1} , \ldots , \overline{g}_{1,N}$ generate a non-solvable subgroup of $\dif$ having no finite orbit.

    \item Consider the sequence $\widetilde{S}_1 (k)$ defined in Section~4 by means of the set
    $\widetilde{S}_1(0) = S_1(0) =S_1 = \{ \overline{g}_{1,1} , \ldots , \overline{g}_{1,N} \}$ so that
    $\widetilde{S}_1(k) = F_1^{-kn} \circ S_1(k) \circ F_1^{kn}$ for every $k \in \N$ and a certain fixed $n \in \N^{\ast}$.
    Then every sequence of elements $\{ g_{1,k} \}$ with $g_{1,k} \in \widetilde{S} (k)$ converges to the identity in the
    $C^2$-topology on the interval $I$.
  \end{itemize}

  \item In fact, if $\{ g_{1,k} \} \subset G_1$ is such that $g_{1,k} \in \widetilde{S}_1 (k)$, $k \in \N$, then
  for every $\Theta \in \R^{\ast}_+$, we have
  \begin{equation}
  \lim_{k \rightarrow \infty} \left[ \frac{\Vert g_{1,k} - {\rm id} \Vert_{2, [-a,a]}}{\Theta^k} \right] = 0 \, . \label{vanishinglimit}
  \end{equation}

\end{enumerate}

Next recall that a point $p \in S^1$ is said to be {\it expandable}\, for a given group $G \subset \dif$ if there is
$g \in G$ such that $g' (p) >1$. Since our diffeomorphisms preserve the orientation of $S^1$, the reader will note
that the conditions $g' (p) >1$ and $\vert g'(p) \vert >1$ are equivalent. With this terminology, we state:

\begin{teo}
\label{ConjugatignEpxanding}
Assume that $G_1$ satisfy all the conditions 1--4 above. Assume also that every point
$p \in S^1$ is expandable for $G_2$. Then every homeomorphism $h: S^1 \rightarrow S^1$ conjugating $G_1$ to $G_2$
coincides with an element of $\dif$.
\end{teo}

The following simple lemma clarifies the connection between Theorem~A and Theorem~\ref{ConjugatignEpxanding}.

\begin{lemma}
\label{FinishingJob-11.1}
Assume that $G \subset \dif$ is a locally $C^2$-non-discrete group satisfying conditions 1--4 above.
Then $G$ leaves no probability measure on $S^1$ invariant. Moreover, every point
$p \in S^1$ is expandable for $G$.
\end{lemma}

\noindent {\it Proof}. Since $G$ has all orbits dense, every probability measure invariant by $G_1$ must
be supported on all of $S^1$. As previously seen, up to parameterizing $S^1$ by the corresponding Radon-Nikodym
derivative, the group $G_1$ becomes conjugate to a group of rotations. This is impossible since $G$ contains
elements exhibiting hyperbolic fixed points.

To establish the second part of the statement, we proceed as follows. Consider first the case of a point $p$ in the
interval $I = [-a,a]$ where $G$ contains a diffeomorphism $F$ satisfying $F (0) =0$ and $F' (0) \in (0,1)$. Owing
to the discussion in Section~2.1, $I$ is equipped with a nowhere zero vector field $X$ contained in the $C^1$-closure of
$G$. Choose then $t_0$ so that the local flow $\phi^t$ of $X$ satisfies $\phi^{t_0} (p) =0$. Therefore the
diffeomorphism $\overline{f} = \phi^{-t_0} \circ F \circ \phi^{t_0}$ satisfies $\overline{f} (p) =p$ and
$\overline{f} (p) >1$. Since $X$ lies in the $C^1$-closure of $G$, there follows that
$\phi^{t_0}$ is the $C^1$-limit of a sequence $\tilde{f}_r$ of elements in $G$ restricted to some small
neighborhood of $p$. Thus, for $r$ large enough, we conclude that $(\tilde{f}_r^{-1} \circ F \circ \tilde{f}_r)' (p)
> 1$ proving the statement for points in $I$. To finish the proof of the lemma, just note that the minimal character
of $G$ enables us to find a finite covering of $S^1$ by intervals satisfying the same conditions used above for
the interval~$I$. The lemma then follows.\qed

\vspace{0.1cm}

\noindent {\it Proof of Theorem~A}. It follows at once from the combination of Theorem~\ref{ConjugatignEpxanding}
with Lemma~\ref{FinishingJob-11.1}.\qed

\vspace{0.1cm}

The rest of this section, will be devoted to the proof of Theorem~\ref{ConjugatignEpxanding}.
To begin with, let us state Proposition~\ref{FinishingJob-11.3}.
For this, first note that diffeomorphisms in $G_1$ having a hyperbolic fixed point in $I$ are far from unique. We have
fixed one of them, namely $F_1$. The element $F_2$ of $G_2$ verifying $F_2 = h^{-1} \circ F_1 \circ h$ has therefore
a fixed point in the interval $J = h^{-1} (I)$, namely the point $q = h^{-1} (0)$. However, since $h$ is only
a homeomorphism, we cannot immediately conclude that $q$ is hyperbolic for $F_2$. In fact, whereas
$F_2$ certainly realizes a ``topological contraction'' on a neighborhood of $q$, it may happen that $F_2' (q) =1$.
The possibility of having $F_2' (q) =1$ is a bit of an inconvenient since it would require us to work with iterations of
a ``parabolic map'' in a context similar to the one discussed in Section~2.1. This type of difficulty, however, can
be overcome with the help of Proposition~\ref{FinishingJob-11.3} (below). The proof of this proposition will be deferred
to the next section and it relies heavily on the methods of \cite{dkn}.

\begin{prop}
\label{FinishingJob-11.3}
Without loss of generality, we can assume that $F_2' (q) <1$ where $q = h^{-1} (0)$.
\end{prop}

Now consider a $C^1$-diffeomorphism
$f : S^1 \rightarrow S^1$. Given an interval $U \subset \R \subset S^1$, the {\it distortion of $f$ in $U$}\, is defined
as
\begin{equation}
\varpi \, (f,U) = \sup_{x,y \in U} \, \log \frac{\vert f'(x) \vert}{\vert f'(y) \vert} = \sup_{x \in U} \log
(\vert  f'(x) \vert) - \inf_{y \in U} \log (\vert f'(y) \vert) \label{definitiondistortion}
\end{equation}
where $\vert \, . \, \vert$ stands for the absolute value. Furthermore, assuming that the map $x \mapsto \log (\vert D_x f \vert)$
has a Lipschitz constant $C_{\rm Lip}$, the estimate
\begin{equation}
\varpi \, (f, U) \leq C_{\rm Lip} \,  \Length (U) \label{Lipschitzestimate}
\end{equation}
holds (where $\Length (U)$ stands for the length of the interval~$U$ with respect to the Euclidean metric
for which the length of $S^1$ is equal to~$2\pi$).
Note also that the mentioned Lipschitz condition is satisfied
provided that $f$ is of class $C^2$ on $U$. On the other, given two diffeomorphisms $f_1, f_2 : S^1 \rightarrow S^1$ as above,
the estimate
\begin{equation}
\varpi \, (f_1 \circ f_2 , U) \leq \varpi \, (f_1 , f_2 (U)) + \varpi \, (f_2, U) \label{subadditivitydistortion}
\end{equation}
holds provided that both sides are well defined.

Let us now go back to the sequence of sets $\{ \widetilde{S}_1 (k) \} \subset G_1$ fixed in the beginning of the section.
Every element in $\widetilde{S}_1 (k)$ is defined on the interval $I = [-a,a]$ and this holds for every $k \in \N$.
Next recall that this sequence was obtained as indicated
in Section~4 by means of the finite set $\{ \overline{g}_{1,1} , \ldots , \overline{g}_{1,N} \} \subset G_1$ and of
the diffeomorphism $F_1$. In particular $\widetilde{S}_1(k) = F_1^{-kn} \circ S_1(k) \circ F_1^{kn}$. From now on,
we fix a sequence $\{ g_{1,k} \} \subset G_1$ of diffeomorphisms such that $g_{1,k} \neq {\rm id}$ belongs to $\widetilde{S}_1(k)$
for every $k \in \N$. Consider also
the corresponding sequence $\{ \widetilde{S}_2 (k) \} \subset G_2$ defined by means of
$\{ \overline{g}_{2,1} , \ldots , \overline{g}_{2,N} \} \subset G_2$ and of
the diffeomorphism $F_2$. More precisely, we set $\overline{g}_{2,j} = h^{-1} \circ \overline{g}_{1,j} \circ h$
for every $j=1, \ldots, N$ and $F_2 = h^{-1} \circ F_1\circ h$ where $F_2$ is assumed to have
a contractive hyperbolic fixed point at $j=h^{-1} (0)$ (cf. Proposition~\ref{FinishingJob-11.3}).
Thus, for every $k \in \N$, we have
$g_{2,k} = h^{-1} \circ g_{1,k} \circ h$. Finally we also pose $J = h^{-1} (I)$.

Next, for every $k \in \N$, let $\mathcal{P}_{I,k}$ denote the partition of the interval $I$ into $5^k$ sub-intervals having
the same size and let $\mathcal{P}_{I,k} = \{ I_{1,k} , \ldots , I_{5^k,k} \}$. By means of the homeomorphism $h$,
the partitions $\mathcal{P}_{I,k}$ induce partitions $\mathcal{P}_{J,k}  = \{ J_{1,k} , \ldots , J_{5^k,k} \}$ of the interval $J$
where $J_{j,k} = h^{-1} (I_{j,k})$ for every $k \in \N$ and for every $j\in \{1, \ldots ,5^k\}$.
Now we have:

\begin{lemma}
\label{FinishingJob-22.1}
Denote by $\varpi \, (g_{2,k}, J_{l_k,k})$ the distortion of $g_{2,k}$ in the interval $J_{l_k,k}$. Then
to each $k \in \N$ there corresponds $l_k \in \{ 1, \ldots ,5^k \}$ such that
the resulting sequence $k \mapsto \varpi \, (g_{2,k}, J_{l_k,k})$ converges to zero.
\end{lemma}

\noindent {\it Proof}. Consider the set formed by the diffeomorphisms $\overline{g}_{2,1} , \ldots , \overline{g}_{2,N}$, $F_2$
along with their inverses. The semigroup generated by this set of diffeomorphisms coincides
with the group generated by $\overline{g}_{2,1} , \ldots , \overline{g}_{2,N}$, and $F_2$
(the set formed by $\overline{g}_{2,1} , \ldots , \overline{g}_{2,N}$,
$F_2$ and their inverses is symmetric in the sense that whenever a diffeomorphism belongs to it so does the
inverse of the diffeomorphism in question). Every element in the group in question
can be represented as a word in the
alphabet whose letters are the diffeomorphisms in the initial (symmetric) set.
If $f$ represents an element in this alphabet, i.e. a letter,
the map $x \mapsto \log (\vert D_x f \vert)$ is well defined on
all of $S^1$ (since $D_x f \neq 0$ for all $x \in S^1$). These maps are also Lipschitz on all of $S^1$ since $f$
is, in any event, a $C^2$-diffeomorphism. Fix then a positive constant $C$ greater than the maximum among the Lipschitz
constants of all the maps $x \mapsto \log (\vert D_x f \vert)$ with $f$ belonging to the alphabet in question.

The explicit construction of the sequences $\{ g_{1,k} \}$ and $\{ g_{2,k} \}$ makes it clear that every diffeomorphism
$g_{2,k}$ belongs to the semigroup generated by $\overline{g}_{2,1} , \ldots , \overline{g}_{2,N}$, $F_2$ and their
inverses. Moreover, as element in this semigroup, the diffeomorphism $g_{2,k}$ can be spelled out in at most
$4^k +2nk$ letters. Next let $c_1$ be a constant such that $c_1 \Length (J) > 2\pi$ (where $\Length (J)$ stands
for the length of $J$). Note also that every diffeomorphism $f$ of the circle must satisfy $\Length (f (J)) < 2\pi$.

Now fixed $k \in \N$, let $g_{2,k} = f_l \circ \cdots \circ f_1$ denote the above mentioned
spelling of $g_{2,k}$. Thus
$l \leq 4^k + 2nk$ and $f_i$ belongs to $\{ F_2^{\pm 1} , \overline{g}_{2,1}^{\pm 1} , \ldots , \overline{g}_{2,N}^{\pm 1} \}$
for every $i \in \{ 1, \ldots , l \}$. The subadditivity relation expressed by~(\ref{subadditivitydistortion})
combined to estimate~(\ref{Lipschitzestimate}) yields
$$
\varpi \, (g_{2,k}, J) \leq C \sum_{i=1}^l \Length ( f_i \circ \cdots \circ f_1 (J) ) + C \, \Length (J)
\leq c_1 C \, \Length (J)  (4^k +2nk) \, .
$$
On the other hand, given an sub-interval $J_{j,k}$ in the partition $\mathcal{P}_{J,k}$
(so that $j \in \{ 1 , \ldots , 5^k \}$). The preceding argument ensures that
$$
\varpi \, (g_{2,k}, J_{j,k}) \leq C \sum_{i=1}^l \Length (f_i \circ \cdots \circ f_1 (J_{j,k}) ) + C \, \Length (J_{j,k})
\leq c_1 C \, \Length (J_{j,k} ) (4^k +2nk) \, .
$$
However, we clearly have
$$
\sum_{j=1}^{5^k} \left[ \sum_{i=1}^l \Length (f_i \circ \cdots \circ f_1 (J_{j,k}) ) \right] + \sum_{j=1}^{5^k}
\Length (J_{j,k} ) = \sum_{i=1}^l \Length (f_i \circ \cdots \circ f_1 (J))  + \Length (J ) \, .
$$
Hence there follows that
$$
\sum_{j=1}^{5^k} \varpi \, (g_{2,k}, J_{j,k}) \leq c_1 C \, \Length (J ) (4^k +2nk) \, .
$$
Finally, if $j(k)$ realizes the minimum of $j \mapsto \varpi \, (g_{2,k}, J_{j,k})$ over the set $\{ 1, \ldots ,5^k \}$,
we conclude that
$$
\varpi \, (g_{2,k}, J_{j(k),k}) \leq \frac{c_1 C \, \Length (J ) (4^k +2nk)}{5^k}
$$
which goes to zero as $k \rightarrow \infty$. The proof of the lemma is completed.\qed

As $k$ increases, we know that $g_{2,k} (y) -y$ converges uniformly to zero on all of $J$. However, when we consider the
sequence of sub-intervals $J_{j(k),k}$ their diameters go to zero as well. A comparison between
$\sup_{y \in J_{j(k),k}} \vert g_{2,k} (y) -y \vert$ and the length $\Length (J_{j(k),k} )$ of $J_{j(k),k}$
will however be needed. In particular,
we would like to claim that the sequence of quotients $\sup_{y \in J_{j(k),k}} \vert g_{2,k} (y) -y \vert /
\Length ( J_{j(k),k} )$ converges to zero as $k \rightarrow \infty$. At this moment, our results are not sufficient
to derive this conclusion since we have no control of the ratio between the lengths of two interval
$J_{j_1(k),k}$ and $J_{j_2(k),k}$. The desired convergence, however, will follow at once from
Proposition~\ref{hisHoldercontinuous} below.

The next step consists of {\it magnifying}\, the intervals $I_{j(k),k}$ and $J_{j(k),k}$ into intervals with diameters
bounded from below by a strictly positive constant. To do this, we shall resort to a slightly more straightforward version
of the celebrated ``Sullivan's expansion strategy'' as expounded in \cite{navas} and \cite{shubsullivan}. The main
difficulty in applying Sullivan's type of argument to our situation lies in the fact that the procedure needs to
be simultaneously applied to both groups $G_1$ and $G_2$. To achieve this, we shall first establish that the
conjugating homeomorphism $h :S^1 \rightarrow S^1$ is H\"older continuous for a suitable exponent
(Proposition~\ref{hisHoldercontinuous} below). Proposition~\ref{hisHoldercontinuous} will subsequently be combined
with the several estimates involving convergence rates for the sequences $\{ g_{1,k} \}_{k \in \N}$ and
$\{ g_{2,k} \}_{k \in \N}$ (restricted to the intervals $I_{j(k),k}$ and $J_{j(k),k}$, respectively) to yield
Theorem~A.

Recall that a map $f : U \subset S^1 \rightarrow S^1$ is said to be {\it $\alpha$-H\"older continuous}\, on the
interval $U$ if the supremum
\begin{equation}
\sup_{x,y \in U \\ x\neq y} \frac{\vert f(x) - f(y) \vert}{\vert x-y \vert^{\alpha}}
\end{equation}
is finite (where the bars $\vert \, . \, \vert$ stand for the fixed Euclidean metric on $S^1$). The above definition
is local in the sense that the length of $U$ and of $f(U)$ are assumed to be smaller than $\pi$ so that the above indicated
distances are well defined. We shall say that $f$ is $\alpha$-H\"older continuous if its restriction to every interval
$U$ satisfying $\max \{ \Length (U) , \; \Length (f(U)) \} < \pi$ is $\alpha$-H\"older continuous on $U$.
With this terminology, we have:

\begin{prop}
\label{hisHoldercontinuous}
There is $\alpha >0$ such that the homeomorphisms $h : S^1 \rightarrow S^1$ and $h^{-1} : S^1 \rightarrow S^1$
are both $\alpha$-H\"older continuous.
\end{prop}

We begin the proof of Proposition~\ref{hisHoldercontinuous} by recalling that
to each point $p \in S^1$ there corresponds a diffeomorphism $f_{1,p} \in G_1$ with $f_{1,p} (p) >1$ (Lemma~\ref{FinishingJob-11.1}).
Owing to the compactness of $S^1$, there is a finite covering
$\mathfrak{U}_1 = \{ U_{1,1}, \ldots , U_{1,s} \}$ of $S^1$ by open connected intervals $U_{1,i}$, $i=1,\ldots ,s$,
satisfying the following conditions:
\begin{itemize}
  \item For $2 \leq i \leq s-1$, the interval $U_{1,i}$ intersects only the intervals $U_{1,i-1}$ and $U_{1,i+1}$. The interval
  $U_{1,1}$ (resp. $U_{1,s}$) intersects only the intervals $U_{1,2}$ and $U_{1,s}$ (resp. $U_{1,s-1}$ and $U_{1,1}$).

  \item To each interval $U_{1,i}$ there corresponds a diffeomorphism $f_{1,i} \in G_1$ such that $f_{1,i}' (x) >1$ for
  every $x \in U_{1,i}$ (recall that $G_1$ and $G_2$ preserve the orientation of $S^1$).
\end{itemize}
Let $m_1 >1$ be given as
$$
m_1 = \min_{i \in \{1, \ldots ,s\}}  \{  \, \inf_{U_{1,i}}  f_{1,i}'  \, \} \, .
$$
Similarly let $M_1 = \max_{i \in \{1, \ldots ,s\}}  \{ \, \sup_{U_{1,i}} f_{1,i}' \, \}$. Next let $L > 0$ denote
the minimum of the lengths of the sets $U_{1,1} \cap U_{1,s}$ and $U_{1,i} \cap U_{1,i+1}$ (for $i=1, \ldots ,s-1$) so that
every interval $[a,b] \subset S^1$ of length less than $L$ is contained in some interval $U_{1,i_1}$
($i_1 \in \{ 1, \ldots ,s \}$). For $[a,b]$ as indicated,
the derivative of $f_{1,i_1}$ is not less than~$m_1 > 1$ at every point in $[a,b]$
and the length $\Length (f_{1,i_1} ([a,b]))$ of $f_{1,i_1} ([a,b])$ is at least $m_1 \, \Length ([a,b]) >
\Length ([a,b])$. When $\Length (f_{1,i_1} ([a,b]))$ is still less than $L$, $f_{1,i_1} ([a,b])$
is again contained in some interval $U_{1,i_2}$. Thus $f_{1,i_2} (f_{1,i_1} ([a,b]))$ has length greater than or equal to
$m_1^2 \, \Length ([a,b])$ and the procedure can be continued provided that
$\Length (f_{1,i_2} \circ f_{1,i_1} ([a,b])) < L$. Thus we have proved the following:

\begin{lemma}
\label{addedinMarch15-1.1}
To every interval $[a,b] \subset S^1$ whose length ($b-a$) is less than $L$, we can assign an element $F_{[a,b]} \in G_1$
satisfying the following conditions:
\begin{enumerate}
  \item $F_{1,[a,b]} = f_{1,i_r} \circ \cdots \circ f_{1,i_1}$ where each $i_l$ belongs to $\{ 1, \ldots ,s\}$.

  \item For every $l \in \{ 1, \ldots ,r\}$, $f_{1,i_{l-1}} \circ \cdots \circ f_{1,i_1} ([a,b])$ is contained
  in $U_{1,i_l}$ (where $f_{1,i_{l-1}} \circ \cdots \circ f_{1,i_1} ([a,b]) =[a,b]$ if $l=1$).

  \item We have
  $$
  L \leq \Length ( F_{1,[a,b]} ([a,b])) \leq L \, M_1 \; .
  $$
  \mbox{}\qed
\end{enumerate}
\end{lemma}

Recalling that $\mathfrak{U}_1 = \{ U_{1,1}, \ldots , U_{1,s} \}$, we define a new covering $\mathfrak{U}_2$ of
$S^1$ by letting $\mathfrak{U}_2 = \{ U_{2,1}, \ldots , U_{2,s} \}$ where $U_{2,i} = h^{-1} (U_{1,i})$ for
every $i=1, \ldots ,s$. To every diffeomorphism $F_{1,[a,b]} = f_{1,i_r} \circ \cdots \circ f_{1,i_1} \in G_1$
as above, we also assign the corresponding diffeomorphism
$$
F_{2,[a,b]} = f_{2,i_r} \circ \cdots \circ f_{2,i_1} =  h^{-1} \circ F_{1,[a,b]} \circ h
$$
where $f_{2,i_l} = h^{-1} \circ f_{1,i_l} \circ h$ for every $l \in \{ 1, \ldots ,r\}$. Clearly
the diffeomorphism $F_{2,[a,b]}$ takes the (small) interval $h^{-1} ([a,b])$ to the interval
$h^{-1} (F_{1,[a,b]} ([a,b]))$ whose diameter is bounded from below by a positive constant since $h$
is uniformly continuous ($S^1$ is compact). Moreover we can still define
$$
M_2 = \max_{i \in \{1, \ldots ,s\}}  \{ \, \sup_{U_{2,i}} f_{2,i}' \, \}
$$
so that $M_2 >1$. However at this point we cannot ensure that $\inf_{U_{2,i}}  f_{2,i}' >1$ for a given
$i\in \{1, \ldots ,s\}$.

Before proving Proposition~\ref{hisHoldercontinuous}, we note that the roles of $G_1$ and $G_2$ in the above
construction can be reversed since, by assumption, every point in $S^1$ is expandable for $G_2$ as well
(see the Theorem~\ref{ConjugatignEpxanding}).

\vspace{0.1cm}

\noindent {\it Proof of Proposition~\ref{hisHoldercontinuous}}. By using the above introduced coverings
$\mathfrak{U}_1$ and $\mathfrak{U}_2$ of $S^1$, we are going to show the existence of $\alpha >0$ so that
$h$ is $\alpha$-H\"older continuous. By reversing the roles of $G_1$ and $G_2$ as indicated above, the same argument
will also imply the $\alpha$-H\"older continuity of $h^{-1}$ as well (up to reducing $\alpha >0$).

To prove that $h$ is $\alpha$-H\"older continuous, we first observe that the statement has a local character. More
precisely, considering points $c\neq d$ in $S^1$, we need to find constants $C \in \R_+^{\ast}$ and
$\alpha >0$ such that
$$
\vert h(d) -h(c) \vert \leq C \, \vert d-c \vert^{\alpha}
$$
provided that $\vert d -c \vert$ is uniformly small at some level. Here the vertical bars $\vert . \vert$
stand for the distance between the corresponding points for the fixed Euclidean metric (i.e.
$\vert d -c \vert = \Length ([c,d])$). Owing to the previous discussion and to the fact that both $h$ and
$h^{-1}$ are uniformly continuous since $S^1$ is compact, there easily follows the existence of a
uniform $\tau >0$ so that all the considerations below are well defined provided that $\vert d -c\vert < \tau$.
By means of these considerations, we shall then establish
that $h$ is $\alpha$-H\"older continuous on intervals whose length does not exceed $\tau$ and this suffices for
the proposition.

We therefore consider $c,d$ as before and let $[a,b] = h([c,d])$. Without loss of generality, $h$ preserves
the orientation so that we set $a = h(c)$ and $b = h(d)$. The next step consists of expanding the interval $[a,b]$
by means of the procedure summarized by Lemma~\ref{addedinMarch15-1.1}. With the notation used in this lemma,
we find $F_{1,[a,b]} = f_{1,i_r} \circ \cdots \circ f_{1,i_1} \in G_1$ such that
\begin{equation}
L \leq \Length ( F_{1,[a,b]} ([a,b])) \leq M_1 L \, . \label{encadrementdelongueur}
\end{equation}
Consider now the corresponding element $F_{2,[a,b]} = h^{-1} \circ F_{1,[a,b]} \circ h$ in $G_2$.
We also set $F_{2,[a,b]} = f_{2,i_r} \circ \cdots \circ f_{2,i_1}$ as previously indicated. There exists a uniform
$\delta >0$ so that
$$
\Length (F_{2,[a,b]} ([c,d])) \geq \delta >0 \, .
$$
Indeed, just note that $F_{2,[a,b]} ([c,d]) = h^{-1} \circ F_{1,[a,b]} ([a,b])$ so that the claim
follows from the uniform continuity of $h^{-1}$ since $\Length ( F_{1,[a,b]} ([a,b])) \geq L >0$.

We now consider the number $r$ of diffeomorphisms $f_{1,i}$ ($i \in \{1 , \ldots ,s\}$) appearing in the above
indicated spelling of $F_{1,[a,b]}$. By construction, at each iteration of $f_{1,i}$ the corresponding interval
is expanded by a factor bounded from below by $m_1 >1$. Hence we obtain
\begin{equation}
\frac{\vert b-a \vert \, m_1^r}{L M_1} \leq 1 \, . \label{firstestimateaboutR}
\end{equation}
On the other hand, considering $F_{2,[a,b]} = f_{2,i_r} \circ \cdots \circ f_{2,i_1}$, there also follows that
at each iteration of $f_{2,i}$ the interval in question cannot be expanded by a factor exceeding $M_2$. Hence,
we similarly obtain $\vert d-c \vert M_2^r \geq \delta$ so that
\begin{equation}
1 \leq \frac{\vert d-c \vert \, M_2^r}{\delta}  \, . \label{secondestimateaboutR}
\end{equation}
Putting together Estimates~(\ref{firstestimateaboutR}) and~(\ref{secondestimateaboutR}), we conclude that
$$
\frac{\vert d-c \vert \, M_2^r}{\delta} \geq \frac{\vert b-a \vert \, m_1^r}{L M_1} \, .
$$
Thus
\begin{equation}
\vert b - a \vert \leq \frac{LM_1}{\delta} \, \vert d-c \vert \, \left( \frac{M_2}{m_1} \right)^r \, . \label{justaboveHolder1}
\end{equation}
Without loss of generality, we can assume $M_2 > m_1$ for otherwise the preceding estimate implies at
once that $h$ is Lipschitz. Note, however, that Estimate~(\ref{firstestimateaboutR}) also yields
\begin{equation}
r \leq \frac{1}{\ln m_1} \left( \ln (LM_1) - \ln \vert b-a \vert \right) \, . \label{logboundonr}
\end{equation}
Therefore Estimate~(\ref{justaboveHolder1}) becomes
$$
\vert b - a \vert \leq \frac{LM_1}{\delta} \left( \frac{M_2}{m_1} \right)^{\ln (LM_1) /\ln m_1} \,
\vert d-c \vert \left( \frac{M_2}{m_1} \right)^{-ln \vert b-a \vert /\ln m_1} \, .
$$
Let now
$$
C_1 = \frac{LM_1}{\delta} \left( \frac{M_2}{m_1} \right)^{\ln (LM_1) /\ln m_1} \; \; \; \; {\rm and} \; \; \; \;
\overline{c} = - \ln \left[ \left( \frac{M_2}{m_1} \right)^{-1/\ln m_1} \right] \, .
$$
Note that $\overline{c} >0$ since $M_2 \geq m_1$ and $\ln m_1 >0$ since $m_1 > 1$. Hence there follows
\begin{eqnarray*}
\vert b-a \vert & \leq & C_1 \vert d-c \vert \exp (\ln \vert b-a \vert^{-\overline{c}}) \\
& \leq & C_1 \vert d-c \vert \vert b-a \vert^{-\overline{c}} \, .
\end{eqnarray*}
The proposition now results by choosing $\alpha =1 /(1+ \overline{c})$ and $C= C_1^{\alpha}$.\qed

We are almost ready to prove Theorem~\ref{ConjugatignEpxanding}. The last ingredient needed in
our proof consists of a simple estimate for the second derivatives of $F_{1,[a,b]}$ and $F_{2,[a,b]}$.
This is as follows. Keep the preceding notation and fix again intervals $[a,b]$ and $[c,d]$ such that
$h ([c,d]) = [a,b]$. Then we have:

\begin{lemma}
\label{estimatesecondderivativeExpansion}
There are constants $\overline{C} >0$ and $\beta \in (0,1)$ such that
$$
\max \left\{ \sup_{x \in [a,b]} \vert D^2 F_{1,[a,b]} (x) \vert \; \; ; \; \;
\sup_{y \in [c,d]} \vert D^2 F_{2,[a,b]} (y) \vert \right\} \leq \overline{C} \vert b-a \vert^{\ln \beta} \; .
$$
\end{lemma}

\noindent {\it Proof}. Let $\overline{M}$ be a constant satisfying
$$
\max_{i=1, \ldots ,s} \left\{  \sup_{U_{1,i}} \vert D^2 f_{1,i} \vert \; \; ; \; \;
\sup_{U_{2,i}} \vert D^2 f_{2,i} \vert \right\} < \overline{M} \, .
$$
First we will show the existence of constants $\overline{C}_1$ and $\beta_1$ for which
$\sup_{x \in [a,b]} \vert D^2 F_{1,[a,b]} (x) \vert \leq \overline{C}_1 \vert b-a \vert^{\ln \beta_1}$.
We begin by recalling that $F_{1,[a,b]} = f_{1,i_r} \circ \cdots \circ f_{1,i_1}$. For $x_0 \in [a,b]$
and $l \in \{ 1, \ldots , r-1\}$, let $x_l = f_{1,i_l} \circ \cdots \circ f_{1,i_1} (x_0)$. Thus we have
$F_{1,[a,b]}'(x_0) = f_{1,i_r}' (x_{r-1}) \cdots  f_{1,i_1}' (x_0)$ and
$$
D^2 F_{1,[a,b]}(x_0) = \prod_{l=1}^r f_{1,i_l}' (x_{l-1})  \left[
\sum_{j=1}^r \left( \frac{D^2 f_{1,i_j} (x_{j-1})}{f_{1,i_j}' (x_{j-1})}
f_{1,i_{j-1}}' (x_{j-2}) \cdots  f_{1,i_1}' (x_0) \right) \right] \, .
$$
Hence
\begin{equation}
\vert D^2 F_{1,[a,b]}(x_0) \vert \leq \overline{M} M_1^{2r} \, . \label{simplesestimate-1}
\end{equation}
On the other hand, recall that $r \leq (\ln LM_1 - \ln \vert b-a \vert )/\ln m_1$ (Estimate~(\ref{logboundonr})).
Setting $\overline{C_1} = \overline{M} M_1^{2\ln LM_1 /\ln m_1}$, there follows that
$$
\vert D^2 F_{1,[a,b]}(x_0) \vert \leq \overline{C}_1 \vert b-a \vert^{-2\ln M_1/\ln m_1} \, .
$$
Since $M_1 \geq m_1 >1$, the exponent $-2\ln M_1/\ln m_1$ is negative and hence has the form
$\ln \beta_1$ for some $\beta_1 \in (0,1)$. This proves the first assertion. To complete the proof
of the lemma it only remains to show that a similar estimate holds for $\vert D^2 F_{2,[a,b]} \vert$ on
$[c,d]$. However, a repetition word-by-word of the above argument yields constants $\overline{C}_2$
and $\beta_2 \in (0,1)$ such that
$$
\vert D^2 F_{2,[a,b]}(y_0) \vert \leq \overline{C}_2 \vert d-c \vert^{\ln \beta_2} \,
$$
for every $y_0 \in [c,d]$. The desired estimate is then an immediate consequence of Proposition~\ref{hisHoldercontinuous}.
The lemma is proved.\qed

\vspace{0.2cm}

\noindent {\it Proof of Theorem~\ref{ConjugatignEpxanding}}. In what follows we keep all the notation introduced in
the course of this section. Consider the interval $I \subset S^1$ (resp. $J=h^{-1} (I) \subset S^1$) and the
sequence of partitions $\mathcal{P}_{I,k}$ on $I$ (resp. $\mathcal{P}_{J,k}$ on $J$). More precisely, consider
the sequences of intervals $k \mapsto I_{l_k,k}$ and $k \mapsto J_{l_k,k}$ where $J_{l_k,k}$ is as in
Lemma~\ref{FinishingJob-22.1}.

Next set $I_{l_k,k} = [a_k,b_k]$ and $J_{l_k,k} = [c_k,d_k]$ so that $a_k = h(c_k)$ and $b_k = h(d_k)$. Also
$\alpha >0$ is fixed so that both homeomorphisms $h$ and $h^{-1}$ are $\alpha$-H\"older continuous
(Proposition~\ref{hisHoldercontinuous}). Now, for each $k \in \N$ fixed, let $F_{1,[a_k,b_k]}$ be the
element of $G_1$ obtained by means of Lemma~\ref{addedinMarch15-1.1}. Thus
we have $F_{1,[a_k,b_k]} = f_{1,i_{r_k}} \circ \cdots \circ f_{1,i_1}$ where each $i_l$ belongs to $\{ 1, \ldots ,s\}$.
Analogously we define $F_{2,[a_k,b_k]} \in G_2$ so that $F_{2,[a_k,b_k]} = h^{-1} \circ F_{1,[a_k,b_k]} \circ h$.
In particular
$$
F_{2,[a_k,b_k]} = f_{2,i_{r_k}} \circ \cdots \circ f_{2,i_1}
$$
with $i_l \in \{ 1, \ldots ,s\}$.

By construction, all the intervals of the form $\{ F_{1,[a_k,b_k]} ([a_k,b_k]) \} \subset S^1$ have length
comprised between $L$ and $L M_1$. Hence, up to passing to a subsequence, we can assume that these intervals
converge towards an open interval $\widetilde{I} =(\widetilde{a}, \widetilde{b}) \subset S^1$ with
$\widetilde{a} \neq \widetilde{b}$. More precisely, this convergence of intervals simply means that
$F_{1,[a_k,b_k]} (a_k) \rightarrow \widetilde{a}$ and $F_{1,[a_k,b_k]} (b_k) \rightarrow \widetilde{b}$.
We also set $\widetilde{J} = h^{-1} (\widetilde{I}) = (\widetilde{c}, \widetilde{d}) \subset S^1$
so that $F_{2,[a_k,b_k]} (c_k) \rightarrow \widetilde{c}$ and $F_{2,[a_k,b_k]} (d_k) \rightarrow \widetilde{d}$.

Consider now the sequences of diffeomorphisms $\{ \widetilde{f}_{1,k} \} \subset G_1$
and $\{ \widetilde{f}_{2,k} \} \subset G_2$ which are defined by
$$
\widetilde{f}_{1,k} = F_{1,[a_k,b_k]} \circ g_{1,k} \circ F_{1,[a_k,b_k]}^{-1} \; \; \;
{\rm and} \; \; \;
\widetilde{f}_{2,k} = F_{2,[a_k,b_k]} \circ g_{2,k} \circ F_{2,[a_k,b_k]}^{-1} \; .
$$

\noindent {\it Claim~1}. The sequence $\{ \widetilde{f}_{1,k} \} \subset G_1$ (resp. $\{ \widetilde{f}_{2,k} \} \subset G_2$)
converges to the identity in the $C^0$-topology on compact parts of $\widetilde{I}$ (resp. $\widetilde{J}$).

\noindent {\it Proof of Claim~1}. Consider first the sequence $\{ \widetilde{f}_{1,k} \}$ and a point $x \in
\widetilde{I}$. By construction the point $y = F_{1,[a_k,b_k]}^{-1} (x)$ lies in $I_{l_k,k} = [a_k,b_k]$
provided that $k$ is large enough. Therefore
\begin{eqnarray*}
\vert \widetilde{f}_{1,k} (x) -x \vert & = & \vert F_{1,[a_k,b_k]} \circ g_{1,k} \circ F_{1,[a_k,b_k]}^{-1} (x) -x \vert \\
& \leq & \sup_{[a_k,b_k]} \vert D^1 F_{1,[a_k,b_k]} \vert \; \vert g_{1,k} (y) -y \vert \\
& \leq & M_1^{r_k} \,  \vert g_{1,k} (y) -y \vert \, .
\end{eqnarray*}
However $r_k$ is bounded by Estimate~(\ref{logboundonr}). In turn, this yields
$M_1^{r_k} \leq {\rm const} \vert b_k - a_k \vert^{-\ln M_1 /\ln m_1}$ for some constant ${\rm const}$. In turn,
up to a multiplicative constant, $\vert b_k -a_k \vert$ equals $5^{-k}$. Thus, for a new suitable constant
${\rm Const}$, we obtain
$$
\vert \widetilde{f}_{1,k} (x) -x \vert \leq {\rm Const} \, 5^{k\ln M_1 /\ln m_1}  \vert g_{1,k} (y) -y \vert
$$
which converges to zero as $k \rightarrow \infty$ by virtue of condition~4 in the beginning of the section.

It remains to show the same holds for the sequence $\{ \widetilde{f}_{2,k} \}$. Setting $z =h^{-1} (x)$ and
$w =h^{-1} (y)$, the same argument used above yields
$$
\vert \widetilde{f}_{2,k} (z) -z \vert \leq {\rm Const}' \, 5^{k\ln M_2 /\ln m_1}  \vert g_{2,k} (w) -w \vert
$$
for a new constant ${\rm Const}'$. However the $\alpha$-H\"older continuity of $h^{-1}$ ensures
that $\vert g_{2,k} (w) -w \vert \leq  \vert g_{1,k} (y) -y \vert^{\alpha}$ so that the claim follows again
from condition~4 in the beginning of the section.\qed

We now consider the problem of $C^1$-convergence for the sequences $\{ \widetilde{f}_{1,k} \}$
and $\{ \widetilde{f}_{2,k} \}$. We begin by recalling that the restriction of $\{ g_{1,k} \}$ to
$I_{l_k,k}$ converges $C^2$ (in particular $C^1$) to the identity. On the other hand, the restriction
of $g_{2,k}$ to $J_{l_k,k}$ is known to satisfy the following conditions:
\begin{itemize}
  \item[(A)]
  $$
  \frac{\sup_{w \in J_{l_k,k}} \vert g_{2,k} (w) -w \vert}{\Length (J_{l_k,k} )} \longrightarrow 0 \, .
  $$

  \item[(B)] The sequence $\{ \varpi \, (g_{2,k}, J_{l_k,k}) \}$ formed by the distortion of $g_{2,k}$ on $J_{l_k,k}$
converges to zero.
\end{itemize}
Tee reader will note that item~(B) is nothing but Lemma~\ref{FinishingJob-22.1}. In turn,
item~(A) follows from the above discussion since the analogous statement holds for $\{ g_{1,k} \}$
(condition~4) and both $h$ and $h^{-1}$ are $\alpha$-H\"older continuous. In fact, the $\alpha$-H\"older continuity
of $h$ ensures that $\Length (J_{l_k,k} ) \geq \Length (I_{l_k,k} )^{1/\alpha}$ while the
$\alpha$-H\"older continuity of $h^{-1}$ yields $\sup_{w \in J_{l_k,k}} \vert g_{2,k} (w) -w \vert \leq
\sup_{x \in I_{l_k,k}} \vert g_{1,k} (x) -x \vert^{\alpha}$ so that condition~4 establishes the desired limit.

Owing to Proposition~\ref{FinishingJob-11.3} and to the fact $G_1$ acts minimally on $S^1$, we choose a point
$p \in \widetilde{I}$ such that the following condition holds: there are conjugate elements $\widetilde{F}_1 \in G_1$ and
$\widetilde{F}_2 \in G_2$ ($\widetilde{F}_2 = h^{-1} \circ \widetilde{F}_1 \circ h$) such that $\widetilde{F}_1$
has a hyperbolic fixed point in $p$ whereas $\widetilde{F}_2$ has a hyperbolic fixed point in $q =h^{-1} (p)$.
For the reasons already explained, we can assume without loss of generality that $\widetilde{f}_{1,k} (p) \neq p$ for
every~$k \in \N$ (which also implies that $\widetilde{f}_{2,k} (q) \neq q$).

The next step consists of estimating the derivative of $\widetilde{f}_{1,k}$ at a point $x \in \widetilde{I}$.
For $y = F_{1,[a_k,b_k]}^{-1} (x)$, we clearly have $\widetilde{f}_{1,k}' (x) = D^1_{g_{1,k} (y)} F_{1,[a_k,b_k]}
g_{1,k}' (y) D^1_{x} F_{1,[a_k,b_k]}^{-1}$. Thus,
\begin{eqnarray*}
\vert \widetilde{f}_{1,k}' (x) \vert & \leq & \vert D^1_{g_{1,k} (y)} F_{1,[a_k,b_k]} -D^1_{y} F_{1,[a_k,b_k]} \vert
\, \vert g_{1,k}' (y) D^1_{x} F_{1,[a_k,b_k]}^{-1} \vert + \vert g_{1,k}' (y) \vert \\
& \leq & \sup_{[a_k,b_k]} \vert D^2 F_{1,[a_k,b_k]} \vert \, \vert g_{1,k} (y) -y \vert \, g_{1,k}' (y) + g_{1,k}' (y) \, .
\end{eqnarray*}
On the other hand, $\vert b_k - a_k \vert$ is bounded by a uniform constant times $5^{-k}$. Thus
Lemma~\ref{estimatesecondderivativeExpansion} yields
$$
\sup_{[a_k,b_k]} \vert D^2 F_{1,[a_k,b_k]} \vert \leq {\rm const}\, 5^{-k \ln \beta} \, .
$$
Therefore condition~4 ensures that
$\sup_{[a_k,b_k]} \vert D^2 F_{1,[a_k,b_k]} \vert \, \vert g_{1,k} (y) -y \vert$ converges to zero as $k$ goes
to infinity.  Since $\{ g_{1,k} \}$ converges $C^1$ to the identity, there follows that the restriction of
$\widetilde{f}_{1,k}$ to every compact part of $\widetilde{I}$ converges $C^1$ to the identity as well.
The claim below shows that a similar phenomenon holds for the sequence $\{ \widetilde{f}_{2,k} \}$ as well.

\noindent {\it Claim~2}. The sequence $\{ \widetilde{f}_{2,k} \}$ converges $C^1$ to the identity on
$\widetilde{J}$.

\noindent {\it Proof of Claim~2}. The argument is more subtle and builds on the previous discussion. Recalling
that $q = h^{-1} (p)$, we set $q_k = F_{2,[a_k,b_k]}^{-1} (q)$. Let also $\lambda_k = g_{2,k}' (q_k)$. We also
immediately note that Lemma~\ref{estimatesecondderivativeExpansion} still yields
$\sup_{[c_k,d_k]} \vert D^2 F_{2,[a_k,b_k]} \vert \leq {\rm const}\, 5^{-k \ln \beta}$ for a suitable constant
${\rm const}$. For $z \in \widetilde{J}$ and $w = F_{2,[a_k,b_k]}^{-1} (z)$, the argument used above now provides
$$
\vert \widetilde{f}_{2,k}' (x) \vert \leq
\sup_{[c_k,d_k]} \vert D^2 F_{2,[a_k,b_k]} \vert \, \vert g_{2,k} (w) -w \vert \, g_{2,k}' (w) + g_{2,k}' (w) \, .
$$
Again $\sup_{[c_k,d_k]} \vert D^2 F_{2,[a_k,b_k]} \vert \, \vert g_{2,k} (w) -w \vert \, g_{2,k}' (w)$ converges
to zero so that $\vert \widetilde{f}_{1,k}' (x) \vert$ becomes arbitrarily close to $g_{2,k}' (w)$. In turn,
owing to Lemma~\ref{FinishingJob-22.1} the derivative $g_{2,k}' (w)$ becomes arbitrarily close to $\lambda_k$.
Finally we can assume that $\lambda_k$ converges to some $\overline{\lambda} \in \R$ for $\lambda_k$ is uniformly
bounded since the lengths of the intervals $\widetilde{f}_{2,k} (\widetilde{J})$ are clearly so. Summarizing what
precedes, the sequence of maps $\{ \widetilde{f}_{2,k}' \}$ converges uniformly on $\widetilde{J}$ to the constant
$\overline{\lambda}$. To conclude that $\overline{\lambda} =1$, just note that the sequence of primitives
$\{ \widetilde{f}_{2,k} \}$ converges uniformly to the identity on $\widetilde{J}$ (Claim~1). This ends the
proof of Claim~2.

To finish the proof of Theorem~\ref{ConjugatignEpxanding} we proceed as follows. We consider again the sequences
of maps $\{ \widetilde{f}_{1,k} \} \subset G_1$ and $\{ \widetilde{f}_{2,k} \} \subset G_2$. By construction,
we have $\widetilde{f}_{2,k} = h^{-1} \circ \widetilde{f}_{1,k} \circ h$ for every~$k \in \N$. Furthermore
$\{ \widetilde{f}_{1,k} \}$ (resp. $\{ \widetilde{f}_{2,k} \}$) converges $C^1$ to the identity
on $\widetilde{I}$ (resp. $\widetilde{J}$). From this point, the standard argument relies on {\it synchronized
vector fields}\, (see \cite{reb1}). This is as follows.

Recall that $\widetilde{f}_{1,k} (p) \neq p$ (resp. $\widetilde{f}_{2,k} (q) \neq q$) for
every~$k \in \N$. Moreover there are conjugate elements $\widetilde{F}_1 \in G_1$ and $\widetilde{F}_2 \in G_2$
which have hyperbolic fixed points in $p$ and $q$, respectively. In suitable local coordinates around
$p \simeq 0$ (resp. $q \simeq 0$), $\widetilde{F}_1$ becomes a homothety $x \mapsto \Lambda_1 x$ (resp. $\widetilde{F}_1$,
$z \mapsto \Lambda_2 z$. Here both $\Lambda_1$ and $\Lambda_2$ belong to $(0,1)$.
Consider the effect of the conjugations $\widetilde{F}_1^{-j} \circ \widetilde{f}_{1,k} \circ
\widetilde{F}_1^{j}$ on $\widetilde{f}_{1,k}$ for $k$ fixed and $j \in \N$. As explained in Section~2.1,
if $j(k)$ is a suitably chosen sequence with $j(k) \rightarrow \infty$, the conjugate diffeomorphisms
$\widetilde{F}_1^{-j(k)} \circ \widetilde{f}_{1,k} \circ \widetilde{F}_1^{j(k)}$ and
$\widetilde{F}_2^{-j(k)} \circ \widetilde{f}_{2,k} \circ \widetilde{F}_2^{j(k)}$ converge in the $C^1$-topology, respectively
on $\widetilde{I}$ and $\widetilde{J}$, to non-trivial translations. Thus, we actually obtain non-zero
constant vector fields $\overline{X}_1$ and $\overline{X}_2$ contained in the $C^1$-closures
of $G_1$ and $G_2$, respectively, and whose flows $\phi_{1}^t$ and $\phi_{2}^t$
satisfy the equation
$$
h \circ \phi_{2}^t (z) =  \phi_{1}^t \circ h (z) \,
$$
whenever both sides are well defined. By fixing $z$  and letting $t$ takes values around $0 \in \R$, we conclude that
$h$ is of class $C^1$ on a neighborhood of $z \in \widetilde{J}$.
The fact that the dynamics of $G_1$ and $G_2$ are minimal
then implies that $h$ is of class~$C^1$ on the entire circle. The proof of Theorem~\ref{ConjugatignEpxanding} is completed.\qed

\section{Ergodic theory and conjugate groups}

This section is devoted to certain problems about topologically conjugate groups acting on~$S^1$ whose solutions involve
probabilistic arguments. Indeed, in what follows, Proposition~\ref{FinishingJob-11.3} will be proved and
Theorem~\ref{PoissonBoundaryApplication} will be stated and proved. Throughout the section,
we fix two subgroups $G_1$ and $G_2$ of $\dif$ satisfying conditions 1--4 in the beginning of Section~5.
In particular $G_1$ is locally $C^2$-non-discrete.

\begin{lemma}
\label{FinishingJob-11.1PRIME}
None of the groups $G_1$ and $G_2$ leaves a probability measure on $S^1$ invariant.
\end{lemma}

\noindent {\it Proof}. The statement holds for $G_1$ thanks to Lemma~\ref{FinishingJob-11.1}. The
conclusion concerning $G_2$ then arises from the fact that these two groups are topologically conjugate.\qed

We also know that the each of the topologically conjugate groups $G_1$ an $G_2$ acts minimally on $S^1$
(i.e. all their orbits are dense). Next recall
that a group $G$ acting on $S^1$ is said to be {\it proximal}\, if every closed interval can be mapped to intervals
of arbitrarily small length by means of elements of $G$. Now we have:

\begin{lemma}
\label{FinishingJob-11.5}
The action of $G_1$ (resp. $G_2$) on $S^1$ is proximal.
\end{lemma}

\noindent {\it Proof}. Since these groups are topologically conjugate, it suffices to prove the statement for $G_1$.
Consider then the diagonal action of $G_1$ on $S^1 \times S^1$ leaving invariant the diagonal $\Delta \subset S^1 \times S^1$.
Note that $G_1$ is proximal provided that this action is minimal on the open set $S^1 \times S^1 \setminus \Delta$.
Since $G_1$ is not abelian, the statement then follows from the discussion in
\cite{raderson}. Indeed, whereas in \cite{raderson} the main
theorem only claims a decomposition of $S^1 \times S^1$ into finitely many invariant sets on which $G_1$ acts
minimally, these sets are reduced to a single one provided that the group $G_1$ has no finite orbits.\qed

Concerning Lemma~\ref{FinishingJob-11.5}, there
is an alternative point of view enabling us to avoid entering in the details left implicit in the argument given in \cite{raderson}.
This stems from an observation due to Ghys in \cite{ghyssurveyCircle} (page 362) which suffices for the purposes of
this section (so that Lemma~\ref{FinishingJob-11.5} can be dispensed with). According to Ghys,
given a minimal group acting on $S^1$,
there is a finite quotient of the circle and a proximal $C^0$-action induced from $G_1$ on this quotient. Furthermore
this induced action on the quotient commutes with the projection and the initial action of $G_1$ on $S^1$. As mentioned
the existence of this finite quotient of $S^1$ where $G_1$ induces a proximal action is enough for the use of the proximal condition
that will be made in this section (further detail will be provided as they become needed).

At this point we can already prove Proposition~\ref{FinishingJob-11.3}.

\vspace{0.1cm}

\noindent {\it Proof of Proposition~\ref{FinishingJob-11.3}}. The proof is actually a by-product
of the proof of Theorem~F in \cite{dkn}. The argument will be
summarized below and the reader is referred to \cite{dkn} for fuller detail. Since $G_1$ acts minimally on $S^1$, the lemma
is reduced to proving the existence of $F_1 \in G_1$ having a hyperbolic fixed point at $p \in S^1$ such that
the corresponding element $F_2 = h^{-1} \circ F_1 \circ h$ in $G_2$ has a hyperbolic fixed point in $q = h^{-1} (p)$.

Consider a finite generating set $L_1$ for $G_1$
containing elements and their inverses (i.e. $L$ generates $G_1$ as semigroup). Denote by $L_2 =
h^{-1} \circ L_1 \circ h$ the corresponding set in $G_2$. The sets $L_1, \, L_2$ can be put in natural correspondence
with a finite set of letters $\Sigma$ and, through this identifications, we equip $\Sigma$ with a probability measure $\mu$
that is symmetric (i.e. gives the same mass to an element and to its inverse) and non-degenerate
(i.e. every element in $\Sigma$ has strictly positive $\mu$-mass). Denote by $\Omega$ the {\it shift space}\, $\Sigma^{\N}$ equipped
with the standard shift map $\sigma : \Omega \rightarrow \Omega$ and with the probability measure $\mathbb{P} (\Sigma) =
\mu^{\N}$. By a small abuse of notation, we shall identify $\mu$ with measures on $L_1$ and $L_2$. Similarly
$\mathbb{P} (\Sigma)$ (resp. $\sigma$) will also be thought of as a measure (resp. shift map) in either $L_1^{\N}$
or $L_2^{\N}$. Finally, we define maps $T_1$ and $T_2$ from $\Omega \times S^1$ to $\Omega \times S^1$
by letting $T_1 (\omega ,x) = (\sigma (\omega) , \tilde{f}_1^1 (x))$ and
$T_2 (\omega ,x) = (\sigma (\omega) , \tilde{f}_1^2 (x))$ where $\tilde{f}_1^1$ (resp. $\tilde{f}_1^2$) is the projection of
$\omega$ in the first copy of $\Sigma$ viewed with the identifications corresponding to $G_1$ (resp. $G_2$).

Next denote by $\nu_1$ (resp. $\nu_2$) the stationary measure of $G_1$ (resp. $G_2$) defined with respect to $\mu$.
In other words, $\nu_1$ (resp. $\nu_2$) is a probability measure on $S^1$ whose value on a Borel set
$\mathcal{B} \subset S^1$ is given by
$$
\nu_1 (\mathcal{B}) = \sum_{g \in G_1} \mu (g) \nu (g^{-1} (\mathcal{B}))
$$
(resp. $\nu_2 (\mathcal{B}) = \sum_{g \in G_2} \mu (g) \nu (g^{-1} (\mathcal{B}))$). These stationary measures
$\nu_1$ and $\nu_2$ are unique after \cite{dkn} complemented by Lemma~\ref{FinishingJob-11.1PRIME}. From
the uniqueness of these stationary measures, there follows that
$h^{\ast} \nu_1 = \nu_2$. According to Furstenberg \cite{furst},
given a continuous function $\psi$, the sequence of random variables
$$
\xi_{1,l} (\omega) = \int_{S^1} \psi d (\tilde{f}_1^1 \cdots \tilde{f}_l^1 (\nu_1))
$$
(resp. $\xi_{2,l} (\omega) = \int_{S^1} \psi d (\tilde{f}_1^2 \cdots \tilde{f}_l^2 (\nu_2))$) is a martingale so that
both limits
$$
\omega (\nu_1) = \lim_{l\rightarrow \infty} \tilde{f}_1^1 \cdots \tilde{f}_l^1 (\nu_1) \; \; \, {\rm and} \; \; \,
\omega (\nu_2) = \lim_{l\rightarrow \infty} \tilde{f}_1^2 \cdots \tilde{f}_l^2 (\nu_2)
$$
exist for a subset of full $\mathbb{P} (\Sigma)$-measure of $\Sigma$.

Now recall that $G_1$ and $G_2$ are proximal (Lemma~\ref{FinishingJob-11.5}). Alternatively, recall that we can work
with a finite quotient of $S^1$ where $G_1$ induces a proximal action. We can then work with this proximal action on
the quotient and then lift back the result to the initial action of $G_1$. We leave the detail of this construction
to the reader while observing that this latter argument is also used in \cite{dkn}, \cite{deroinETDS}.

The importance of the proximal character of $G_1$ (resp. $G_2$) lies in the fact that the resulting measure
$\omega (\nu_1)$ (resp. $\omega (\nu_2)$) becomes a Dirac mass as originally proved in \cite{antonov};
see also \cite{kleptsynnalski} and Proposition~5.2 of \cite{dkn}.
A topological analogue of the last assertion can be obtained as follows.
Define the {\it contraction coefficient}\, $c(g)$ of a diffeomorphism (homeomorphism) $g$ of $S^1$ as the infimum over
$\epsilon >0$ for which there are closed intervals $U$ and $V$ of sizes not greater than
$\epsilon$ and such that $g (\overline{S^1 / U}) = V$. With this definition, the preceding argument on Dirac masses
also implies that the contraction coefficients $c_l^1 (\tilde{f}_l^1 \cdots \tilde{f}_1^1)$ and
$c_l^2 (\tilde{f}_l^2 \cdots \tilde{f}_1^2)$ converge to zero for a set of full $\mathbb{P} (\Sigma)$-measure of $\Sigma$
(Proposition~5.3 of \cite{dkn}).

The rest of the proof consists of repeating word-by-word the argument detailed in Section~4.4 of \cite{dkn} (aimed at the proof
of Theorem~F in the mentioned paper). Indeed, for a generic choice of $\omega \in \Sigma$, there are a sequence
of intervals $U_l^1, \, V_l^1$ (resp. $U_l^2, \, V_l^2$) whose sizes converge to zero and such that
$$
\tilde{f}_l^1 \cdots \tilde{f}_1^1 (U_l^1) \subset V_l^1 \; \; \, {\rm and} \; \; \,
\tilde{f}_l^2 \cdots \tilde{f}_1^2 (U_l^2) \subset V_l^2 \, .
$$
When $U_l^1, \, V_l^1$ are disjoint then the fixed points of $\tilde{f}_l^1 \cdots \tilde{f}_1^1$ are contained
in these intervals (and the analogous conclusion holds for $\tilde{f}_l^2 \cdots \tilde{f}_1^2 \in G_2$ since
$G_1$ and $G_2$ are topologically conjugate). The argument in \cite{dkn} then continues by showing first that
$U_l^1, \, V_l^1$ are often disjoint. In a second moment, the authors use
techniques of Lyapunov exponents to control the contraction rate
so as to conclude that the fixed points are of hyperbolic nature. This ends the proof of the proposition.\qed

The rest of this section is devoted to the proof of Theorem~\ref{PoissonBoundaryApplication}. This theorem
concerns the potential existence of topologically conjugate groups $G_1$ and $G_2$ acting on~$S^1$ with
$G_1$ being locally $C^2$-non-discrete whereas $G_2$ is locally $C^2$-discrete. This discussion will lead to the
proof of Theorem~B in the introduction. We begin by stating
Theorem~\ref{PoissonBoundaryApplication}. For the rest of this section, $\Gamma$ will always denote
an abstract hyperbolic group which is neither finite nor a finite extension of~$\Z$.

\begin{teo}
\label{PoissonBoundaryApplication}
For $\Gamma$ as above, let $\rho_1 : \Gamma \rightarrow \dif$ be a faithful representation of $\Gamma$ in $\dif$ and set $G_1 =
\rho_1 (\Gamma)$. Suppose that $G_2 \subset \dif$ is another subgroup of $\dif$ which is topologically conjugate
to $G_1$. If $G_1$ is locally $C^2$-non-discrete then so is $G_2$
\end{teo}

The proof of Theorem~\ref{PoissonBoundaryApplication} relies on the combination of a few deep results
including Theorem~1.1 of \cite{deroinETDS}, Kaimanovich's theorem in \cite{KaimanovichPoisson} and
Connell-Muchnik construction in \cite{connel}; cf. Proposition~\ref{regularstationary-addedinApril}.
For suitable background on hyperbolic groups and on measure theoretic methods in group theory,
the reader is referred to \cite{KaimanovichPoisson}, \cite{VershikSurvey}, and \cite{GhysHarpe}.

In the sequel we assume by way of contradiction that the statement of
Theorem~\ref{PoissonBoundaryApplication} is false. Thus there are two topologically conjugate
subgroups $G_1=\rho_1 (\Gamma)$ and $G_2$ of $\dif$ such that $G_1$ is locally $C^2$-non-discrete whereas $G_2$
is locally $C^2$-discrete. By post-composing $\rho_1$ with a conjugating homeomorphism $h$, we obtain another faithful representation
$\rho_2 : \Gamma \rightarrow \dif$ satisfying
$$
\rho_2 (\gamma) = h^{-1} \circ \rho_1 (\gamma ) \circ h
$$
for every $\gamma \in \Gamma$ and where $G_2 = \rho_2 (\Gamma)$.
In other words, the representations $\rho_1$ and $\rho_2$ are topologically
conjugated by~$h$. We begin with the following lemma.

\begin{lemma}
\label{C1discrete-addedinApril}
Without loss of generality we can assume that the group $G_2$ is locally $C^1$-discrete.
\end{lemma}

\noindent {\it Proof}. The proof relies on Theorem~\ref{ConjugatignEpxanding}. In fact, according to
this theorem, the only possibility for $G_1$ being locally $C^2$-discrete occurs when $G_2$ has a non-expanding
point. Hence to prove the lemma it suffices to check that a locally $C^1$-non-discrete subgroup
of $\dif$ (having all orbits dense and leaving no probability measure invariant) expands every point
$p \in S^1$.

Consider then a diffeomorphism $F_2 \in G_2$ having a hyperbolic fixed point $q \in S^1$. In local coordinates
around $q \simeq 0$, we then have $F_2 (x) = \lambda x$ for some $\lambda \in (0,1)$. Next, suppose
that $G_2$ is locally $C^1$-non-discrete. By using the minimal character of $G_2$, we then obtain a sequence
$g_{2, j}$ of diffeomorphisms in $G_2$ ($g_{2, j} \neq {\rm id}$ for all $j \in \N$) which converges to the identity
on a small interval $(-\varepsilon , \varepsilon)$ around $q \simeq 0$ (for some $\varepsilon >0$). Again the discussion
in Section~2.1 allows us to assume that $g_{2, j} (0) \neq 0$ for every $j \in \N$. Thus, as shown in the end
of the proof of Theorem~\ref{ConjugatignEpxanding}, there is a sequence of positive integers $m(j) \rightarrow \infty$ such that
the corresponding diffeomorphisms $F_2^{-m(j)} \circ
g_{2, j} \circ F_2^{m(j)}$ converge in the $C^1$-topology on $(-\varepsilon , \varepsilon)$ to a non-trivial translation.
There also follows that the vector field $\partial /\partial x$ on $(-\varepsilon , \varepsilon)$ is contained in
the $C^1$-closure of $G_2$. Since $q \simeq 0 \in (-\varepsilon , \varepsilon)$ is clearly expanding for $G_2$,
there immediately follows that every point in $(-\varepsilon , \varepsilon)$ must be expanding for $G_2$. The lemma
follows since $G_2$ acts minimally on $S^1$.\qed

Now it is convenient to revisit the notion of stationary measures in fuller detail.
Consider a finite generating set $A = \{ \overline{\gamma}_{1} ,\ldots , \overline{\gamma}_{r},
\overline{\gamma}_{1}^{-1} , \ldots , \overline{\gamma}_{r}^{-1} \}$ for $\Gamma$ containing elements and their
inverses (so that $A$ generates $\Gamma$ as semigroup). Fix some non-degenerate, probability
measure $\mu$ on $\Gamma$ so that $\mu$ gives strictly positive mass to every element of $A$.
Note that the measure $\mu$ {\it is not}\, required to be symmetric and, in addition,
the set $A$ can be strictly contained in the support of $\mu$. Next define the {\it entropy of $\mu$}\, by
\begin{equation}
H (\mu) = - \sum_{\gamma \in \Gamma} \mu (\gamma ) \ln \mu (\gamma) . \label{entropydefinition}
\end{equation}
For the time being let $\mu$ be a measure as above having finite entropy. A specific choice of $\mu$
will be made only later in connection with Proposition~\ref{regularstationary-addedinApril} and with
the construction in \cite{connel}.

Now denote by $\partial \Gamma$ the geometric boundary of the hyperbolic group
$\Gamma$, see \cite{GhysHarpe}. The boundary $\partial \Gamma$ is a compact metric space which is effectively acted upon
by the group $\Gamma$ itself. Thus we often identify an element $\gamma \in \Gamma$ with the corresponding automorphism
of $\partial \Gamma$ which is still denoted by $\gamma$.

Since $\Gamma$ is endowed with the measure $\mu$, a unique {\it stationary measure}\, $\nu_{\Gamma}$ on $\partial \Gamma$
is associated to the action of $\Gamma$ on $\partial \Gamma$ (cf. \cite{KaimanovichPoisson}). In other words, for every
Borel set $\mathcal{B} \subset \partial \Gamma$, we have
$$
\nu_{\Gamma} (\mathcal{B}) = \sum_{\gamma \in \Gamma} \, \mu (\gamma) \nu_{\Gamma} (\gamma^{-1} (\mathcal{B}))
$$
where $\gamma_{i} (\mathcal{B})$ refers to the identification of $\gamma \in \Gamma$ with the corresponding automorphism
of $\partial \Gamma$.

Next let $g_{1,i} \in G_1$ (resp. $g_{2,i} \in G_2$) be defined as $g_{1,i} = \rho_1 (\overline{\gamma}_{i})$
(resp. $g_{2,i} = \rho_2 (\overline{\gamma}_{i})$), $i=1, \ldots ,r$.
We also pose $A_1 = \{ g_{1,1} , \ldots , g_{1,r}, g_{1,1}^{-1},
\ldots , g_{1,r}^{-1} \}$ and $A_1 = \{ g_{2,1} , \ldots , g_{2,r}, g_{2,1}^{-1}, \ldots , g_{2,r}^{-1} \}$. Since
both representations $\rho_1$ and $\rho_2$ from $\Gamma$ to $\dif$ are one-to-one, the groups
$G_1$ and $G_2$ become equipped with the probability measure $\mu$ up to the obvious identifications.

Going back to the action of $G_1$ on $S^1$, Lemma~\ref{FinishingJob-11.1} allows us to apply
the main theorem of \cite{dkn} to ensure the existence of a unique stationary
measure $\nu_1$ for $G_1$ (with respect to $\mu$). The support of $\nu_1$ is all of
$S^1$ since $G_1$ is minimal. It is also well known that $G_1$ gives no mass to points. Analogous conclusions hold
for the stationary measure $\nu_2$ on $S^1$ arising from $G_2$ and $\mu$. Now the combination of \cite{deroinETDS}
with \cite{KaimanovichPoisson} yields the following.

\begin{lemma}
\label{FinishingJob-11.7}
There is a measurable isomorphism $\theta_2$ from $(\partial \Gamma , \nu_{\Gamma})$ to $(S^1, \nu_2)$.
\end{lemma}

\noindent {\it Proof}. Whereas $G_2$ was initially assumed to be locally $C^2$-discrete,
Lemma~\ref{C1discrete-addedinApril} shows that $G_2$ is, in fact, locally $C^1$-discrete.
Recalling that the measure $\mu$ is assumed to have finite entropy, we apply
Theorem~1.1 of \cite{deroinETDS} to the action of $G_2$ on $S^1$. Since $G_2$ is locally $C^1$-discrete
and $\mu$ has finite entropy, all the conditions required by the theorem in question are satisfied so that
the Poisson boundary of $G_2$ coincides with its $(G_2, \mu)$-boundary (see \cite{deroinETDS}, \cite{connel} for terminology).

In turn, Kaimanovich theorem in \cite{KaimanovichPoisson} ensures that the Poisson boundary of $G_2$ can be identified
with $(\partial \Gamma , \nu_{\Gamma})$ (recall that $G_2$ is isomorphic to the fixed hyperbolic group $\Gamma$).
Thus, to complete the proof of the lemma, it suffices to show that $(G_2, \mu)$-boundary of $G_2$ can be identified
with $(S^1, \nu_2)$. For $G_2$ proximal (and leaving no probability measure invariant, see Lemma~\ref{FinishingJob-11.1}),
this is exactly the contents of \cite{antonov} and \cite{kleptsynnalski}. The reader willing to avoid using
Lemma~\ref{FinishingJob-11.5} may proceed as follows. Consider the finite quotient of $S^1$ where $G_1$ induces a proximal action.
This quotient in endowed with a unique stationary measure $\nu_2'$. The pair $(S^1, \nu_2')$ is the
$(G, \mu)$-boundary of the quotient owing to the result of Antonov and Kleptsyn-Nal'ski. Finally, the $(G_2, \mu)$-boundary
of $G_2$ can then be identified with $S^1$ equipped with the pull-back (still denoted by $\nu_2$) of
$\nu_2'$ by the projection map. This completes the proof of the lemma.\qed

It is implicitly understood in the statement of Lemma~\ref{FinishingJob-11.7} that $\theta_2$ is $\Gamma$-equivariant in
the sense that $\theta_2^{\ast} \nu_2 = \nu_{\Gamma}$ and
\begin{equation}
\theta_2 \circ \gamma  (x) = \rho_2 (\gamma) \circ \theta_2 (x) \label{measurableconjugation}
\end{equation}
for every $\gamma \in \Gamma$ and $\nu_{\Gamma}$-almost all point $x \in \partial \Gamma$.

We can now further specify the choice of the measure $\mu$. This choice is summarized
by Proposition~\ref{regularstationary-addedinApril} below which is based on the construction
developed in \cite{connel}. To state Proposition~\ref{regularstationary-addedinApril} let $G \subset \dif$ be a finitely
generated locally $C^2$-non-discrete group leaving no probability measure on $S^1$ invariant
(in particular $G$ is not solvable and its action on $S^1$ is minimal).

\begin{prop}
\label{regularstationary-addedinApril}
Let $G \subset \dif$ be a group satisfying the above conditions. Then
there exists a measure $\mu$ on $G$ so that the following holds:
\begin{itemize}
  \item $\mu$ is non-degenerate and has finite entropy.

  \item The resulting stationary measure on $S^1$ is absolutely continuous. In fact, this stationary measure
  coincides with the Haar measure on $S^1$.
\end{itemize}
\end{prop}

\noindent {\it Proof}. As mentioned, in the wake of the work of Connell and Muchnik in \cite{connel}, the content
of this proposition seems to have become known to some experts.
Also, compared to the general setting of \cite{connel}, the case of a group acting on the circle is rather particular
and ideas from Furstenberg \cite{furst} might as well be sufficient to establish the statement in question.
In any event, we shall content ourselves
of explaining the main points of the proof while referring to \cite{connel} for additional information.
Consider the a finitely generated group $G \subset \dif$ as in the statement of the proposition.

Denote by ${\rm Leb}$ the Haar measure on $S^1$. We look for a (non-degenerate) measure $\mu$ on $G$ satisfying
$\mu \ast {\rm Leb} = {\rm Leb}$ or, equivalently, such that
$$
{\rm Leb}\, (\calb) =  \sum_{g \in G} \mu (g) \, {\rm Leb}\, (g^{-1} (\calb))
$$
for every Borel set $\calb \subseteq S^1$. As will follow
from the discussion, the construction of $\mu$ has a certain degree of flexibility. This allows us to ensure that $\mu$
is non-degenerate. In any event, the subsequent construction makes it clear that $\mu$ can be chosen so as to take
strictly positive values on a finite set of $G$ generating a locally
 $C^2$-non-discrete hyperbolic group. Ultimately this weaker statement would be sufficient for our purposes.

Our strategy consists of showing how the
Basis Theorem of \cite{connel} (Theorem~6.2 in page 736 of \cite{connel}) can be used to find solutions $\mu$ for
the preceding equation. The fact that solutions with finite entropy can also be found is just an additional
elaboration which is carried out in Section~7 of \cite{connel}
(cf. Theorem~7.1, page 745 of \cite{connel}); this elaboration can be omitted from our discussion.

Recall that a probability measure $\widetilde{\nu}$ on $S^1$ is said to have {\it $(Q,\theta)$-decay}\, if
there is a constant ${\rm Const} >0$ such that
$$
\int_{S^1 \setminus B(x,r)} \frac{1}{d\, (y,x)^{Q+\theta}} d\widetilde{\nu} (y) \leq \frac{\rm Const}{r^Q}
$$
for every $x \in S^1$ and $1 \geq r >0$. The distance ``$d$'' considered here is again the Euclidean distance giving
length one to $S^1$. Similarly $B (x,r) \subset S^1$ stands for the ball (interval) of radius~$r$ around $x \in S^1$.
Finally, in case $Q=0$, the right hand side in the above estimate should be replaced by ${\rm Const} (1 + \vert \ln r \vert)$.
There follows from this definition that the Haar measure ${\rm Leb}$ on $S^1$ has $(Q,\theta)$-decay for every pair
$(Q,\theta)$ with $Q \in \R$ and $\theta \geq 1$.

In the sequel we will always focus on the Haar measure ${\rm Leb}$ on $S^1$. Fix then a pair
$(Q,\theta)$ with $Q \geq 0$ and $\theta \geq 1$ so that ${\rm leb}$ has $(Q,\theta)$-decay. Adapting the terminology
of \cite{connel} to the present context, a {\it continuous spike}\, is defined as a $6$-tuple
$(\zeta (x), r,a,Q, \theta, C)$ where $\zeta$ is a continuous positive function on $S^1$, $r>0$, $a \in S^1$,
and $C >1$ which satisfies the following conditions:
\begin{enumerate}
  \item $\zeta(x) \geq \Vert \zeta \Vert_{L^{\infty}}/C$ on $B(a,r)$ (where $B(a,r)$ denotes the interval
  of radius~$r$ around $a \in S^1$).

  \item For each $y \in S^1 \setminus B(a,r)$, we have
  $$
  0 \leq \zeta (y) \leq \zeta (a) r^Q \int_{B(a,r)} \frac{C}{d\, (y,x)^{Q+\theta}} {\rm Leb}\, (x) \, .
  $$

  \item For every two points $y, y' \in S^1$ such that $d(y,y') \leq r$, we have $\zeta (y') < C \zeta (y)$.
\end{enumerate}
A (continuous) spike is said to be an {\it unit (continuous) spike}\, if $\Vert \zeta \Vert_{L^{\infty}} =1$.

Keeping $Q$ and $\theta$ fixed, consider a family $\fol_{\alpha} = \{ (\zeta_{\alpha}, r_{\alpha}, a_{\alpha},
Q, \theta, C_{\alpha}) \}$ of continuous spikes indexed by $\alpha \in \mathcal{A}$. For every $C >1$,
let $S_C = \{ \alpha \in \mathcal{A} \; ; \; C_{\alpha} < C \}$. Now we pose
$$
B_C (r) = \bigcup_{\alpha in S_C \\ r_{\alpha} < r} B (a_{\alpha} , r_{\alpha}) \; \; \; {\rm and} \; \; \;
B_C = \bigcap_{r >0} B_C (r) \, .
$$
With this notation, assume that the family $\fol_{\alpha}$ satisfies the following condition:

\noindent Condition~($\ast$) For $C$ large enough, we have $S^1 \subseteq B_C$.

Consider now a family $\fol_{\alpha}$ of {\it unit}\, continuous spikes satisfying Condition~($\ast$).
For this family, Theorem~6.2 (Basis Theorem) in \cite{connel}
implies the existence of countably many indices $\alpha_i$ along with associated nonnegative numbers $c_i$ such that
$$
\sum_{i=1}^{\infty} c_i \zeta_{\alpha_i} (x) = 1
$$
with uniform convergence on $S^1$.

The proof of Proposition~\ref{regularstationary-addedinApril} will hinge from this theorem. To exploit it, we first
note that finite sets can always be eliminated from the family $\fol_{\alpha}$. In other words, if $\fol_{\alpha}$
satisfies Condition~($\ast$) and $\fol_{\alpha}'$ is a new family of spikes obtained from $\fol_{\alpha}$ by eliminating
finitely many spikes, then $\fol_{\alpha}'$ satisfies Condition~($\ast$) as well. Another simple observation is that
spikes can naturally be renormalized to become unit spikes in the following sense: if
$(\zeta_{\alpha}, r_{\alpha}, a_{\alpha}, Q, \theta, C_{\alpha})$ is a (continuous) spike, then
$(\zeta_{\alpha} /\Vert \zeta_{\alpha} \Vert_{L^{\infty}} , r_{\alpha}, a_{\alpha}, Q, \theta, C_{\alpha})$ is an unit
(continuous) spike. We shall return to these points later.

Next, for every $g \in G$, let $\zeta_g$ denote the Radon-Nikodym derivative $d(g_{\ast} \Leb )/d \, \Leb$ which coincides
with the usual derivative of $g^{-1}$. Now we have:

\vspace{0.1cm}

\noindent {\it Claim~1}. Assume there is
a subset $\{ g_{\alpha} \} \subset G$ of elements in $G$ yielding a set of function $\{ \zeta_{\alpha} \}$ to which it is
possible to associate spikes $\{ (\zeta_{\alpha}, r_{\alpha}, a_{\alpha}, Q, \theta, C_{\alpha}) \}$ so that the
resulting family $\fol_{\alpha}$ of spikes satisfies Condition~($\ast$). Then Proposition~\ref{regularstationary-addedinApril}
holds.

\noindent {\it Proof of the Claim~1}. Up to renormalizing the functions $\zeta_{\alpha}$, which become
$\zeta_{\alpha}/ \Vert \zeta_{\alpha} \Vert_{L^{\infty}}$, the family $\fol_{\alpha}$ can
be regarded as a family of unit continuous spikes. Therefore there are nonnegative constants $c_{\alpha}$ such that
$\sum_{\alpha} c_{\alpha} \zeta_{\alpha} (x) / \Vert \zeta_{\alpha} \Vert_{L^{\infty}} =1$.
Define then the auxiliary measure $\mu$ by letting $\mu (g_{\alpha}) = c_{\alpha}/ \Vert \zeta_{\alpha} \Vert_{L^{\infty}}$.
Since we have uniform convergence and the integral of $\zeta_{\alpha}$ on $S^1$ equals~$1$, we conclude that
$$
\sum_{g \in G} m(g) = \sum_{\alpha} \frac{c_{\alpha}}{\Vert \zeta_{\alpha} \Vert_{L^{\infty}}} =1 \,.
$$
Proposition~\ref{regularstationary-addedinApril} now follows from observing that
$$
\int_{\calb}\sum_{g \in G} \mu (g) \, d (g_{\ast} {\rm Leb}) =
\int_{\calb} \sum_{\alpha} \frac{c_{\alpha}}{\Vert \zeta_{\alpha} \Vert_{L^{\infty}}}
d (g_{\ast} {\rm Leb}) = \int_{\calb} \sum_{\alpha} \frac{c_{\alpha} \zeta_{\alpha}}{\Vert \zeta_{\alpha} \Vert_{L^{\infty}}}
d{\rm Leb} = {\rm Leb}\, (\calb) \, .
$$
\mbox{ }\qed

The next step consists of a natural relaxation on the definition of spikes. Consider a family $\fol_{\alpha}$
as in Claim~1 and, up to take the renormalization, assume that the spikes in $\fol_{\alpha}$ are already unit spikes.
Fix an element $(\zeta_{\alpha}, r_{\alpha}, a_{\alpha}, Q, \theta, C_{\alpha})$ in $\fol_{\alpha}$.
Conditions numbers~$2$ and~$3$ in the definition of a spike impose restrictions on the
behavior of $\zeta_{\alpha}$ away from $a_{\alpha}$ and these conditions can significantly be relaxed. Here however,
a more accurate discussion is needed. To begin with, we fix a sequence of elements in $\fol_{\alpha}$ with
$C_{\alpha}$ uniformly bounded and such that $r_{\alpha} \rightarrow 0$. Note that Condition~($\ast$) ensures that
$\fol_{\alpha}$ must, in fact, contain sequences as above. We can also assume that $a_{\alpha}$ converges towards
some point $a_{\infty} \in S^1$ so as to fix a small uniform neighborhood $U \subset S^1$ of $a_{\infty}$. If $y$ lies
in $S^1 \setminus U$, we must have (by assumption)
$$
\zeta_{\alpha} (y) \leq \zeta_{\alpha} (a_{\alpha}) r_{\alpha}^Q \int_{B(a_{\alpha},r_{\alpha})} \frac{C_{\alpha}}{
d(y,x)^{Q+\theta}} dx \leq {\rm Const}\, . \, r_{\alpha}^{Q+1}
$$
for a uniform constant ${\rm Const}$ and provided that $r_{\alpha}$ is sufficiently small. In particular $\zeta_{\alpha} (y)$
decays as $r_{\alpha}^{Q+1}$. Our purpose is to weaken this condition as well as condition number~$3$ away from
$a_{\alpha}$. For this we are led to define what we call a {\it family of local (unit, continuous) spikes
satisfying Condition~($\ast$)}. We consider again a family $\fol_{\alpha}' =
\{ (\zeta_{\alpha}, r_{\alpha}, a_{\alpha}, Q, \theta, C_{\alpha}) \}$ of $6$-tuples as before, except that the
conditions on the functions $\zeta_{\alpha}$ will be weaker. In particular we can still consider whether or not
Condition~($\ast$) is satisfied by $\fol_{\alpha}'$ since Condition~($\ast$) does not depend on the functions $\zeta_{\alpha}$.
Concerning the functions $\zeta_{\alpha}$, we only impose the existence of a uniform $\delta >0$ such that the
following holds:
\begin{itemize}
  \item[(A)] On the ball $B (a_{\alpha}, \delta)$, the function $\zeta_{\alpha}$ satisfies all the conditions
required by the definition of spikes (with $\Vert \zeta_{\alpha} \Vert_{L^{\infty}}$ bounded from below by
some positive constant or simply equal to~$1$).

  \item[(B)] The functions $\zeta_{\alpha}$ converge uniformly to zero on $S^1 \setminus B (a_{\alpha}, \delta)$.
  More precisely, for $r_{\alpha}$, $a_{\alpha}$, and $C_{\alpha}$ fixed, we can choose $\zeta_{\alpha}$ arbitrarily
  close to zero on $S^1 \setminus B (a_{\alpha}, \delta)$.
\end{itemize}
Families of local (unit, continuous) spikes are easier to be constructed since we only need to ``fine
tune'' the functions $\zeta_{\alpha}$ on a small interval, as opposed to on all of $S^1$. Yet, these families can be
turned into families of actual spikes by adding to $\zeta_{\alpha}$ a positive function. More precisely, consider
a family $\{ (\zeta_{\alpha}, r_{\alpha}, a_{\alpha}, Q, \theta, C_{\alpha}) \}$ of local (unit, continuous) spikes,
satisfying Condition~($\ast$). Up to passing to a subset of indices, we can find
a sub-family $\{ (\overline{\zeta}_{\alpha}, r_{\alpha}, a_{\alpha}, Q, \theta, C_{\alpha}) \}$ of
spikes satisfying Condition~($\ast$) and such that
$\overline{\zeta}_{\alpha} = \zeta_{\alpha} + R_{\alpha}$, where $R_{\alpha}$ is a positive function.
This can be done as follows. Consider an element $(\zeta_{\alpha}, r_{\alpha}, a_{\alpha}, Q, \theta, C_{\alpha})$
in the initial family. Owing to Condition~(B), we can assume that the restriction of $\zeta_{\alpha}$ to
$S^1 \setminus B (a_{\alpha}, \delta)$ is arbitrarily small and, in particular, small enough to satisfy the corresponding
decay condition (depending only on $r_{\alpha}$, $Q$, and $C_{\alpha}$). Then we first set
$R_{\alpha} = \sup_{x \in S^1 \setminus B (a_{\alpha}, \delta)} \zeta_{\alpha} (x)$. Thus we just need to make
$R_{\alpha}$ decay to zero on an one-sided neighborhood of the two points corresponding to the boundary of $B (a_{\alpha}, \delta)$.
It is clear that this can be done while keeping the conditions verified by $\overline{\zeta}_{\alpha}$ on
$B (a_{\alpha}, \delta)$. In view of this, we now obtain:

\vspace{0.1cm}

\noindent {\it Claim~2}. To prove Proposition~\ref{regularstationary-addedinApril} it suffices to find
a subset $\{ g_{\alpha} \} \subset G$ of elements in $G$ yielding a set of function $\{ \zeta_{\alpha} \}$ to which it is
possible to associate a family of local spikes $\fol_{\alpha} = \{ (\zeta_{\alpha}, r_{\alpha}, a_{\alpha}, Q, \theta,
C_{\alpha}) \}$ satisfying Condition~($\ast$).

\noindent {\it Proof of the Claim~2}. Consider the associated family $\overline{\fol}_{\alpha} = \{
(\overline{\zeta}_{\alpha}, r_{\alpha}, a_{\alpha}, Q, \theta, C_{\alpha}) \}$ of spikes, with
$\overline{\zeta}_{\alpha} = \zeta_{\alpha} + R_{\alpha}$ for some $R_{\alpha} \geq 0$. Now
Theorem~6.2 of \cite{connel} ensures the existence of a set $\mathcal{A}^1$
such that $\sum_{\alpha \in \mathcal{A}^1} c_{\alpha} \overline{\zeta}_{\alpha} = 1$ for suitable reals $c_{\alpha} >0$.
Let now $\varepsilon >0$ be fixed. Since $c_{\alpha}$ and $\overline{\zeta}_{\alpha}$ are positive for all $\alpha \in
\mathcal{A}^1$, the preceding implies the existence of a finite set $\mathcal{A}_1^1 \subset \mathcal{A}^1$ such that
$\overline{L}_1 = 1 - \sum _{\alpha \in \mathcal{A}_1^1} c_{\alpha} \overline{\zeta}_{\alpha}$ is still a positive function
and, furthermore, satisfies
$$
\sup_{x \in S^1} \vert \overline{L}_1 (x) \vert = \sup_{x \in S^1} \overline{L}_1 (x) \leq \varepsilon/2.
$$
In addition, the functions $R_{\alpha}$ can be assumed arbitrarily small so that the convergence scheme of \cite{connel}
(Section~6) allows us to obtain
\begin{equation}
L_1 = 1- \sum _{\alpha \in \mathcal{A}_1^1} c_{\alpha} \zeta_{\alpha} = \overline{L}_1 +
\sum_{\alpha \in \mathcal{A}_1^1} c_{\alpha} R_{\alpha} \leq \varepsilon \, . \label{L1withRalphapositive}
\end{equation}
Clearly $L_1$ is positive. Furthermore, as already mentioned, the family of spikes $\overline{\fol}_{\alpha}$ still
form a basis even after eliminating finitely many terms from it. Thus, we can now apply Theorem~6.2 of \cite{connel}
to find a set $\mathcal{A}^2$, $\mathcal{A}^2 \cap \mathcal{A}_1^1 = \emptyset$, such that
$$
\sum _{\alpha \in \mathcal{A}^2} c_{\alpha} \overline{\zeta}_{\alpha} = L_1
$$
for certain constants $c_{\alpha} >0$.

We can now repeat the above argument to find a finite set $\mathcal{A}_1^2 \subset \mathcal{A}^2$ such that
$\overline{L}_2 = L_1 - \sum _{\alpha \in \mathcal{A}_1^2} c_{\alpha} \overline{\zeta}_{\alpha} \leq \varepsilon/4$
with $\sum _{\alpha \in \mathcal{A}_1^2} c_{\alpha} R_{\alpha} \leq \varepsilon /4$. Therefore
$$
L_2 = 1- \sum _{\alpha \in \mathcal{A}_1^1} c_{\alpha} \zeta_{\alpha} - \sum _{\alpha \in \mathcal{A}_1^2} c_{\alpha}
\zeta_{\alpha} \leq \varepsilon /2
$$
is still a positive function. The procedure can then be continued by induction to derive the existence of a
set $\widetilde{\mathcal{A}}$ such that $\sum _{\alpha \in \widetilde{\mathcal{A}}} c_{\alpha} \zeta_{\alpha} =1$.
The remainder of the proof of Claim~2 is totally analogous to the proof of Claim~1.\qed

Proposition~\ref{regularstationary-addedinApril} is now reduced to finding a family of local spikes
$\fol_{\alpha} = \{ (\zeta_{\alpha}, r_{\alpha}, a_{\alpha}, Q, \theta, C_{\alpha}) \}$ satisfying Condition~($\ast$)
and obtained from elements $g_{\alpha} \in G$ by the above indicated procedure. The remainder of the proof
is devoted to constructing these elements in $G$. To fix terminology, we will say that a family
contains {\it enough (local) spikes}\, if it satisfies Condition~($\ast$). Similarly, $G$ is said to contain
enough (local) spikes if it yields (through the already mentioned procedure)
a family containing enough (local) spikes. Finally, fixed a point
$p \in S^1$, we say that $G$ contains enough (local) spikes {\it at $p$} if $G$ yields a family of spikes
centered at $p$ which, once enlarged by the natural effect of Euclidean rotations leads to a family
containing enough (local) spikes.

To abridge notation we can assume that $G$ is not conjugate to a finite covering of a subgroup of ${\rm PSL}\, (2,\R)$.
Otherwise the existence of the desired elements is well known (and established in a much more general setting in
\cite{connel}). Therefore Theorem~3.4 of \cite{rebproceedings} yields a finite covering $\{ J_1, \ldots, J_r \}$
of $S^1$ such that the following holds: for every compactly contained interval $(a,b) \subset J_i$ and every
$C^{\infty}$-diffeomorphism $\phi : (a,b) \rightarrow \phi (a,b) \subset J_i$, there is a sequence $\{ g_k \}$
of elements in $G$ whose restrictions to $(a,b)$ converge to $\phi$ in the $C^{\infty}$-topology.

On the other hand, as previously mentioned, the group $G$ contains a diffeomorphism $F$ possessing exactly
two fixed points $p,q$ in $S^1$. Furthermore these two fixed points are hyperbolic with multipliers $\lambda_1$
and $\lambda_2$ (say $0 < \lambda_1 <1 < \lambda_2$ so that $p$ is an attracting fixed point while $q$ is
a repelling one, cf. Th\'eor\`eme~F in \cite{dkn}). Naturally $F$ is linearizable around both $p$ and $q$.

Fix $\delta >0$ so small that every ball of radius $\delta$ around a point $x \in S^1$ is fully contained in
one of intervals $J_1, \ldots, J_r$. We shall first prove:

\vspace{0.1cm}

\noindent {\it Claim~3}. For $\delta >0$ as above, the group $G$ contains enough local spikes centered
at $p$.

\noindent {\it Proof of Claim~3}. Without loss of generality we assume that $B (p, \delta)$ is contained
in $J_1$. The existence of $F$ allows us to make the interval $J_1$ ``global'' in the following sense.
Consider an interval $\overline{J} = (\overline{a}, \overline{b})$ compactly contained in
$S^1 \setminus \{ q\}$. Clearly there is a $n_0 \in \N^{\ast}$
such that $F^{n_0} (\overline{J}) \subset J_1$. Thus, by virtue of the above stated property of $J_1$ and $G$,
we can find a sequence of elements $\{ g_k \} \subset G$ such that the following holds:
\begin{enumerate}
  \item $g_k (\overline{a}) \rightarrow \overline{a}$ and $g_k (\overline{b}) \rightarrow \overline{b}$
  \item The restriction of each element in $\{ g_k' \}$ to $\overline{J}$
  yields a sequence of suitable spikes centered around $p$ in $\overline{J}$.
\end{enumerate}
The complete the proof of Claim~3 we just need to consider the behavior of the sequence of derivatives
$\{ g_k' \}$ on $S^1 \setminus \overline{J}$. Since $\overline{J}$ is ``almost fixed'' by $g_k$, it is clear
that the $L^1$-norm of $g_k'$ on $S^1 \setminus \overline{J}$ is small in the sense that it is comparable to the length of
$S^1 \setminus \overline{J}$. To see that the $C^0$-norm of $g_k'$ can be made small as well, we proceed
as follows. We can assume that $S^1 \setminus \overline{J}$ is contained in the domain of linearization of
$F$ around $q$. Thus, for each $k$ fixed, consider diffeomorphisms of the form $F^{-1} \circ g_k$, $F^{-2} \circ g_k$,
and so on. Since $(F^{-1})' (x) < 1$ on $S^1 \setminus \overline{J}$, the derivative of $F^{-N_k} \circ g_k$
is small on $S^1 \setminus \overline{J}$ provided that $N_k$ is very large. Thus we just need to consider
the behavior of $(F^{-N_k})' (g_k) g_k'$ on $\overline{J}$. This however is essentially a re-scaling of
the behavior of $g_k'$ so that the claim follows by renormalizing again these maps into unit spikes.\qed

To finish the proof of Proposition~\ref{regularstationary-addedinApril} we still need to show that
$G$ contains enough (local) spikes from the fact that it contains enough (local) spikes {\it at $p$}.
This however follows from the existence of constant vector fields on the closure of $G$ (cf. Theorem~3.4 of \cite{rebproceedings})
that locally behave as rotations. Proposition~\ref{regularstationary-addedinApril} is proved.\qed

In the sequel we fix $\mu$ on $\Gamma$ as in Proposition~\ref{regularstationary-addedinApril} so
that the resulting stationary measure $\nu_1$ for $G_1$ is the Haar measure.
We are now ready to prove Theorem~\ref{PoissonBoundaryApplication}.

\vspace{0.1cm}

\noindent {\it Proof of Theorem~\ref{PoissonBoundaryApplication}}. Let $h : S^1 \rightarrow S^1$ be a homeomorphism
conjugating $G_1$ to $G_2$. By way of contradiction, we have assumed that $G_1$ is locally $C^2$-non-discrete whereas
$G_2$ is locally $C^2$-discrete. Recall also
that $\nu_1$ (resp. $\nu_2$) is the unique stationary measure for $G_1$ (resp. $G_2$) with respect to $\mu$
(see \cite{dkn}). From this uniqueness, there follows again that $h^{\ast} \nu_1 = \nu_2$. Furthermore $\nu_1$ coincides
with the Haar measure.

Consider the measurable isomorphism $\theta_2 : (\partial \Gamma , \nu_{\Gamma}) \longrightarrow (S^1 ,\nu_2)$
of Lemma~\ref{FinishingJob-11.7} and define
a new measurable isomorphism $\theta_1 : (\partial \Gamma , \nu_{\Gamma}) \longrightarrow (S^1 ,\nu_1)$ by letting
$\theta_1 = h \circ \theta_1$. The equivariant nature of $\theta_2$ expressed by Equation~(\ref{measurableconjugation})
combines with the fact that $h^{\ast} \nu_1 = \nu_2$ to yield
\begin{equation}
\theta_1 \circ \gamma (x) = \rho_1 (\gamma) \circ \theta_1 (x) \label{measurableconjugation-forgroup1}
\end{equation}
for every $\gamma \in \Gamma$ and $\nu_{\Gamma}$-almost all point $x \in \partial \Gamma$. Furthermore $\theta_1^{\ast} \nu_1
=\nu_{\Gamma}$. Up to eliminating null measure
sets, we fix once and for all a Borel set $\mathcal{B} \subset \partial \Gamma$ having full $\nu_{\Gamma}$-measure
and such that Equation~(\ref{measurableconjugation-forgroup1}) holds for every $x \in \mathcal{B}$ and every $\gamma \in \Gamma$
(in particular both sides of this equation are well defined). To complete the proof of the proposition,
we are going to show that the existence of $\theta_1$ is not compatible with the fact that $G_1$ is locally $C^2$-non-discrete.
To do this, we proceed as follows.

Fix an interval $I \subset S^1$ along with a sequence of elements $\{ g_j \} \subset G_1$, $g_j \neq {\rm id}$ for
every $j \in \N$, whose restrictions to $I$ converge to the identity in the $C^2$-topology. The existence of $I$ and
of $\{ g_j \}$ clearly follows from the assumption that $G_1$ is locally $C^2$-non-discrete. Now
Lusin approximation theorem \cite{bartle} ensures the existence of a Cantor set $K$ satisfying the following conditions:
\begin{enumerate}
  \item $K \subset I \cap \theta_1 (\mathcal{B})$, i.e. $K$ is contained in the domain of definition of
  $\theta_1^{-1}$.

  \item The {\it restriction}\, of $\theta_1^{-1}$ to $K$ is continuous from $K$ to $\partial \Gamma$
  (where the reader is reminded that $\partial \Gamma$ is a compact metric space).

  \item $\nu_1 (K) \geq 9 \nu_1 (I) /10$.
\end{enumerate}

\noindent Next, for each~$j$, let $\gamma_j \in \Gamma$ be such that $\rho_1 (\gamma_j) = g_j$.

\vspace{0.1cm}

\noindent {\it Claim}. There is a Cantor set $K_{\Gamma} \subset \partial \Gamma$ such that
the restrictions of the elements $\gamma_j$ to $K_{\Gamma}$ converge uniformly to the identity.

\noindent {\it Proof of the Claim}. Since $\{ g_j \}$ converges to the identity in the $C^1$-topology and $\nu_1$
coincides with the Haar measure, there follows that $\nu_1 (K \cap g_j^{-1} (K))$ converges to $\nu_1 (K)$ as $j \rightarrow \infty$.
Therefore, up to passing to a subsequence, we can assume that
$$
K_{\infty} =  K \cap \bigcap_{j=1}^{\infty} g_j^{-1} (K)
$$
is an actual (non-empty) Cantor set. Furthermore, by construction, $K_{\infty} \subset K$ and $g_j (K_{\infty}) \subset K$
for every $j \in \N^{\ast}$. Finally let $K_{\Gamma} = \theta_1^{-1} (K_{\infty})$.

To complete the proof of the claim, note that the restriction of $\theta_1$ to $K_{\Gamma}$ is continuous
since $\theta_1^{-1}$ is continuous and one-to-one on the Cantor set $K$ (and $K_{\infty} \subset K$). On the other hand,
on $K_{\Gamma}$ we have
$$
\gamma_j = \theta_1^{-1} \circ g_j \circ \theta_1 \,
$$
i.e. the left hand side is well defined on $K_{\Gamma}$.
Since $\theta_1$ is continuous on $K_{\Gamma}$ and $\theta_1^{-1}$ is continuous on
$\circ g_j \circ \theta_1 (K_{\Gamma}) \subset K$, the fact that $g_j$ converges uniformly (and actually $C^1$)
to the identity implies the claim.\qed

We have just found a sequence $\{ \gamma_j \}$ of elements in $\Gamma$, $\gamma_j \neq {\rm id}$ for every $j$,
whose restrictions to a (non-empty) Cantor set $K_{\Gamma} \subset \partial \Gamma$ converge uniformly to the identity.
The theorem now from Lemma~\ref{longhyperbolicgroup} below claiming that such a sequence cannot exist in a finitely generated
hyperbolic group.\qed

To state Lemma~\ref{longhyperbolicgroup} recall that every element $\gamma \in \Gamma$ can be identified with the
corresponding automorphism of $\partial \Gamma$. Naturally
$\gamma$ can equally well be identified with its translation action on $\Gamma$ which happens to be an isometry for the natural
left-invariant metric on $\Gamma$ (see \cite{GhysHarpe}).

\begin{lemma}
\label{longhyperbolicgroup}
Let $\Gamma$ be a hyperbolic group which is neither finite nor a finite extension of $\Z$. Let $K_{\Gamma}$ be a
Cantor set contained in the boundary $\partial \Gamma$ of $\Gamma$ and let $\{ \gamma_j \}$ be a sequence of
elements in $\Gamma$ thought of as automorphisms of $\partial \Gamma$. Assume that the sequence
$\{ \gamma_{j,\vert K_{\Gamma}} \}$ obtained by restricting $\gamma_j$ to $K_{\Gamma}$ converges uniformly to the
identity. Then we have $\gamma_j = {\rm id}$ for large enough $j \in \N$.
\end{lemma}

\noindent {\it Proof}. The lemma is certainly well known to the specialists albeit we have not been able to find
it explicitly stated in the literature. In the sequel, the reader is referred to the chapters~7 and~8 of \cite{GhysHarpe}
for background material.

Assume for a contradiction that $\gamma_j \neq {\rm id}$ for every $j \in \N$. Consider also a base point $w \in \Gamma$
along with the sequence $\gamma_j (w)$. Since $\gamma_j$ acts as an isometry of $\Gamma$, there follows that
the sequence $\{ \gamma_j (w) \}$ leaves every compact part of $\Gamma$. Thus, up to a passing to a subsequence, we
assume that $\gamma_j (w) \rightarrow b \in \partial \Gamma$.

Next fix another point $a \in \partial \Gamma \setminus K_{\Gamma}$, $a \neq b$, and consider the family
of metrics $d_{\varepsilon , a, w'}$ on $\partial \Gamma \setminus \{ a \}$ for a fixed (small) $\varepsilon >0$
and where $w' \in \Gamma$ (see \cite{GhysHarpe}, page 141). Let $\beta_a$ denote the Busemann function relative
to the point $a \in \partial \Gamma$. Since $\gamma_j (w) \rightarrow b$, with $b \neq a$, there follows from the
general properties of Busemann functions that
$$
\beta_a (w, \gamma_j (w)) \longrightarrow - \infty
$$
(cf. \cite{GhysHarpe} page 136).
In particular, there is some uniform constant $C$ such that
$$
\frac{1}{C} \exp (-\varepsilon \beta_a (w, \gamma_j (w))) \leq \frac{d_{\varepsilon , a,\gamma_j (w)} (x,y)}{
d_{\varepsilon , a, w} (x,y)} \leq C \exp (-\varepsilon \beta_a (w, \gamma_j (w))) \, ;
$$
see \cite{GhysHarpe}, page 141. In other words, the metric $d_{\varepsilon , a,\gamma_j (w)}$ is bounded
from below and by above by the metric $d_{\varepsilon , a, w}$ multiplied by suitable constants going to infinity
as $j \rightarrow \infty$. However,
by construction, these metrics also satisfy $d_{\varepsilon , a,\gamma_j (w)} (\gamma_j (x) ,\gamma_j (y))
= d_{\varepsilon , a, w} (x,y)$. Therefore
$$
\frac{d_{\varepsilon , a, w} (x,y)}{d_{\varepsilon , a, w} (\gamma_j (x) ,\gamma_j (y))} \longrightarrow \infty
$$
uniformly for every pair $x\neq y$ in $\partial \Gamma \setminus \{a \}$. The desired contradiction now arises
by choosing $x \neq y \in K_{\Gamma}$ so that $\gamma_j (x) \rightarrow x$ and
$\gamma_j (y) \rightarrow y$. The proof of the lemma is completed.\qed

\section{Appendix: on locally $C^r$-non-discrete groups}

For $r \geq 2$, every subgroup $G$ of $\dif$ that is locally $C^r$-non-discrete is clearly locally $C^l$-non-discrete
for every $l \leq k$. A sort of converse for the above claim also holds in most cases.
This is the content of the theorem below.

\begin{teo}
\label{C2impliesCinfty}
Let $G \subset \dif$ be a non-solvable group and assume that $G$ is locally $C^2$-non-discrete. Then $G$
is locally $C^{\infty}$-non-discrete.
\end{teo}

To prove Theorem~\ref{C2impliesCinfty} we shall use the same technique of regularization (or renormalization)
employed in Section~4. By assumption there is an open (non-empty) interval $I \subset S^1$ and a sequence
$\{ f_j \}$, $f_j \neq {\rm id}$ for every $j \in \N$, of elements in $G$ whose restrictions to $I$ converge to the
identity in the $C^2$-topology. In fact, arguing as in Section~3, we can assume without loss of generality that
the following holds: for every given $\varepsilon >0$, there is a finite set $\overline{f}_1 ,\ldots , \overline{f}_N$
of elements in $G$ satisfying the two conditions below.
\begin{itemize}
  \item The group $G_{1-N} \subset G$ generated by $\overline{f}_1 ,\ldots , \overline{f}_N$ is not solvable.

  \item For every $i \in \{1, \ldots, N\}$, the restriction of $\overline{f}_i$ to the interval $I$ is
  $\varepsilon$-close to the identity in the $C^2$-topology on $I$.
\end{itemize}

First we state:

\begin{prop}
\label{Appendix-ReductionProposition}
If $\varepsilon >0$ is small enough, then the group $G_{1-N}$
is locally $C^r$-non-discrete for every $r \in \N$.
\end{prop}

The proof of Theorem~\ref{C2impliesCinfty} can be derived from Proposition~\ref{Appendix-ReductionProposition} as
follows.

\vspace{0.1cm}

\noindent {\it Proof of Theorem~\ref{C2impliesCinfty}}. We can assume once and for all that $G_{1-N}$ has no finite orbits,
otherwise Theorem~\ref{C2impliesCinfty} follows at once from the discussion in Section~2.1. In turn, it is clearly sufficient
to prove that the subgroup
$G_{1-N}$ is locally $C^{\infty}$-non-discrete provided that $\varepsilon >0$ is small enough. This is equivalent to
finding an open, non-empty interval $I_{\infty} \subset S^1$ on which ``$G_{1-N}$ is locally $C^{r}$-non-discrete for every
$r \in \N$''. More precisely, for every fixed $r \in \N$, there is a sequence $\{ f_{j, C^r} \}_{j \in \N}$,
$f_{j, C^r} \neq {\rm id}$ for every $j \in \N$, of elements in $G_{1-N}$ whose restrictions to $I_{\infty}$ converge
to the identity in the $C^r$-topology.

On the other hand, by assumption, to every $r \in \N$ there corresponds a non-trivial sequence $\{ \tilde{f}_{j, C^r} \}_{j \in \N}$
of elements in $G_{1-N}$ whose restriction to some open, non-empty interval $I_r$ converges to the identity
in the $C^r$-topology on $I_r$. Thus the only difficulty to derive Theorem~A lies in the fact that the intervals $I_r$
depend on $r$. To show that these intervals can be chosen in a uniform way, we proceed as follows.

First recall that $G_{1-N}$ contains and element $F$ exhibiting a hyperbolic fixed point. Furthermore
$S^1$ can be covered by finitely many
intervals $J_1, \ldots ,J_l$ such that each interval $J_i$ is equipped with a constant (non-zero) vector field $X_i$ in the
$C^1$-closure of $G_{1-N}$; cf. Theorem~3.4 of \cite{rebproceedings} (which, in particular, recovers the fact that
all orbits of $G$ are dense in $S^1$). By using these constant vector fields and the diffeomorphism $F$,
we obtain a sequence $F_r$ of elements in $G_{1-N}$ satisfying the following conditions:
\begin{itemize}
  \item The diffeomorphism $F_r$ has an attracting hyperbolic fixed point $p_r$ lying in $I_r$.

  \item The basin of attraction of $p_r$ with respect to $F_r$ has length greater than a certain $\delta >0$
  (in other words, there is $\delta >0$ such that $F_r$ has no other fixed point on a $\delta$-neighborhood
  of $p_r$).
\end{itemize}

Now each interval $I_r$ can be ``re-scaled'' by means of $F_r$ so as to have length bounded from below by $\delta$.
More precisely, fixed $r$ and $n_r \in \N$, the sequence of elements of $G_{1-N}$ given by
$j \mapsto F_r^{-n_r} \circ \tilde{f}_{j, C^r} \circ F_r^{n_r}$ clearly converges to the identity in the $C^r$-topology
on the interval $\tilde{I}_r = F_r^{-n_r} (I_r)$. The above stated conditions on the diffeomorphisms $F_r$ then
ensure that $n_r$ can be chosen
so that $\tilde{I}_r = F_r^{-n_r} (I_r)$ has length bounded from below by $\delta >0$. Up to passing to a subsequence,
the sequence of intervals $\{ \tilde{I}_r \}$ must converge to a uniform interval $I_{\infty}$ satisfying the desired
conditions. The proof of Theorem~\ref{C2impliesCinfty} is completed.\qed

As in Section~4, we consider the sequence of sets $S(k)$ defined by means of the initial set $S=S(0) =
\{ \overline{f}_1 ,\ldots , \overline{f}_N \}$. Since the group generated by $\overline{f}_1 ,\ldots , \overline{f}_N$
is not solvable, none of the sets $S(k)$ is reduced to the identity diffeomorphism.

We can now prove Proposition~\ref{Appendix-ReductionProposition}.

\vspace{0.1cm}

\noindent {\it Proof of Proposition~\ref{Appendix-ReductionProposition}}. The proof is essentially by induction. First we
are going to prove that $G_{1-N}$ is locally $C^3$-non-discrete. To do this, we proceed as follows.
Consider a fixed set $\{ \overline{f}_1 ,\ldots , \overline{f}_N \}$ generating a non-solvable
group $G_{1-N}$ as before. Assume moreover that for every $i=1, \ldots ,N$, both diffeomorphisms $\overline{f}_i$
and $\overline{f}_i^{-1}$ are $\varepsilon$-close to the identity in the $C^2$-topology on $I$ where the value of $\varepsilon >0$ will
be fixed later on.

As already seen, the group $G$ contains an element $F$ exhibiting a hyperbolic fixed point
in $I$. Without loss of generality, we can assume that this fixed point coincides with $0 \in I \subset \R$. Furthermore
in suitable coordinates, $F$ becomes a homothety $x \mapsto \lambda x$ on all of the interval $I$. Still
keeping the notation of Section~4, consider the sequence of sets $\widetilde{S}(k)$ given by
$\widetilde{S}(k) = F^{-kn} \circ S(k) \circ F^{kn}$ for some $n \in \N^{\ast}$ fixed. We will show that the diffeomorphisms in
$\widetilde{S}(k)$ converge to the identity in the $C^3$-topology on $I$ provided that $n$ is suitably chosen.

\vspace{0.1cm}

\noindent {\it Claim}. There is $n \in \N$ such that every non-trivial sequence $\{ \tilde{f}_k \}$, with
$\tilde{f}_k \in \widetilde{S}(k)$, converges to the identity in the $C^3$-topology on $I$.

\noindent {\it Proof of the Claim}. Fix a sequence $\{ \tilde{f}_k \}$ as in the statement.
In Section~4 it was seen that these elements converge to the identity in the $C^2$-topology. More precisely,
we have
\begin{equation}
\Vert \overline{f}_k - {\rm id} \Vert_{2, I} < \frac{\varepsilon}{\sqrt{2^k}} \, \label{generalsecondorderestimate}
\end{equation}
for every diffeomorphism $\overline{f}_k \in \widetilde{S}(k)$ and for a suitable fixed $n$.
To show that convergence takes place in the $C^3$-topology
as well, we first estimate the third derivative $D^3 [f_1, f_2]$ of a commutator
$[f_1, f_2] = f_1 \circ f_2 \circ f_1^{-1} \circ f_2^{-1}$.
For this we shall use the fact that $f_1, f_2$ and their inverses $f_1^{-1}, f_2^{-1}$ are $C^2$-close to the identity.
Recall then that higher order derivatives of a composed function are given by Fa\`a di Bruno formula which, in the present
case, simply means
\begin{equation}
D^3 (f_1 \circ f_2) = D^3_{f_2(x)} f_1 . (D_xf_2)^3 + 3 D^2_{f_2(x)} f_1 . D^2_x f_2 . D^1_x f_2 + D^1_{f_2(x)} . D^3_x f_2 \, .
\label{faadibrunoorder3}
\end{equation}
Thus, if $\varepsilon >0$ is sufficiently small, we have
\begin{eqnarray*}
\vert D^1 (f_1 \circ f_2) -1 \vert & \leq & 3 \max \{ \sup_I \vert D^1 (f_1 - {\rm id} )\vert , \,
\sup_I \vert D^1 (f_2 - {\rm id} )\vert \}; \\
D^2 (f_1 \circ f_2) & \leq & 3 \max \{ \sup_I \vert D^2 f_1 \vert , \, \sup_I \vert D^2 f_2 \vert \}; \\
D^3 (f_1 \circ f_2) & \leq & 3 \max \{ \sup_I \vert D^3 f_1 \vert , \, \sup_I \vert D^3 f_2 \vert \} .
\end{eqnarray*}
Similar estimates also hold for $D^1 (f_1^{-1} \circ f_2^{-1})$, $D^2 (f_1^{-1} \circ f_2^{-1})$, and $D^3 (f_1^{-1} \circ f_2^{-1})$.
Now applying again the preceding estimates to $(f_1 \circ f_2) \circ (f_1^{-1} \circ f_2^{-1})$, we conclude that
\begin{equation}
D^3 [f_1, f_2] \leq 10 \max \{ \sup_I \vert D^3 f_1 \vert , \, \sup_I \vert D^3 f_2 \vert, \, \sup_I \vert D^3 f_1^{-1} \vert,
\, \sup_I \vert D^3 f_2^{-1}  \vert \} \, \label{fundamentalestimatethirdderivative}
\end{equation}
provided that $f_1$, $f_2$, $f_1^{-1}$, and $f_2^{-1}$ are $\varepsilon$-close to the identity (for some small
$\varepsilon >0$ fixed). From Estimate~(\ref{fundamentalestimatethirdderivative}), we conclude that
\begin{eqnarray*}
D^3 (F^{-n} \circ [f_1, f_2] \circ F^n) & = & D^3 (\lambda^{-n} . [f_1, f_2] (\lambda^n x)) \\
& \leq & 10 \lambda^{2n} \max \{ \sup_I \vert D^3 f_1 \vert , \, \sup_I \vert D^3 f_2 \vert,
\, \sup_I \vert D^3 f_1^{-1} \vert, \, \sup_I \vert D^3 f_2^{-1}  \vert \} \, .
\end{eqnarray*}
If $n$ is chosen so that $\lambda^{2n} < 1/10$, there follows that the third order derivatives of elements in $\widetilde{S} (1)$
are smaller than the maximum of the third order derivatives of elements in $S(0)$. This procedure can be iterated to
higher order commutators by virtue of Estimate~(\ref{generalsecondorderestimate}) so that
third order derivatives of elements in $\widetilde{S} (k)$ actually
decay geometrically with $k$. The claim results at once.\qed

The remainder of the proof of Proposition~\ref{Appendix-ReductionProposition} is a straightforward induction step.
By repeating the previous discussion, we just need to prove that a locally $C^r$-non-discrete group is also
locally $C^{r+1}$-non-discrete provided that $r \geq 2$. The argument is totally analogous to the one employed in the proof
of the above claim (the general Fa\`a di Bruno formulas can be used in the context). The detail is left to the reader.\qed

\vspace{0.1cm}

\noindent $\bullet$ {\bf Final comments}. We close this paper by pointing out a couple of specific issues involved in
our regularization scheme for iterated commutators, as explained above and in Section~4. First, the reader will note that the
analytic assumption is not needed in order to ensure the corresponding diffeomorphisms converge to the identity. The
importance of the analytic assumption lies in the fact that the sequence of sets $S(k)$ (and hence $\widetilde{S} (k)$)
does not degenerate into $\{ {\rm id} \}$. As mentioned this result is due to Ghys~\cite{ghys} and has a formal algebraic
nature: it depends on ensuring that a $C^{\infty}$-diffeomorphism $f$ of $S^1$
coincides with the identity so long there is a point in $S^1$ at which $f$ is $C^{\infty}$-tangent to the identity.
It would be nice to know whether or not there are finitely
generated pseudo-solvable, yet non-solvable, groups in ${\rm Diff}^{\infty} (S^1)$.

Finally note also that our regularization technique falls short of working in the $C^1$-case. Therefore, even in
the analytic category, we have not proved that a locally $C^1$-non-discrete subgroup of $\dif$ is also locally
$C^{\infty}$-non-discrete. Although this statement is very likely to hold, the renormalization procedure $x \mapsto
\lambda x$ used here does not decrease the first order derivative of the diffeomorphism and this accounts for the
special nature of locally $C^1$-non-discrete groups. To overcome this difficulty, our iteration scheme must be further
elaborated. This can probably be done by suitably
adding further ``take the commutator'' steps so as to keep control on the growing rate of first order derivatives.

\begin{flushleft}

{\sc Anas Eskif \hspace{0.5cm} \& \hspace{0.5cm} Julio C. Rebelo}\\
Institut de Math\'ematiques de Toulouse\\
Universit\'e de Toulouse\\
118 Route de Narbonne F-31062, Toulouse\\
FRANCE\\
Anas.Eskif@math.univ-toulouse.fr\\
rebelo@math.univ-toulouse.fr

\end{flushleft}

\end{document}